\documentclass[preprint,authoryear,3p]{elsarticle}

\usepackage{graphicx}

\usepackage{amssymb,amsthm,amsmath}
\usepackage{cleveref}
\usepackage{float}
\usepackage{multirow}
\usepackage{longtable,lscape}
\usepackage{cases}
\usepackage{url}
\usepackage{enumitem}
\usepackage[usenames, dvipsnames]{color}

\usepackage{fancyhdr}
\theoremstyle {plain}
\newtheorem{thm}{Theorem}[section]

\newtheorem{prop}[thm]{Proposition}

\theoremstyle{definition}

\theoremstyle{remark}

\newtheorem*{note}{Note}

\usepackage{titlesec}
\titleformat{\paragraph}
{\normalfont\normalsize\bfseries}{\theparagraph}{1em}{}
\titlespacing*{\paragraph}
{0pt}{3.25ex plus 1ex minus .2ex}{1.5ex plus .2ex}

\usepackage{tikz}
\usetikzlibrary{arrows.meta}
\usepackage{subcaption}

\DeclareMathOperator*{\argmin}{argmin}

\numberwithin{equation}{section}

\journal{ISTTT24} 

\begin{document}
{\vspace*{-4cm} \small \color{red} Please cite this paper as: Huang, K., Chen, X., Di, X., \& Du, Q. (2021). Dynamic driving and routing games for autonomous vehicles on networks: A mean field game approach. Transportation Research Part C: Emerging Technologies, 128, 103189.}

\begin{frontmatter}



\title{Dynamic Driving and Routing Games for Autonomous Vehicles on Networks: 
A Mean Field Game Approach}



\author[cuAPAM]{Kuang Huang}
\author[cu]{Xu Chen}
\author[cu,dsi]{Xuan Di\corref{cor}}
\ead{sharon.di@columbia.edu}

\author[cuAPAM,dsi]{Qiang Du}

\cortext[cor]{Corresponding author. Tel.: +1 212 853 0435;}

\address[cuAPAM]{Department of Applied Physics and Applied Mathematics, Columbia University}
\address[cu]{Department of Civil Engineering and Engineering Mechanics, Columbia University}
\address[dsi]{Data Science Institute, Columbia University}

\begin{abstract}

This paper aims to answer the research question as to optimal design of decision-making processes for autonomous vehicles (AVs), including dynamical selection of driving velocity and route choices on a transportation network.
Dynamic traffic assignment (DTA) has been widely used to model travelers' route choice or/and departure-time choice and predict dynamic traffic flow evolution in the short term. 
However, the existing DTA models do not explicitly describe one's selection of driving velocity on a road link. 
Driving velocity choice may not be crucial for modeling the movement of human drivers
but it is a must-have control to maneuver AVs.
In this paper, we aim to develop a game-theoretic model to solve for AVs' optimal driving strategies of velocity control in the interior of a road link and route choice at a junction node. 
To this end, we will first reinterpret the DTA problem as an $N$-car differential game and show that this game can be tackled with a general mean field game-theoretic framework. 
The developed mean field game is challenging to solve because of the forward and backward structure for velocity control and the complementarity conditions for route choice. 
An efficient algorithm is developed to address these challenges.    
The model and the algorithm are illustrated on the Braess network and the OW network with a single destination. 
On the Braess network, we first compare the LWR based DTA model with the proposed game and find that the driving and routing control navigates AVs with overall lower costs. 
We then compare the total travel cost without and with the middle link and find that the Braess paradox may still arise under certain conditions. 
We also test our proposed model and solution algorithm on the OW network.
\end{abstract}


\end{frontmatter}

\section{Motivation}

When all the human-driven vehicles (HVs) on public roads are replaced by autonomous vehicles (AVs), 
how should we design \emph{decision-making processes for AVs such that they dynamically select driving velocity and route choices on a transportation network to minimize some travel costs over a predefined planning horizon?}
The challenge is, however, when a large number of AVs with different origin-destination pairs navigate a transportation network, traffic evolves as time elapses and congestion arises, which could increase AVs' travel costs and in turn influence their decision-making. 
We are interested in understanding the \emph{equilibrium state} in which each AV cannot minimize its travel cost by unilaterally switching driving velocity and route choices. 

The traffic assignment problem models one's route choice behavior while interacting with other travelers (which is the travel demand) in order to predict network-wide traffic congestion. 
While static traffic assignment is proposed for the long-term planning purpose \citep{di2013boundedlyTRB,di2014braess,di2016boundedly}, 
dynamic traffic assignment (DTA) is the most popular prescriptive and normative approach to 
model travelers' route choice or/and departure-time choice and predict dynamic traffic flow evolution in a short term \citep{friesz1993variational,peeta2001foundations,ban2012continuous}.  
DTA integrates both notions of travel demand (consistent with static traffic assignment) and traffic flow (i.e., dynamic traffic evolution)
and is thus appropriate for modeling dynamic movement of human drivers across a network \citep{friesz2013dynamic}. 
DTA problems bear many variants and interested readers can refer to \cite{peeta2001foundations} and \cite{nie2005comparative} for a comprehensive overview of DTA models. 

The existing DTA models do not explicitly depict one's selection of driving velocity on a road link. 
Travel speed is either determined by traffic density when macroscopic traffic flow models are used 
or is simplified when bottleneck or queuing models are used. 
For example, in the LWR model \citep{lighthill1952sound,richards1956shock}, one's vehicle velocity on a link is determined by traffic density on that link via fundamental diagrams \citep{friesz2013dynamic,han2016continuity,han2019computing}. 
On a link with bottlenecks, one vehicle moves in the interior of the link without incurring any travel time and only joins a queue at the end of the link \citep{ban2012continuous,osorio2011dynamic}. 
Driving behavioral simplification mentioned above may not be so crucial for the existing models that were primarily developed for human drivers. 
However, when it comes to design a network-wide control scheme for AVs to navigate a transportation network, 
optimal selection of driving velocity is a must-have control to maneuver AVs on a link, in addition to departure-time choice at origins and route choice at intermediate junctions.  

In this paper, we aim to develop a game-theoretic model to solve for AVs' optimal driving strategies in longitudinal control and route choice.
To this end, we will first reinterpret the DTA problem 
as an $N$-car differential game 
and show that this game can be tackled with a general mean field game-theoretic framework. 
Then we will generalize the existing DTA problem to a broader context when one can select driving velocity while moving on a network. 

In the remainder of the paper, 
we will first review related literature on DTA and mean field games in Section~\ref{sec:lit}. 
We will then define one AV's optimal control problem to navigate a transportation network, and then extend it to a differential game when multiple AVs solve their individual optimal control problems while interacting among one another on a congested network in Section~\ref{sec:mc_one}. 
Due to complexity of the differential game, we will apply the mean field approximation and derive a mean field game for a large number of AVs in Section~\ref{sec:MFG}. We will demonstrate the connection between the mean field game and the classical dynamic user equilibrium concept in Section~\ref{sec:connection}.
In Section~\ref{sec:algo}, an efficient algorithm is developed to solve the mean field equilibrium. 
In Section~\ref{sec:numerical}, we will do some numerical experiments to demonstrate the solved mean field equilibrium on Braess networks and the OW network. 
Conclusions follow in Section~\ref{sec:conclud} with future research directions. 

\section{Literature review} \label{sec:lit}

\subsection{Dynamic user equilibrium} 

Dynamic user equilibrium (DUE) or dynamic user optimal 
describes the equilibrium state in which one cannot minimize his or her generalized travel cost by unilaterally switching
route choices or departure times simultaneously \citep{friesz2019mathematical}.  
Solving DUE usually requires to solve two subproblems: a dynamic network loading (DNL) procedure and a route and/or departure-time choice model. 

\textbf{The DNL procedure} propagates aggregate traffic flow dynamics and congestion in space and time (depicted by link volume, exit flow, link/path travel time or delay, queue length) given link inflows or path departure rates. 
It is composed of a link model and a junction model.
A \emph{link model} captures three phenomena: flow conservation, flow behavior (to compute travel delay), and flow propagation (how link volume propagates to exit \citep{bliemer2017genetics, yu2020day}).   
These phenomena are modeled with two types of functions: the delay function model and the exit-flow function model \citep{nie2005comparative}. 
The former computes link delays or queues explicitly using a linear delay function \citep{friesz1993variational,xu1999advances}, 
while the latter implicitly derives link traversal time from link exit flow. 
All the existing exit-flow functions model link flow on a macroscopic scale, including:   
M-N model \citep{merchant1978model,merchant1978optimality}, 
whole link model \citep{xu1999advances,carey2002behaviour,carey2003comparing,nie2005delay}, 
traffic flow models without spillback \citep{friesz2013dynamic} or with spillback \citep{gentile2007spillback,han2016continuity}, 
discrete-space traffic flow models including cell transmission models (CTM) \citep{lo1999dynamic,lo2002cell,szeto2004cell,zhu2015dta} and link transmission models (LTM) \citep{ukkusuri2012dta,gentile2015using,batista2021dta}, 
point-queue or bottleneck models such as single queue \citep{ban2012continuous,han2013partial1} or double queues \citep{osorio2011dynamic},  
and physical or spatial queue models \citep{kuwahara2001dynamic,zhang2013modelling,doan2015queue}.  
Path delay can be computed using an implicit path delay operator \citep{friesz2013dynamic,han2016continuity}. 
A \emph{junction model} determines how flows are distributed at intersections. 
In traffic flow models on a network, boundary conditions at a junction can become complex and may lead to ill-posed models.  
Fixed turning ratios or entropy optimization is imposed for a unique flow distribution in LWR models \citep{han2016continuity}, 
and demand-supply functions are proposed with sending and receiving flows as functions of traffic states \citep{lebacque1999modelling}. 

\textbf{A choice model} stipulates how an individual traveler makes route or departure-time choice, given the prevailing or predictive traffic condition. 
If the prevailing or present traffic information is used for decision-making, we call it the instantaneous UE (IUE) principle \citep{wie1990dynamic,boyce1995solving,lam1995dynamic,ban2012continuous}, 
otherwise it is ideal DUE \citep{ran1996link1,ran1996link2,friesz2013dynamic}. 

The inherent integration of dynamic traffic evolution and route choices 
enables
the DUE problem to be analytically formulated as
differential variational inequalities (DVI) \citep{friesz2019mathematical}
or differential complementarity systems (DCS) \citep{ban2012continuous}. 
Depending on spatial and temporal resolutions,   
DUE is described by 
algebraic difference equations or  
differential-algebraic equations (DAEs), which is a system of (time-delayed) ordinary differential equations (ODEs) or partial differential equations (PDEs) with algebraic constraints. 
Due to computational intractability of the analytical DUE problem, 
a large body of literature has resorted to simulations or heuristics \citep{mahmassani2001dynamic,lo2002cell,chiu2013technical} to compute the traffic equilibrium state. 
\cite{levin2016multiclass} is one among a few that designs networked traffic controls for AVs using the DUE framework with CTMs. 
To model the mixed traffic comprised of AVs and human drivers, a multi-class LWR model is developed with a fundamental diagram assuming that AVs and human drivers have different reaction times but the same driving speed. 
The maximum capacity varies with the penetration rate of AVs. 
At junctions, intelligent traffic management 
policy is developed to compute nodal delay for each vehicle class. 
Simulations are performed to compute optimal controls for AVs and the resultant equilibrium.  

\subsection{Multi-agent differential games and mean field games}

The spatio-temporal evolution of aggregate traffic dynamics essentially arises from a large number of individuals' dynamic travel choices and their complex interactions on a network. 
The fundamental tool to model individuals' decision-making processes in a multi-agent dynamic system is 
the non-cooperative multi-agent differential game \citep{dockner2000differential}. 
In a differential game, each agent solves its own control from an optimal control problem when every other agent in the system does so. 
Accordingly, one's optimal control problem is coupled with all others through the system state.  

Assuming agents are anonymous, 
instead of solving a long list of highly coupled optimal control problems in a multi-agent differential game, 
mean field approximation can be applied to exploit the ``smoothing" effect of large numbers of interacting individuals. 
At equilibrium, each agent interacts and reacts only to a ``mass" resulting from the aggregate effect of all agents. 
The mean field game (MFG) is a micro-macro model that allows one to define individual drivers on a microscopic level as rational, utility-optimizing agents
while translating their rich microscopic behaviors to a macroscopic scale. 
It is composed of two coupled components: 
agent dynamic and mass dynamic. 
The \emph{agent dynamic} models individuals' dynamics with an optimal control in the form of a backward Hamilton-Jacobi-Bellman (HJB) equation, 
while the \emph{mass dynamic} models system evolution arising from individuals' choices in the form of a forward Fokker-Planck-Kolmogorov (FPK) equation.

MFGs 
have become an increasingly popular tool to design new decision-making processes in 
finance \citep{gueant2011mean,lachapelle2010computation}, 
engineering \citep{djehiche2016mean,couillet2012electrical}, 
learning \citep{gummadi2012mean,iyer2014mean}, 
social science \citep{degond2014large}, 
and pedestrian crowd movement \citep{lachapelle2011mean,burger2013mean}.
\cite{djehiche2016mean} built a connection between the Wardrop equilibrium and the mean field equilibrium
by reinterpreting path flows as the mean field of individuals' route choices.  
In traffic modeling, there exist only a few studies that employ MFGs on  
velocity control \citep{kachroo2016inverse,chevalier2015micro,kachroo2017multiscale,chen2016optimal} and
lane-change \citep{festa2017mean}. 
The authors' recent studies 
have employed this tool to solve for optimal velocity control of AVs in pure AV traffic \citep{huang2019game} or mixed AV-HV traffic \citep{huang2019mixed,huang2019stable}. 
However, these studies are primarily focused on a single road rather than on a network. 
On a transportation network, 
\cite{bauso2016density} used the notion of MFG for routing games where static traffic density flow is solved on each link but its temporal evolution is not captured. 

This paper aims to fill the above gaps by formulating the AVs' dynamic control problem on a network as a MFG. 
Due to the integration of driving velocity control and route choices, 
the formulation and computation of the MFG on a network is expected to be more challenging. 

\section{$N$-car driving and routing differential game on a network}\label{sec:mc_one} 

\subsection{Problem statement} 
		
There are $N$ cars indexed by $n\in\{1,2,\dots,N\}$ who need to move from their initial positions to a destination on a congested transportation network during a predefined time horizon. 
Their travel mode is autonomous vehicle (AV). 
We assume these AVs as intelligent agents who select optimal control strategies throughout the network, in order to minimize their pre-programmed travel cost functionals over the predefined planning horizon. 
In other words, each AV aims to select a minimum-cost driving profile to navigate the transportation network.    
En-route, one AV needs to select its driving speed continuously; at any junction, the AV needs to choose the next-go-to link. 
When one AV selects its own driving speed and routing policy while everybody else does so simultaneously, a non-cooperative \emph{differential game} forms.  
We are interested in an equilibrium control strategy at which no AV can improve its travel cost by unilaterally switching its driving speed or routing policy. 

\textbf{Assumptions.} To clarify, we first make the following modeling assumptions.
\begin{enumerate}
	\item Each AV receives the complete information of prevailing traffic, including all other AVs' speeds and route choices.
	\item Each AV anticipates future traffic on the predefined time horizon $[0,T]$ where $T>0$.
	\item Each AV controls its driving speed in the interior of a link and next-go-to link choice at a junction
	to minimize a predefined travel cost functional over a network.
	\item There is a vertical queue at each junction. AVs passing through a junction experience a queuing delay when the queue at that junction is non-empty.
	\item The ideal DUE principle is employed to solve the traffic equilibrium, which is, the actual traffic information is used for decision-making. In other words, AVs make optimal choices according to what would actually happen in the future.
\end{enumerate}	


\textbf{Notation.} To present the mathematical formulation of the differential game, we will first introduce the notation that will be used in this paper. 
A transportation network is represented by a directed graph $\mathcal{G}=\left\lbrace \mathcal{N},\mathcal{L}\right\rbrace$ that includes a set of nodes $\mathcal{N}$ and a set of links $\mathcal{L}\subset\mathcal{N}\times\mathcal{N}$.
For each link $l=(i,j)\in\mathcal{L}$ where $i,j\in\mathcal{N}$, denote its starting point as $\text{START}(l)=i$ and its end point as $\text{END}(l)=j$. 
Accordingly, one can write $l=\left(\text{START}(l),\text{END}(l)\right)$ for each link $l\in\mathcal{L}$. We will use symbols $l$ or $(i,j)$ interchangeably to represent links in $\mathcal{L}$.
The length of link $l\in\mathcal{L}$ is denoted as $\text{len}(l)\geq0$.
For each node $i\in\mathcal{N}$, denote $\text{IN}(i)\subset\mathcal{L}$ the set of links whose end point is node $i$ and $\text{OUT}(i)\subset\mathcal{L}$ the set of links whose starting point is node $i$. 
Denote $\mathcal{N}_O\subset\mathcal{N}$ the set of origins and $s\in\mathcal{N}$ the destination. 
The intermediate node set is ${\cal N}_I={\cal N}\setminus{\cal N}_O\setminus \{s\}$.

We will first formulate one AV's hybrid optimal control problem on a network in Section~\ref{sec:optm_control}. Then the optimal control problem will be extended to the $N$-car differential game in Section~\ref{sec:diff_game}. Throughout this paper, we will use AVs or cars interchangeably to refer to the intelligent agents who navigate a transportation network with optimal driving strategies.

\subsection{One AV's hybrid optimal control on a network}\label{sec:optm_control}

To formulate one AV's optimal control problem on a network, we first introduce the state and control variables in the interior of a link and at a node.

\textbf{State variable.} Consider the $n^{\text{th}}$ car ($n=1,\dots,N$) driving on the network, the car's state is described in two scenarios:

\begin{itemize}
    \item[(i)] The car is driving in the interior of a link at time $t\in[0,T]$. We denote the link the car is on as $l^{(n)}(t)$ and parameterize the car's position on this link as a real number $x^{(n)}(t)\in\left(0,\text{len}(l^{(n)}(t))\right)$. In this case, the car's state is specified by:
    \begin{align}
    	X^{(n)}(t)=\left(x^{(n)}(t),l^{(n)}(t)\right)\in\mathbb{R}\times\mathcal{L},\quad 0< x^{(n)}(t)<\text{len}(l^{(n)}(t)).
    \end{align}
    \item[(ii)] The car is in the queue at a node at time $t\in[0,T]$. We denote the node as $i^{(n)}(t)$. 
    In this case, the car's state is specified by:
    \begin{align}
        X^{(n)}(t)=i^{(n)}(t)\in\mathcal{N}.
    \end{align}
\end{itemize}



\textbf{Control variable.} Correspondingly, the car has two types of decisions to make during its travel on the network: speed control and route choice. 
\begin{itemize}
    \item[(i)] The $n^{\text{th}}$ car's driving speed is denoted as $u^{(n)}(t)$ when the car is in the interior of any link at time $t\in[0,T]$.
    Its feasible region is $u^{(n)}(t)\in[u_{\text{min}},u_{\text{max}}]$ where $u_{\text{min}}$ and $u_{\text{max}}$ are the minimum and maximum speeds for all cars.
    \item[(ii)] The $n^{\text{th}}$ car's next-go-to link choice is denoted as $\gamma^{(n)}(i,t)$ when the car exits the queue at node $i\in\mathcal{N}$ at time $t\in[0,T]$.
    Its feasible region is the outgoing link set $\text{OUT}(i)$ of node $i$.
\end{itemize}
Summarizing these two types of decisions, we have a driving strategy vector:
\begin{align}
    \alpha^{(n)}=[u^{(n)}(\cdot),\gamma^{(n)}(\cdot,\cdot)],\quad n=1,\dots,N.
\end{align}

Because the speed control is implemented as continuous-in-time driving speed profile while the route choice is implemented as discrete-in-time next-go-to link choices, 
the optimal control problem to solve these two types of decisions is a \emph{hybrid optimal control problem}. 
With state and control variables defined, we formulate the AV's state dynamics under a control in Section~\ref{sec:state_dynamics} and define the AV's travel cost functional in Section~\ref{sec:cost_functional}.

\subsubsection{One AV's state dynamics}\label{sec:state_dynamics}

Knowing the control variables for the $n^{\text{th}}$ car, we are ready to formulate a dynamical system that stipulates the evolution of the car's state under a hybrid control.
There are two types of car state evolution: link dynamic and node transition.

\textbf{Link dynamic.} When the $n^{\text{th}}$ car drives in the interior of a link $l^{(n)}(t)\in {\cal L}$ at time $t$, its position $x^{(n)}(t)$ 
changes continuously according to its controlled speed $u^{(n)}(t)$. The car's movement on a link is accordingly modeled by an ODE:  $dx^{(n)}(t)/dt=u^{(n)}(t)$. 
Since the car keeps driving on the same link, the link index of the car's state remains the same and the dynamic of the link index becomes: $dl^{(n)}(t)/dt=0$. 

\textbf{Node transition.} 
When the car reaches the end point of the link it is driving on, which is $x^{(n)}(t-)=\text{len}(l^{(n)}(t-))$, 
a transition of state happens at time $t$. 
Here 
$t-$ represents the left time limit. 
With the transition, the car moves to the node $i^{(n)}(t)=\text{END}(l^{(n)}(t-))$ and enters the queue at node $i^{(n)}(t)$. We denote the car's exit time from the queue as $e^{(n)}(t)\in [t,T]$. During the time interval $t\leq \tau<e^{(n)}(t)$, the car's state keeps unchanged, i.e.,  $di^{(n)}(\tau)/d\tau=0$. 
Once the car exits the queue, it leaves the node at time $e^{(n)}(t)$.
According to the car's next-go-to link choice, it switches to the link $\gamma^{(n)}\left(i^{(n)}(t), e^{(n)}(t)\right)$
and its initial position on the new link is reset to zero at time $e^{(n)}(t)$. We denote the set of all time instants of node transition during the $n^{\text{th}}$ car's trip over the network as $R^{(n)}\subset[0,T]$.

\textbf{Initital State.} Suppose the $n^{\text{th}}$ car enters the network through the origin $o^{(n)}\in\mathcal{N}_O$ at time $t_0^{(n)}\geq0$, the car's initial state is given by $i^{(n)}(t_0^{(n)})=o^{(n)}$. The car's initial state corresponds to its first node transition at time $t_0^{(n)}$. That is, the car switches to the link
$\gamma^{(n)}\left(o^{(n)},
e^{(n)}(t_0^{(n)})\right)$ when it exits the queue at the origin $o^{(n)}$ at time $e^{(n)}(t_0^{(n)})$.


Summarizing both link dynamic and node transition for the $n^{\text{th}}$ car ($n=1,\dots,N$), we now obtain its dynamical system on a network as:

\begin{subequations}
\begin{align}
\mbox{(Link dynamic)}\quad&\frac{dx^{(n)}(t)}{dt}=u^{(n)}(t),\quad \frac{dl^{(n)}(t)}{dt}=0, \quad \text{if } 0<x^{(n)}(t)<\text{len}(l^{(n)}(t));\\
\mbox{(Node transition)}\quad& i^{(n)}(t)=\text{END}(l^{(n)}(t-)),\quad t\in R^{(n)},\quad \text{if } x^{(n)}(t-)=\text{len}(l^{(n)}(t-));\\
&\frac{di^{(n)}(\tau)}{d\tau}=0 \text{ for } t<\tau<e^{(n)}(t),\notag\\ &x^{(n)}(e^{(n)}(t))=0,\quad l^{(n)}(e^{(n)}(t))=\gamma^{(n)}\left(i^{(n)}(t),e^{(n)}(t)\right), \quad \text{if } t\in R^{(n)};\\
\mbox{(Initial state)}\quad &i^{(n)}(t_0^{(n)})=o^{(n)},\quad t_0^{(n)}\in R^{(n)}.
\end{align}\label{eq:one_car_dynamics}
\end{subequations}


Eq.~\eqref{eq:one_car_dynamics} is a hybrid dynamical system. The link dynamic is modeled by continuous ODEs, 
while the node transition
happens only at discrete time instants in $R^{(n)}$. We denote $R^{(n)}=\left\{t_0^{(n)},t_1^{(n)},\cdots,t_{|R^{(n)}|-1}^{(n)}\right\}$ where $|R^{(n)}|$ is the number of times of the $n^{\text{th}}$ car's node transition and $t_0^{(n)}<t_1^{(n)}<\cdots<t_{|R^{(n)}|-1}^{(n)}$ are the transition time instants.
By integrating the $n^{\text{th}}$ car's dynamical system, we can obtain: (i) the car's trajectory $x^{(n)}(t)$ in the interior of a link when $e^{(n)}(t_k^{(n)})\leq t<t_{k+1}^{(n)}$ for $k=0,\cdots,|R^{(n)}|-1$; (ii) the car's queuing delay $e^{(n)}(t_k^{(n)})-t_k^{(n)}$ at a node for $k=0,\cdots,|R^{(n)}|-1$.
If the $n^{\text{th}}$ car arrives at the destination before time $T$, we denote this arrival time as $t_f^{(n)}<T$; otherwise the car is still on the network at the end of the planning horizon $[0,T]$ and we denote $t_f^{(n)}=T$. In both scenarios, the $n^{\text{th}}$ car's actual travel time horizon on the network is $[t_0^{(n)}, t_f^{(n)}]$. Note that in the latter scenario, the car is either on a link or in the queue at a node at time $t_f^{(n)}=T$. If the car is in the queue at a node and its exit time $e^{(n)}(t_{|R^{(n)}|-1}^{(n)})>T$, we will make the convention that $e^{(n)}(t_{|R^{(n)}|-1}^{(n)})=T$. We will also make the convention that $t_{|R^{(n)}|}=t_f^{(n)}$.

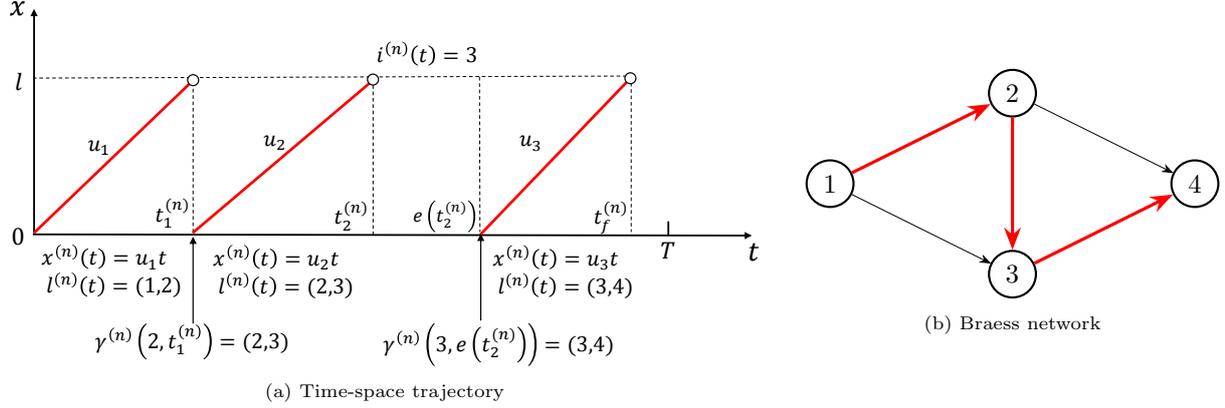
\begin{figure}[htbp]
	\centering
	\begin{subfigure}{.6\textwidth}
	    \centering
    	\includegraphics[width=\textwidth]{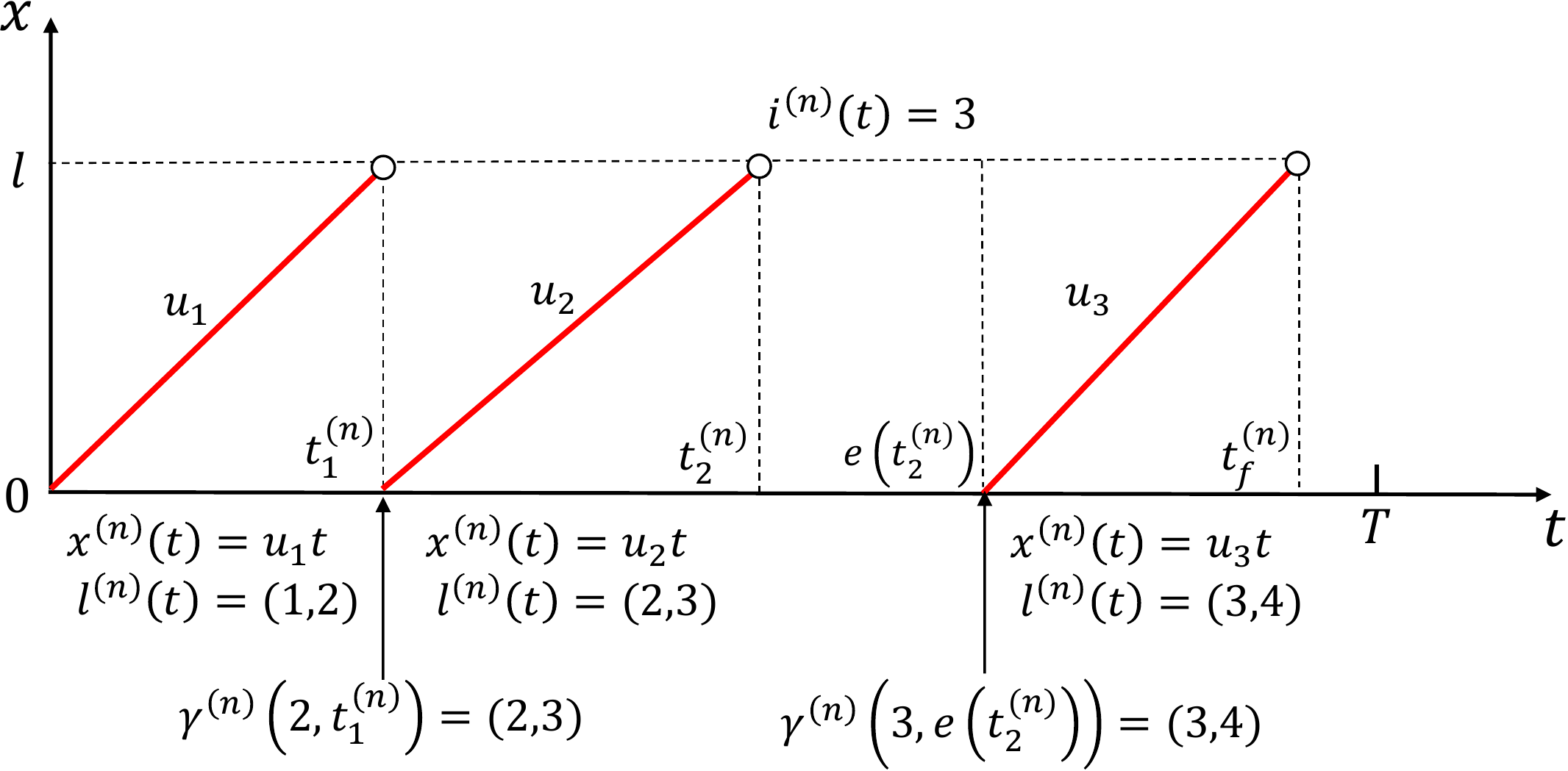}
    	\caption{Time-space trajectory}
    	\label{fig:one_car_dyn_traj}
	\end{subfigure}
	\begin{subfigure}{.39\textwidth}
        \centering
        \begin{tikzpicture}[scale=.4]
            \begin{scope}[every node/.style={circle,thick,draw}]
                \node (s) at (-6,0) {1};
                \node (d) at (6,0) {4};
                \node (1) at (0,3) {2};
                \node (2) at (0,-3) {3};
            \end{scope}
            \begin{scope}[>=Stealth,auto,
                          every edge/.style={draw=red,very thick}]
                \path [->] (s) edge node {} (1);
                \path [->] (2) edge node {} (d);
                \path [->] (1) edge node {} (2);
            \end{scope}
            \begin{scope}[>=Stealth,auto,
                          every edge/.style={draw=black,very thin}]
                \path [->] (s) edge node {} (2);
                \path [->] (1) edge node {} (d);
            \end{scope}
        \end{tikzpicture}
        \caption{Braess network}
        \label{fig:one_car_dyn_net}
    \end{subfigure}
	\caption{One AV's state dynamics across the Braess network}\label{fig:one_car_dyn}
\end{figure}

To demonstrate one trajectory realization solved from the above dynamical system, we plot the $n^{\text{th}}$ car's time-space trajectory (Figure~\ref{fig:one_car_dyn_traj}) on the Braess network (Figure~\ref{fig:one_car_dyn_net}). 

Suppose the car enters the network through the origin node $1$ at time $t_0^{(n)}=0$ and travels along the path $1\to2\to3\to4$. It is assumed that all links $(1,2)$,  $(2,3)$, $(3,4)$ have the same length $l$. 
In the interior of a link, the position of the car, denoted as $x$, changes continuously with respect to time $t$ under the influence of its speed control. 
Suppose the car keeps a constant speed on each link, its time-space trajectory consists of slanted straight lines.
The time instants when the car arrives at a non-destination node are $R^{(n)}=\left\{t_0^{(n)}, t_1^{(n)}, t_2^{(n)}\right\}$.
Suppose that the queues at node $1$ and $2$ are always empty but the queue at node $3$ is non-empty at time $t_2^{(n)}$, the car exits those queues at time $e(t_0^{(n)})=t_0^{(n)}=0$, $e(t_1^{(n)})=t_1^{(n)}$ and $e(t_2^{(n)})>t_2^{(n)}$. Exiting a queue, the car switches to the next-go-to link and its position jumps back to zero on the new link. The car reaches the destination $s=4$ at time $t_f^{(n)}<T$ and ends its trip there.

%
%
%

\subsubsection{One AV's travel cost functional}\label{sec:cost_functional}

To select the optimal driving speed and routing policy, the $n^{\text{th}}$ car solves a hybrid optimal control problem by minimizing a predefined travel cost functional.
The travel cost functional consists of three parts: link running cost, node queuing cost and terminal cost.

\textbf{Link running cost.} Within each time interval $e^{(n)}(t_k^{(n)})\leq t<t_{k+1}^{(n)}$ for $k=0,\cdots,|R^{(n)}|-1$, the car drives in the interior of a link. We define the car's running cost on the link as 
\begin{align*}
    \int_{e^{(n)}(t_k^{(n)})}^{t_{k+1}^{(n)}}f^{(n)}_{\text{run}} \left(u^{(n)}(t),X^{(n)}(t)\right)\,dt,
\end{align*}
which is an integral of the car's instantaneous running cost $f^{(n)}_{\text{run}} \left(u^{(n)}(t),X^{(n)}(t)\right)$ over the time interval. Here $f^{(n)}_{\text{run}}(\cdot,\cdot)$ is the $n^{\text{th}}$ car's \emph{running cost function} that quantifies driving objectives such as efficiency and safety.

\textbf{Node queuing cost.} Within each time interval $t_k^{(n)}\leq t<e^{(n)}(t_k^{(n)})$ for $k=0,\cdots,|R^{(n)}|-1$, the car stays in the queue at a node. We define the car's queuing cost at the node as $f^{(n)}_{\text{que}}\left(e^{(n)}(t_k^{(n)})-t_k^{(n)}\right)$, in which $e^{(n)}(t_k^{(n)})-t_k^{(n)}$ is the car's queuing delay at the node. Here $f^{(n)}_{\text{que}}(\cdot)$ is the $n^{\text{th}}$ car's \emph{queuing cost function} that quantifies the impact of queuing delay.

\textbf{Terminal cost.} The car terminates its trip at the state $X^{(n)}(t_f^{(n)})$ at time $t_f^{(n)}$. We define the car's terminal cost as  $V^{(n)}_{\text{ter}}\left(t_f^{(n)}, X^{(n)}(t_f^{(n)})\right)$. Here $V^{(n)}_{\text{ter}}(\cdot,\cdot)$ is the $n^{\text{th}}$ car's \emph{terminal cost function} that represents the car's preference on its terminal state. 
There are three types of terminal state:
\begin{itemize}
    \item[(i)] The car arrives at the destination $s$ at time $t_f^{(n)}<T$. In this case the terminal cost is written as $V^{(n)}_{\text{ter}}\left(t_f^{(n)}, s\right)$ and we will always assume $V^{(n)}_{\text{ter}}\left(t_f^{(n)}, s\right)=0$.
    \item[(ii)] The car is in the interior of a link at time $t_f^{(n)}=T$. In this case the terminal cost is written as $V^{(n)}_{\text{ter}}\left(T, x^{(n)}(T), l^{(n)}(T)\right)$, where $l^{(n)}(T)$ is the link the car is on and $x^{(n)}(T)$ is the car's position on the link at time $T$.
    \item[(iii)] The car is in the queue at a node at time $t_f^{(n)}=T$. In this case the terminal cost is written as $V^{(n)}_{\text{ter}}\left(T, i^{(n)}(T)\right)$, where $i^{(n)}(T)$ is the node at which the car is in the queue at time $T$.
\end{itemize}


Summing up the three types of costs along the $n^{\text{th}}$ car's trip over the network, the general form of the $n^{\text{th}}$ car's travel cost functional is defined as:

\small{
\begin{align}
J_n(\alpha^{(n)}) 
=&\sum_{k=0}^{|R^{(n)}|-1}\left[\underbrace{\int_{e^{(n)}(t_k^{(n)})}^{t_{k+1}^{(n)}}f^{(n)}_{\text{run}} \left(u^{(n)}(t),X^{(n)}(t)\right)\,dt}_{\text{link running cost}}+\underbrace{f^{(n)}_{\text{que}}\left(e^{(n)}(t_k^{(n)})-t_k^{(n)}\right)}_{\text{node queuing cost}}\right]
+\underbrace{V^{(n)}_{\text{ter}}\left(t_f^{(n)},X^{(n)}(t_f^{(n)})\right)}_{\text{terminal cost}}.\label{eq:Jn}
\end{align}
}


Summarizing both state dynamics defined in Eq.~\eqref{eq:one_car_dynamics} and the travel cost functional defined in Eq.~\eqref{eq:Jn}, the $n^{\text{th}}$ car's hybrid optimal control problem on a network is defined as:

\small{
\begin{subequations}
\begin{align}
&\alpha^{*(n)}=\argmin_{\alpha^{(n)}} J_n(\alpha^{(n)}),\label{eq:one_car_optm_ctrl_cost} \\
\text{s.t. }&\begin{cases}
\frac{dx^{(n)}(t)}{dt}=u^{(n)}(t),\ \frac{dl^{(n)}(t)}{dt}=0, \quad \text{if } 0< x^{(n)}(t)<\text{len}(l^{(n)}(t));\\
i^{(n)}(t)=\text{END}(l^{(n)}(t-)),\ t\in R^{(n)}, \quad \text{if } x^{(n)}(t-)=\text{len}(l^{(n)}(t-));\\
\frac{di^{(n)}(\tau)}{d\tau}=0 \text{ for } t<\tau<e^{(n)}(t),\
x^{(n)}(e^{(n)}(t))=0,\ l^{(n)}(e^{(n)}(t))=\gamma^{(n)}\left(i^{(n)}(t),e^{(n)}(t)\right),\ \text{if } t\in R^{(n)};\\
i^{(n)}(t_0^{(n)})=o^{(n)},\ t_0^{(n)}\in R^{(n)}.
\end{cases}\label{eq:one_car_optm_ctrl_dynamic}
\end{align}\label{eq:one_car_optm_ctrl}
\end{subequations}
}


\subsection{N-car differential game on a network}\label{sec:diff_game}

The $n^{\text{th}}$ car ($n=1,\dots,N$) solves the optimal control problem defined by Eq.~\eqref{eq:one_car_optm_ctrl}. 
Remember that there are $N$ cars moving on the same transportation network.  
For any $n=1,2,\dots,N$, suppose the $n^{\text{th}}$ car knows others' driving strategies:
\begin{align}
	\alpha^{(-n)}=[\alpha^{(1)},\cdots,\alpha^{(n-1)},\alpha^{(n+1)},\cdots,\alpha^{(N)}]. 
\end{align}

By integrating Eq.~\eqref{eq:one_car_optm_ctrl_dynamic} for all other cars under their driving strategies, the $n^{\text{th}}$ car 
is able to predict the system state:
\begin{align}
	X^{(-n)}(t)=[X^{(1)}(t),\cdots,X^{(n-1)}(t),X^{(n+1)}(t),\cdots,X^{(N)}(t)],\label{eq:state_vec}
\end{align}
on the planning horizon $t\in[0,T]$. We make the convention that Eq.~\eqref{eq:state_vec} only includes the cars on the network at time $t$.


Due to interactions between cars, the $n^{\text{th}}$ car's 
travel cost functional $J_n(\cdot)$ defined in Eq.~\eqref{eq:Jn} now depends on others' states and actions. 
Accordingly, the $n^{\text{th}}$ car's travel cost functional is 
reformulated as follows:
\small{
\begin{align}
J_n^N(\alpha^{(n)},\alpha^{(-n)}) 
=&\sum_{k=0}^{|R^{(n)}|-1}\left[\underbrace{\int_{e^{(n)}(t_k^{(n)})}^{t_{k+1}^{(n)}}f^{(n)}_{\text{run}} \left(u^{(n)}(t),X^{(n)}(t),X^{(-n)}(t)\right)\,dt}_{\text{link running cost}}+\underbrace{f^{(n)}_{\text{que}}\left(e^{(n)}\left(t_k^{(n)},X^{(-n)}(t_k^{(n)})\right)-t_k^{(n)}\right)}_{\text{node queuing cost}}\right]\notag\\
&+\underbrace{V^{(n)}_{\text{ter}}\left(t_f^{(n)},X^{(n)}(t_f^{(n)})\right)}_{\text{terminal cost}},\label{eq:optmctrlcost_game}
\end{align}
}
where the $n^{\text{th}}$ car's running cost function $f_{\text{run}}^{(n)}(\cdot)$ and exit time function $e^{(n)}(\cdot)$ also depend on other cars' states.


A Nash equilibrium of the game is a tuple of controls $\alpha^*=[\alpha^{*(1)},\alpha^{*(2)},\dots,\alpha^{*(N)}]$ satisfying:
\begin{align} 
	J_n^N(\alpha^{*(n)},\alpha^{*(-n)}) \leq J_n^N(\alpha^{(n)},\alpha^{*(-n)}),\quad \forall \alpha^{(n)},\quad n=1,\dots,N.\label{eq:DGE}
\end{align}
At equilibrium, no car can improve its travel cost by unilaterally switching its driving speed or routing policy. 

When $N$ cars interact with one another while moving through a transportation network, we need to solve a total of $N$ optimal control problems that are coupled through states of all cars. This brings challenges in solving the equilibrium of the $N$-car differential game as $N$ becomes large \citep{cardaliaguet2010notes}.
To simplify this system, 
this paper will develop a scalable framework to solve an approximate equilibrium by resorting to the mean field approximation. 

\section{Mean field driving and routing game on a network} \label{sec:MFG}

In this section we will derive a mean field driving and routing game from the $N$-car differential game on a network. This derivation follows from the approach in \cite{huang2019game} and the derived mean field game can be viewed as the continuum limit of the $N$-car differential game as the number of cars $N\to\infty$.

We will use the \emph{mean field approximation} to bridge between the discrete differential game and continuum mean field game. To apply the mean field approximation, we need the following assumptions:\\
(A1)\quad All AVs are homogeneous in the sense that they have the same form of travel cost functional.\\
(A2)\quad All AVs are indistinguishable on the road.\\
(A3)\quad Each AV's instantaneous link running cost only depends on its driving speed and nearby traffic density; Each AV's instantaneous node queuing cost only depends on its queuing delay.\\
(A4)\quad Each AV stays in a queue until all the AVs entering the queue before him leave the queue, i.e., the first-in-first-out (FIFO) principle.\\
(A5)\quad The queue at each node has a bottleneck capacity. The outgoing flow rate from the queue reaches this capacity as long as the queue is non-empty.\\
(A6)\quad AVs implement their next-go-to link choices in a mixed strategy.

The mean field approximation is applied on both the AVs' speed control and routing policy. The derived mean field game is formulated as the coupled PDE system of a Hamilton-Jacobi-Bellman (HJB) equation and a continuity equation on a network.

\begin{figure}[htbp]
	\centering
	\includegraphics[width=.8\textwidth]{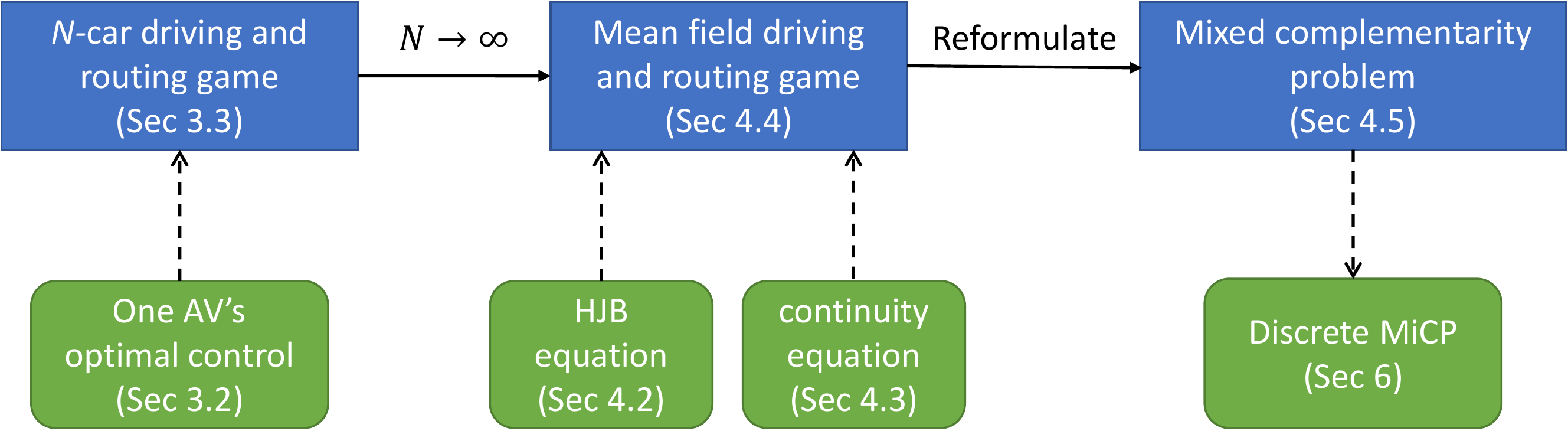}
	\caption{Connection between different problem formulations}\label{fig:diagram}
\end{figure}

To present how a MFG on a network is formulated, 
we will first apply the mean field approximation to the $N$-car differential game 
and derive a generic car's optimal control problem in Section~\ref{sec:mf_approx}. The HJB equation will be derived from the generic car's optimal control problem in Section~\ref{sec:hjb} and the continuity equation will be derived from the aggregate traffic evolution of all cars' dynamics in Section~\ref{sec:continuity}. The derived MFG system is presented in Section~\ref{sec:mfg_system} and reformulated as a mixed complementarity problem (MiCP) in Section~\ref{sec:reform}. The problem can be solved numerically on discrete space-time grids, which will be discussed in Section~\ref{sec:algo}.
Figure~\ref{fig:diagram} shows the whole diagram and illustrates connection between different problem formulations.

\subsection{Mean field approximation}\label{sec:mf_approx}

In this section, we will apply the mean field approximation to the $N$-car differential game introduced in Section~\ref{sec:diff_game} and derive a generic car's optimal control problem. The two steps are:
\begin{enumerate}
    \item Instead of tracking a total of $N$ states for every car, we define an aggregate traffic state on the network using the mean field approximation.
    \item Instead of solving a total of $N$ coupled optimal control problems for every car, we derive a generic car's optimal control problem by reformulating the travel cost functional defined in Eq.~\eqref{eq:optmctrlcost_game} and the equilibrium condition defined in Eq.~\eqref{eq:DGE} under the aggregate traffic state.
\end{enumerate}

\textbf{Mean field approximation.}
By the assumption (A2) that all AVs are indistinguishable, each AV only sees the aggregate distribution of all other AVs but cannot distinguish their identity numbers. As a corollary, each AV's travel cost functional only depends on the distribution of all other AVs on the network. When $N$ is very large, this distribution is well approximated by:
\begin{itemize}
    \item[(i)] The link traffic density $\rho_l(x,t)$, where $l\in\mathcal{L}$ is the link index, $x\in\left(0,\text{len}(l)\right)$ and $t\in[0,T]$ are the space and time coordinates; 
    \item[(ii)] The node queue size $Q_i(t)$, where $i\in\mathcal{N}\backslash\{s\}$ is the node index, $t\in[0,T]$ is the time coordinate.
\end{itemize}
We will define the aggregate traffic state on the network as the whole of all link traffic densities and node queue sizes. 
In other words, the traffic state represents a ``mean field'' of all AVs on the network.
 The complex interactions between individual AVs are then decoupled into two components: a generic car's optimal control under the traffic state, and the evolution of the traffic state from all cars' movements. 

\textbf{A generic car's optimal control.} 

In the $N$-car differential game, each AV aims to select its optimal control to minimize its travel cost functional defined in Eq.~\eqref{eq:optmctrlcost_game} in the existence of all other AVs. 
By the assumption (A1) that all AVs have the same form of travel cost functional, we can consider a generic car's optimal control problem under the traffic state.
Let us drop the car index $n$ for all functions and variables in Eq.~\eqref{eq:optmctrlcost_game}.
Then we reformulate the travel cost functional defined in Eq.~\eqref{eq:optmctrlcost_game} under the traffic state as follows:
\begin{enumerate}
    \item By the assumption (A3), we can replace the arguments $X^{(n)}(t),X^{(-n)}(t)$ in the generic car's running cost function $f_{\text{run}}(\cdot)$ by the
    local traffic density $\rho_{l(t)}(x(t),t)$. 
    \item By the assumption (A4), the generic car entering the queue at node $i(t)$ at time $t$ exits the queue at the time when all cars existed in the queue at time $t$ exit. By the assumption (A5), the total time for those cars to exit is $Q_{i(t)}(t)/M_{i(t)}$, where $M_{i(t)}>0$ denotes the bottleneck capacity of the queue at node $i(t)$, and the generic car's exit time is $e(t)=t+Q_{i(t)}(t)/M_{i(t)}$. When $t+Q_{i(t)}(t)/M_{i(t)}>T$, we make the convention that $e(t)=T$. Hence the generic car's exit time function is given by:
    \begin{align}
        e(t)=\min\left\{t+\frac{Q_{i(t)}(t)}{M_{i(t)}},T\right\}.\label{eq:delay_time}
    \end{align}
\end{enumerate}

Denote $\alpha=[u(\cdot),\gamma(\cdot,\cdot)]$ as the generic car's hybrid control. It is composed of the speed control $u(t)$ for $t\in[t_0,t_f]$ and the routing policy
$\gamma(i,t)$ for $i\in\mathcal{N}$ and $t\in[t_0,t_f]$. 
The generic car's optimal control problem under the traffic state is formulated as:
\begin{subequations}
\begin{align}
	\alpha^*=&\argmin_{\alpha} J(\alpha),\notag\\
	=& \argmin_{\alpha} \sum_{k=0}^{|R|-1}\left[\underbrace{\int_{e(t_k)}^{t_{k+1}}f_{\text{run}} \left(u(t),\rho_{l(t)}(x(t),t)\right)\,dt}_{\text{link running cost}}+\underbrace{f_{\text{que}}\left(e(t_k)-t_k\right)}_{\text{node queuing cost}}\right]
+\underbrace{V_{\text{ter}}\left(t_f,X(t_f)\right)}_{\text{terminal cost}}, \label{eq:generic_car_optm1}\\
	\text{s.t. }&\begin{cases}
    \frac{dx(t)}{dt}=u(t),\ \frac{dl(t)}{dt}=0, \quad \text{if } 0<x(t)<\text{len}(l(t));\\
    i(t)=\text{END}(l(t-)),\ t\in R, \quad \text{if } x(t-)=\text{len}(l(t-));\\
    \frac{di(\tau)}{d\tau}=0 \text{ for } t<\tau<e(t),\
    x(e(t))=0,\ l(e(t))=\gamma\left(i(t),e(t)\right),\ \text{if } t\in R;\\
    i(t_0)=o,\ t_0\in R.
    \end{cases} \label{eq:generic_car_optm2}
\end{align}\label{eq:generic_car_optm}
\end{subequations}
where:\\
$[t_0,t_f]$ is the generic car's travel time horizon on the network;\\
$R=\{t_0,t_1,\cdots,t_{|R|-1}\}$ is the set of all time instants of node transition of the generic car, and we make the convention that $t_{|R|}=t_f$;\\
$e(t)$ is the generic car's exit time function under the traffic state, which is defined by Eq.~\eqref{eq:delay_time};\\
$f_{\text{run}}(\cdot,\cdot)$, $f_{\text{que}}(\cdot)$ and $V_{\text{ter}}(\cdot,\cdot)$ are the generic car's running, queuing and terminal cost functions, respectively.


Given the traffic state on the network in the planning horizon $[0,T]$, we can solve a generic car's optimal control problem defined in Eq.~\eqref{eq:generic_car_optm} with any initial condition. 
When all cars follow their individual optimal controls, the traffic state evolves as the aggregate behavior of all cars' dynamics. Assuming the traffic state is known, one HJB equation is derived from all cars' optimal controls. 
Assuming all cars' optimal controls are known, one continuity equation is derived from their dynamic motions.
The mean field game (MFG) is then formulated as the coupled system of these two equations.


\subsection{Backward Hamilton-Jacobi-Bellman (HJB) equation}\label{sec:hjb}


With a generic car's optimal control problem  introduced, now we move on to derive the HJB equation. 
Suppose that: (i) the link traffic density $\rho_l(x,t)$ is known for all $l\in\mathcal{L}$, $x\in(0,\text{len}(l))$ and $t\in[0,T]$; (ii) the queue size $Q_i(t)$ is known for all $i\in\mathcal{N}$ and $t\in[0,T]$. The HJB equation determines a set of optimal velocity fields and next-go-to link distributions characterizing a generic car's optimal control. 

Our derivation is based on the dynamic programming principle. That is, we consider the following subproblem of the optimal control problem defined in Eq.~\eqref{eq:generic_car_optm}: suppose the generic car starts from any state $X(t_\ast)$ at any time $t_\ast\in[t_0,t_f]$ and aims to minimize its travel cost during $t\in[t_\ast,t_f]$, what is the generic car's optimal control strategy and optimal cost value at time $t_\ast$?

We will solve the above subproblem in three cases: speed choice on a link, route choice at a node, and the terminal state, which give the HJB equation on a link, its spatial boundary condition at a node, and its terminal condition at the terminal state, respectively.


\textbf{Speed choice on a link.}\quad Suppose $X(t_\ast)=(x_\ast, l_\ast)$ where $0< x_\ast<\text{len}(l_\ast)$ and $0<t_\ast<T$. That is, the car is driving in the interior of a link. In this case the car needs to select its optimal speed at time $t_\ast$. We define $V_{l_\ast}(x_\ast,t_\ast)$ to be the optimal cost value of the car's trip from the state $(x_\ast, l_\ast)$ at time $t_\ast$. 

Denote the car's speed at time $t_\ast$ as $u\in[u_{\text{min}}, u_{\text{max}}]$.
We choose a small time step $\Delta t$ such that 
the car arrives at position $x_\ast+u\Delta t<\text{len}(l_\ast)$ on link $l_\ast$ at time $t_\ast+\Delta t<T$. The car's decision process is then decomposed into two stages: during the time interval $[t_\ast, t_\ast+\Delta t]$ the car's link running cost is $f_{\text{run}}\left(u, \rho_{l_\ast}(x_\ast,t_\ast)\right)\Delta t$; after time $t_\ast+\Delta t$, the car selects optimal driving speed and routing policy that give the optimal cost value $V_{l_\ast}(x_\ast+u\Delta t,t_\ast+\Delta t)$. Summarizing the costs at these two stages, we have:
\begin{align}
    V_{l_\ast}(x_\ast,t_\ast)=\min_{u_{\text{min}}\leq u\leq u_{\text{max}}} \left\{V_{l_\ast}(x_\ast+u\Delta t,t_\ast+\Delta t)+f_{\text{run}}\left(u, \rho_{l_\ast}(x_\ast,t_\ast)\right)\Delta t\right\}.\label{eq:derive_hjb_tmp}
\end{align}
Applying the first order Taylor's expansion on $V_{l_\ast}(x_\ast+u\Delta t,t_\ast+\Delta t)$ and letting $\Delta t\to0$, we obtain the following HJB equation from Eq.~\eqref{eq:derive_hjb_tmp}:
\begin{align}
	\partial_tV_{l_\ast}(x_\ast,t_\ast)+\min_{u_{\text{min}}\leq u\leq u_{\text{max}}}\left\{u\partial_xV_{l_\ast}(x_\ast,t_\ast)+f_{\text{run}}(u,\rho_{l_\ast}(x_\ast,t_\ast))\right\}=0.\label{eq:hjb1}
\end{align}
See \cite{huang2019game} for more details.
We denote the car's optimal speed as $u_{l_\ast}(x_\ast,t_\ast)$, it is given by the minimizer in Eq.~\eqref{eq:hjb1}:
\begin{align}
 	u_{l_\ast}(x_\ast,t_\ast)=\argmin_{u_{\text{min}}\leq u\leq u_{\text{max}}} \left\{u\partial_xV_{l_\ast}(x_\ast,t_\ast)+f_{\text{run}}(u,\rho_{l_\ast}(x_\ast,t_\ast))\right\}.\label{eq:hjb2}
\end{align}

Eqs.~\eqref{eq:hjb1}\eqref{eq:hjb2} are defined on all links of the network during the planning horizon $(0,T)$. We will replace $l_\ast$, $x_\ast$ and $t_\ast$ by the general link index $l\in\mathcal{L}$, space coordinate $x\in(0,\text{len}(l))$ and time coordinate $t\in(0,T)$. 


\textbf{Route choice at a node.}\quad Suppose $X(t_\ast)=i_\ast$ and $t_\ast\in R$, where $i_\ast\in\mathcal{N}\backslash\{s\}$ and $0< t_\ast<T$. That is, the car enters the queue at node $i_\ast$ at time $t_\ast$. 
In this case, the car exits the queue at time $e(t_\ast)=\min\{t_\ast+Q_{i_\ast}(t_\ast)/M_{i_\ast}, T\}$ and needs to select its next-go-to link at time $e(t_\ast)$. We define $\pi_{i_\ast}(e(t_\ast))$ to be the optimal cost value of the car's trip after exiting the queue at node $i_\ast$ at time $e(t_\ast)$.

By the assumption (A6), the next-go-to link selection is implemented in a mixed strategy. Remember that $\gamma(i_\ast,e(t_\ast))\in \text{OUT}(i_\ast)$ refers to the car's next-go-to link choice at node $i_\ast$ at time $e(t_\ast)$. But in the mixed strategy, there exists a probability distribution over all feasible choices. We define $\beta(i_\ast,l,e(t_\ast))\in[0,1]$ as the probability that the car selects the next-go-to link $l$ for all $l\in\text{OUT}(i_\ast)$. The probability distribution should satisfy:
\begin{align}
	\sum\nolimits_{l\in\text{OUT}(i_\ast)}\beta(i_\ast,l,e(t_\ast))=1.\label{eq:hjb3}
\end{align}
Once the car selects a link $l\in\text{OUT}(i_\ast)$, it moves to the starting point of the link at time $e(t_\ast)$ and its optimal cost value is $V_l(0,e(t_\ast))$. Taking the minimum over all $l\in\text{OUT}(i_\ast)$, the car's optimal cost value at the node is given by:
\begin{align}
	\pi_{i_\ast}(e(t_\ast))&=\min_{l\in\text{OUT}(i_\ast)}V_l(0,e(t_\ast)).\label{eq:hjb4}
\end{align}
The equilibrium condition is defined as:
\begin{align}
	\beta(i_\ast,l,e(t_\ast))\geq0,\quad\pi_{i_\ast}(e(t_\ast))=V_l(0,e(t_\ast))\text{ if }\beta(i_\ast,l,e(t_\ast))>0,\quad l\in\text{OUT}(i_\ast).\label{eq:hjb5}
\end{align}
That is, any used link $l\in\text{OUT}(i_\ast)$ is the optimal choice. This is consistent with the classical dynamic user equilibrium definition. Eqs.~(\ref{eq:hjb3}-\ref{eq:hjb5}) are defined on all non-destination nodes of the network during the planning horizon $(0,T)$. We will replace $i_\ast$ and $e(t_\ast)$ by the general node index $i\in\mathcal{N}\backslash\{s\}$ and time coordinate $t\in(0,T)$.

When $i_\ast$ is an intermediate node, the car moves to node $i_\ast$ from the end point of some link $h\in\text{IN}(i_\ast)$ at time $t_\ast$. With the existence of queuing delay at node $i_\ast$, the car's optimal cost value at the end point of link $h$ at time $t_\ast$, which is $V_h\left(\text{len}(h), t_\ast\right)$, is the nodal optimal cost value $\pi_{i_\ast}(e(t_\ast))$ plus the car's queuing cost $f_{\text{que}}(e(t_\ast)-t_\ast)$. 
Note that $e(t_\ast)=\min\{t_\ast+Q_{i_\ast}(t_\ast)/M_{i_\ast}, T\}$, we obtain: 
\begin{align}
    V_h(\text{len}(h),t_\ast)=\pi_{i_\ast}\left(\min\{t_\ast+\frac{Q_{i_\ast}(t_\ast)}{M_{i_\ast}},T\}\right)+f_{\text{que}}\left(\min\{\frac{Q_{i_\ast}(t_\ast)}{M_{i_\ast}},T-t_\ast\}\right).\label{eq:hjb6}
\end{align}
Eq.~\eqref{eq:hjb6} is defined on all intermediate nodes of the network during the planning horizon $(0,T)$. We will replace $i_\ast$ and $t_\ast$ by the general node index $i\in\mathcal{N}_I$ and time coordinate $t\in(0,T)$.

\textbf{Terminal state.}\quad When $t_\ast=t_f$, the car terminates its trip over the network at time $t_\ast$. There are three cases:
\begin{itemize}
    \item[(i)] The car arrives at the destination $s$ at time $t_\ast<T$. In this case the car's terminal cost gives the optimal cost value at the destination:
    \begin{align}
        V_l(\text{len}(l), t_\ast)=V_{\text{ter}}\left(t_\ast, s\right)=0,\quad l\in\text{IN}(s).
    \end{align}
    \item[(ii)] The car is in the interior of a link at time $t_\ast=T$. In this case the car's terminal cost gives the optimal cost value on a link:
    \begin{align}
        V_{l_\ast}(x_\ast,T)=V_{\text{ter}}\left(T, x_\ast, l_\ast\right),
    \end{align}
    where $l_\ast$ is the link the car is on and $x_\ast$ is the car's position on the link at time $t_\ast=T$.
    \item[(iii)] The car is in the queue at a non-destination node at time $t_\ast=T$. In this case the car's terminal cost gives the optimal cost value at a node:
    \begin{align}
        \pi_{i_\ast}(T)=V_{\text{ter}}\left(T, i_\ast\right),
    \end{align}
    where $i_\ast$ is the node at which the car is in the queue at time $t_\ast=T$.
\end{itemize}
Summarizing the three cases, the terminal conditions of the HJB equation include: (i) $V_l(\text{len}(l),t)=0$ for all $l\in\text{IN}(s)$ and $t\in(0,T)$; (ii) $V_l(x,T)$ for all $l\in\mathcal{L}$ and $x\in(0,\text{len}(l))$; (iii) $\pi_i(T)$ for all $i\in\mathcal{N}\backslash\{s\}$.
We may choose all $V_l(x,T)$ and $\pi_i(T)$ to be the same constant if we do not care about the cars' final positions. Or we can properly define those terminal costs to penalize cars who cannot arrive at the destination in time.
	
The HJB equation defined by  Eqs.~(\ref{eq:hjb1}-\ref{eq:hjb6}) 
characterizes a generic car's optimal speed control and routing policy. The optimal velocity fields $u_l(x,t)$ for $l\in\mathcal{L}$, $x\in\left(0,\text{len}(l)\right)$, $t\in(0,T)$ and next-go-to link distributions $\beta(i,l,t)$ for $i\in\mathcal{N}\backslash\{s\}$, $l\in\text{OUT}(i)$, $t\in(0,T)$ are solved from the HJB equation, 
provided the link traffic densities and node queue sizes.
The HJB equation is solved backward in time. 

\subsection{Forward continuity equation}\label{sec:continuity}

When all cars follow the optimal speed control and routing policy, the traffic evolution on the network follows from their dynamic motions.
Suppose that: (i) the velocity field $u_l(x,t)$ is known for all $l\in\mathcal{L}$, $x\in\left(0,\text{len}(l)\right)$ and $t\in(0,T)$; (ii) the next-go-to link distribution $\beta(i,l,t)$ is known for all $i\in\mathcal{N}\backslash\{s\}$, $l\in\text{OUT}(i)$ and $t\in(0,T)$.
We derive the continuity equation characterizing the evolution of traffic state on the network.

Our derivation is based on the conservation of cars for both link flow propagation and node flow assignment. We will also introduce the initial condition of the continuity equation.

\textbf{Link flow propagation.} In the interior of any link $l\in\mathcal{L}$, 
the conservation of cars is described by the following partial differential equation of traffic density $\rho_l(x,t)$ and flux $q_l(x,t)$:  
\begin{align}
	\partial_t\rho_l(x,t)+\partial_xq_l(x,t)=0,\quad x\in(0,\text{len}(l)),\ t\in(0,T).\label{eq:cont1} 	
\end{align}
The flux $q_l(x,t)$ is computed from traffic density $\rho_l(x,t)$ and velocity field $u_l(x,t)$ as follows:
\begin{align}
	q_l(x,t)=\rho_l(x,t)u_l(x,t),\quad x\in(0,\text{len}(l)],\ t\in(0,T).\label{eq:cont2}
\end{align}
Eqs.~\eqref{eq:cont1}\eqref{eq:cont2} give the continuity equation on a link.


\textbf{Node flow assignment.} 
At any node $i\in\mathcal{N}\backslash\{s\}$, the conservation of cars  is described by the relation:
\begin{align*}
    \text{rate of change of queue size} = \text{total incoming flow} - \text{total outgoing flow}.
\end{align*}
Based on the relation, we are going to compute the total outgoing flow from the rate of change of queue size and the total incoming flow.

Knowing the queue size $Q_i(t)$ at node $i$ for $t\in[0,T]$, the rate of change of the queue size is given by $dQ_i(t)/dt$.
The total incoming flow is composed of two parts: the sum of flows from the node's all  incoming links and the traffic demand if the node is an origin. The contribution from the former part is given by $\sum\nolimits_{h\in\text{IN}(i)}q_h(\text{len}(h),t)$; For the latter part, we denote $d_i(t)$ the traffic demand on node $i$ at time $t$ and make the convention that $d_i(t)=0$ for all $t\in[0,T]$ if $i\notin\mathcal{N}_O$. Then the total incoming flow at node $i$ at time $t$ is given by $\sum\nolimits_{h\in\text{IN}(i)}q_h(\text{len}(h),t)+d_i(t)$. Knowing the rate of change of queue size and the total incoming flow,
the total outgoing flow at node $i$ at time $t$ is given by $\sum\nolimits_{h\in\text{IN}(i)}q_h(\text{len}(h),t)+d_i(t)-\frac{dQ_i(t)}{dt}$.

The total outgoing flow distributes to the node's all outgoing links. Remember that $\beta(i,l,t)$, $l\in\text{OUT}(i)$ is the probability distribution over a single car's next-go-to link choices from node $i$ at time $t$. With a large number of cars, $\beta(i,l,t)$ can be understood as the percentage of cars selecting the outgoing link $l$ among all cars arriving at node $i$ at time $t$.
According to the outgoing flow distribution, we obtain:
\begin{align}
	q_l(0,t)=\beta(i,l,t)\left(\sum\nolimits_{h\in\text{IN}(i)}q_h(\text{len}(h),t)+d_i(t)-\frac{dQ_i(t)}{dt}\right),\quad l\in\text{OUT}(i),\ t\in(0,T).\label{eq:cont3}
\end{align}


By the assumption (A5), the total outgoing flow is no greater than $M_i$ and it equals $M_i$ as long as $Q_i(t)>0$, which gives the following equations characterizing the evolution of the queue size:
\begin{align}
    \sum\nolimits_{h\in\text{IN}(i)}q_h(\text{len}(h),t)+d_i(t)-\frac{dQ_i(t)}{dt}&\leq M_i,\quad Q_i(t)\geq0,\quad t\in(0,T);\label{eq:cont4}\\
    \sum\nolimits_{h\in\text{IN}(i)}q_h(\text{len}(h),t)+d_i(t)-\frac{dQ_i(t)}{dt}&=M_i,\quad\text{if }Q_i(t)>0,\quad t\in(0,T).\label{eq:cont5}
\end{align}
Eqs.~(\ref{eq:cont3}-\ref{eq:cont5}) incorporate the Vickrey junction model \citep{vickrey1969congestion} and give the spatial boundary condition of the continuity equation. 

\textbf{Initial condition.} The network is empty when $t=0$. Hence the initial condition of the continuity equation is given by:
\begin{align}
     \rho_l(x,0)&=0,\quad l\in\mathcal{L},\  x\in(0,\text{len}(l));\\
     Q_i(0)&=0,\quad i\in\mathcal{N}.
\end{align}

The continuity equation defined by Eqs.~(\ref{eq:cont1}-\ref{eq:cont5}) characterizes the evolution of traffic state on the network. The link traffic densities $\rho_l(x,t)$ for $l\in\mathcal{L}$, $x\in(0,\text{len}(l))$, $t\in(0,T)$ and the node queue sizes $Q_i(t)$ for $i\in\mathcal{N}$, $t\in(0,T)$ are solved from the continuity equation, provided the optimal velocity fields and next-go-to link distributions. The continuity equation is solved forward in time.

\subsection{Mean field game system}\label{sec:mfg_system}

Summarizing all the equations (\ref{eq:hjb1}-\ref{eq:hjb6},\ref{eq:cont1}-\ref{eq:cont5}),
we obtain the following mean field game system on a network.

\begin{subequations}
\textbf{[MFGnet]}
\begin{align}
	\mbox{(CE-Link)}\quad & \partial_t\rho_l(x,t)+\partial_xq_l(x,t)=0,\quad l\in\mathcal{L};\label{eq:mfg_cont1}\\
	& q_l(x,t)=\rho_l(x,t)u_l(x,t),\quad l\in\mathcal{L};\label{eq:mfg_cont2}\\
	\mbox{(CE-Node)}\quad& q_l(0,t)=\beta(i,l,t)\left(\sum_{h\in\text{IN}(i)}q_h(\text{len}(h),t)+d_i(t)-\frac{dQ_i(t)}{dt}\right),\ i\in\mathcal{N}\backslash\{s\},\  l\in\text{OUT}(i); \label{eq:mfg_cont3}\\
	& \sum_{h\in\text{IN}(i)}q_h(\text{len}(h),t)+d_i(t)-\frac{dQ_i(t)}{dt}\leq M_i,\quad Q_i(t)\geq0\quad i\in\mathcal{N}\backslash\{s\};\label{eq:mfg_cont4}\\
	& \sum_{h\in\text{IN}(i)}q_h(\text{len}(h),t)+d_i(t)-\frac{dQ_i(t)}{dt}=M_i,\quad \text{if } Q_i(t)>0,\quad i\in\mathcal{N}\backslash\{s\};\label{eq:mfg_cont5}\\
	\mbox{(HJB-Link)}\quad & \partial_tV_{l}(x,t)+\min_{u_{\text{min}}\leq u\leq u_{\text{max}}}\left\{u\partial_xV_{l}(x,t)+f_{\text{run}}(u,\rho_{l}(x,t))\right\}=0,\quad l\in\mathcal{L};\label{eq:mfg_hjb1}\\
	& u_l(x,t)=\argmin_{u_{\text{min}}\leq u\leq u_{\text{max}}}\left\{u\partial_xV_{l}(x,t)+f_{\text{run}}(u,\rho_{l}(x,t))\right\},\quad l\in\mathcal{L};\label{eq:mfg_hjb2}\\
	\mbox{(HJB-Node)}\quad & V_h(\text{len}(h),t)=\pi_{i}\left(\min\{t+\frac{Q_{i}(t)}{M_{i}},T\}\right)+f_{\text{que}}\left(\min\{\frac{Q_{i}(t)}{M_{i}},T-t\}\right),\ i\in\mathcal{N}_I,\  h\in\text{IN}(i);\label{eq:mfg_hjb3}\\
	& \sum_{l\in\text{OUT}(i)}\beta(i,l,t)=1,\quad i\in\mathcal{N}\backslash\{s\};\label{eq:mfg_hjb4}\\
	& \pi_i(t)=\min_{l\in\text{OUT}(i)}V_l(0,t),\quad i\in\mathcal{N}\backslash\{s\};\label{eq:mfg_hjb5}\\
	& \beta(i,l,t)\geq0,\quad\pi_i(t)=V_l(0,t)\text{ if }\beta(i,l,t)>0,\quad  i\in\mathcal{N}\backslash\{s\},\ l\in\text{OUT}(i).\label{eq:mfg_hjb6}
\end{align}
\end{subequations}

[MFGnet] is a coupled system of the forward continuity equation and the backward HJB equation defined on a network. The following conditions are predetermined for [MFGnet]:
\begin{itemize}
    \item The initial traffic densities $\rho_l(x,0)=0$ for all $l\in\mathcal{L}$ and $x\in(0,\text{len}(l))$;
    \item The initial queue sizes $Q_i(0)=0$ for all $i\in\mathcal{N}\backslash\{s\}$;
    \item The terminal costs at the destination $V_l(\text{len}(l),t)=0$ for all $l\in\text{IN}(s)$ and $t\in(0,T)$;
\end{itemize}
and the following conditions need to be specified:
\begin{itemize}
    \item The traffic demands $d_i(t)$ for all $i\in\mathcal{N}_O$ and $t\in[0,T]$.
	\item The terminal costs at the final time $V_l(x,T)$ for all $l\in\mathcal{L}$ and $x\in(0,\text{len}(l))$, and $\pi_i(T)$ for all $i\in\mathcal{N}\backslash\{s\}$.
\end{itemize}

The mean field equilibrium condition is defined as the following: the evolution of link traffic densities and node queue sizes resulting from the optimal velocity fields and route choices is consistent with the actual traffic evolution.
The solution of [MFGnet], which we will refer to as SOL([MFGnet]), gives the mean field equilibrium. The mean field equilibrium contains the following:
\begin{itemize}
    \item $\rho_l(x,t)$, $u_l(x,t)$ and $V_l(x,t)$ for all $l\in\mathcal{L}$, $x\in(0,\text{len}(l))$ and $t\in[0,T]$;
    \item $\pi_i(t)$ and $Q_i(t)$ for all $i\in\mathcal{N}\backslash\{s\}$ and $t\in[0,T]$;
    \item $\beta(i,l,t)$ for all $i\in\mathcal{N}\backslash\{s\}$, $l\in\text{OUT}(i)$ and $t\in[0,T]$.
\end{itemize}

The existence of equilibrium solutions depends on both the network topology and cost functions. The following proposition gives a sufficient condition on the solution existence of [MFGnet].
\begin{prop}
    Suppose the network contains only one link, [MFGnet] becomes a mean field game system defined on a single road with only the speed control. In this case, suppose the AVs' running cost function $f_{\text{run}}(u,\rho)=g(u)+h(\rho)$, where $g(u)$ is a strictly convex function of speed $u$, $h(\rho)$ is a strictly increasing function of density $\rho$, and $g(u)$ and $h(\rho)$ satisfy certain growth conditions. Then there exists a unique weak solution of [MFGnet] \citep{cardaliaguet2015weak}.
\end{prop}

\subsection{Reformulation as a mixed complementarity problem}\label{sec:reform}


To facilitate the algorithm development and give meaning to each equation, in this subsection, we will reformulate [MFGnet] as a mixed complementarity problem.

First, the Vickrey equations \eqref{eq:mfg_cont4}\eqref{eq:mfg_cont5} can be reformulated as the following linear complementarity relation \citep{ban2012continuous}:
\begin{align}
    0\leq \frac{dQ_i(t)}{dt}+M_i-\sum_{h\in\text{IN}(i)}q_h(\text{len}(h),t)-d_i(t)\perp Q_i(t)\geq0.\label{eq:reform1}
\end{align}

Then let us reformulate the equation \eqref{eq:mfg_hjb4} as a complementarity relation. We will choose all cost functions to be positive, then the cars' travel costs are always positive and we have $\pi_i(t)>0$ for all  $i\in\mathcal{N}\backslash\{s\}$ and $t\in[0,T]$. 
In this case,
Eq.~\eqref{eq:mfg_hjb4} is equivalent to the following complementarity relation:
\begin{align}
    0\leq\pi_{i}(t)\perp\sum\nolimits_{l\in\text{OUT}(i)}\beta(i,l,t)-1\geq0,\label{eq:reform2}
\end{align}
because the right hand side term $\sum\nolimits_{l\in\text{OUT}(i)}\beta(i,l,t)-1$ must be zero when $\pi_i(t)>0$.

The equilibrium condition defined by  Eqs.~\eqref{eq:mfg_hjb5}\eqref{eq:mfg_hjb6} can be rewritten as the following complementarity relation:
\begin{align}
    0\leq\beta(i,l,t)\perp V_l(0,t)-\pi_{i}(t)\geq0,\quad \forall l\in\text{OUT}(i).\label{eq:reform3}
\end{align}



Using Eqs.~(\ref{eq:reform1}-\ref{eq:reform3}), [MFGnet] can be reformulated as a new system, denoted as [MFG-MiCP], defined in Eq.~\eqref{eq:mfg-micp} below. 
It is a mixed complementarity problem (MiCP) composed of partial differential equations, algebraic equations, and complementarity relations \citep{facchinei2007finite,ban2008link,ban2012continuous}.

\begin{subequations}
\textbf{[MFG-MiCP]}
\begin{align}
	\mbox{[link flow balance]}\quad &\partial_t\rho_l(x,t)+\partial_xq_l(x,t)=0,\quad l\in\mathcal{L};\label{eq:due_propagation}\\
	\mbox{[link flow propagation]}\quad &q_l(x,t)=\rho_l(x,t)u_l(x,t),\quad l\in\mathcal{L};\label{eq:due_flow}\\
	\mbox{[link influx]}\quad
	&q_l(0,t)=\beta(i,l,t)\left(\sum_{h\in\text{IN}(i)}q_h(\text{len}(h),t)+d_i(t)-\frac{dQ_i(t)}{dt}\right),\notag\\ 
	&i\in\mathcal{N}\backslash\{s\},
	\ l\in\text{OUT}(i); \label{eq:due_assign}\\
	\mbox{[nodal delay]}\quad &0\leq \frac{dQ_i(t)}{dt}+M_i-\sum_{h\in\text{IN}(i)}q_h(\text{len}(h),t)-d_i(t)\perp Q_i(t)\geq0,\quad i\in\mathcal{N}\backslash\{s\};\label{eq:due_vickery}\\
	\mbox{[optimal link cost]}\quad &\partial_tV_{l}(x,t)+\min_{u_{\text{min}}\leq u\leq u_{\text{max}}}\left\{u\partial_xV_{l}(x,t)+f_{\text{run}}(u,\rho_{l}(x,t))\right\}=0,\quad l\in\mathcal{L};\label{eq:due_cost_link}\\
	\mbox{[optimal speed]}\quad & u_l(x,t)=\argmin_{u_{\text{min}}\leq u\leq u_{\text{max}}}\left\{u\partial_xV_{l}(x,t)+f_{\text{run}}(u,\rho_{l}(x,t))\right\},\quad l\in\mathcal{L};\label{eq:due_speed}\\
    \mbox{[optimal nodal cost]} \quad & V_h(\text{len}(h),t)=\pi_{i}\left(\min\{t+\frac{Q_{i}(t)}{M_{i}},T\}\right)+f_{\text{que}}\left(\min\{\frac{Q_{i}(t)}{M_{i}},T-t\}\right),\notag\\ 
    & i\in\mathcal{N}_I,\  h\in\text{IN}(i);\label{eq:due_cost_node}\\
	\mbox{[nodal flow conservation]}\quad &0\leq\pi_i(t)\perp\sum_{l\in\text{OUT}(i)}\beta(i,l,t)-1\geq0,\quad i\in\mathcal{N}\backslash\{s\};\label{eq:due_conservation}\\
	\mbox{[equilibrium turning ratio]}\quad &0\leq\beta(i,l,t)\perp V_l(0,t)-\pi_i(t)\geq0,\quad i\in\mathcal{N}\backslash\{s\},\  l\in\text{OUT}(i).\label{eq:due_equilibrium}
\end{align}\label{eq:mfg-micp}
\end{subequations}




To solve [MFG-MiCP], we need the same conditions (initial traffic densities and queue sizes, terminal costs, and traffic demands) as those for [MFGnet], which are introduced in Section~\ref{sec:mfg_system}.

Now we group the above nine equations and explain their meanings: 
\begin{enumerate}
    \item Flow on links (corresponding to ``dynamic loading" in DUE): 
    Eq.~\eqref{eq:due_propagation} describes the flow balance on each link;
    Eq.~\eqref{eq:due_flow} prescribes the fundamental relation between traffic density, driving speed and traffic flow;
    Eq.~\eqref{eq:due_assign} describes link influx from each node;
    and Eq.~\eqref{eq:due_vickery} describes the nodal delay at each node.
    \item Flow at nodes (corresponding to ``route choice" in DUE): 
    Eq.~\eqref{eq:due_cost_node} describes how cars' optimal travel costs are determined at each node;
    Eq.~\eqref{eq:due_conservation} guarantees flow conservation at each node;
    and Eq.~\eqref{eq:due_equilibrium} is the equilibrium condition of cars' next-go-to link choices at each node.
    \item Speed selection: 
    Eq.~\eqref{eq:due_cost_link} and Eq.~\eqref{eq:due_speed} describe how cars' optimal travel costs and optimal driving speeds are solved from traffic density in the interior of each link. 
\end{enumerate}

\begin{note}
The first two groups of equations are projected to two components of classical DUE, 
which are dynamic loading and route choice. The connection between MFE and DUE will be further established in the next section. 
\end{note}

We denote the solution of [MFG-MiCP] as SOL([MFG-MiCP]), which is the same as SOL([MFGnet]).
In the rest of the paper, we will focus on the system  [MFG-MiCP]. Its solution SOL([MFG-MiCP]) describes the mean field equilibrium (MFE), which is the equilibrium dynamic controls of all AVs.
\color{black}

\section{Connection between MFE and DUE}\label{sec:connection}

\begin{figure}[htbp]
	\centering
	\includegraphics[width=.7\textwidth]{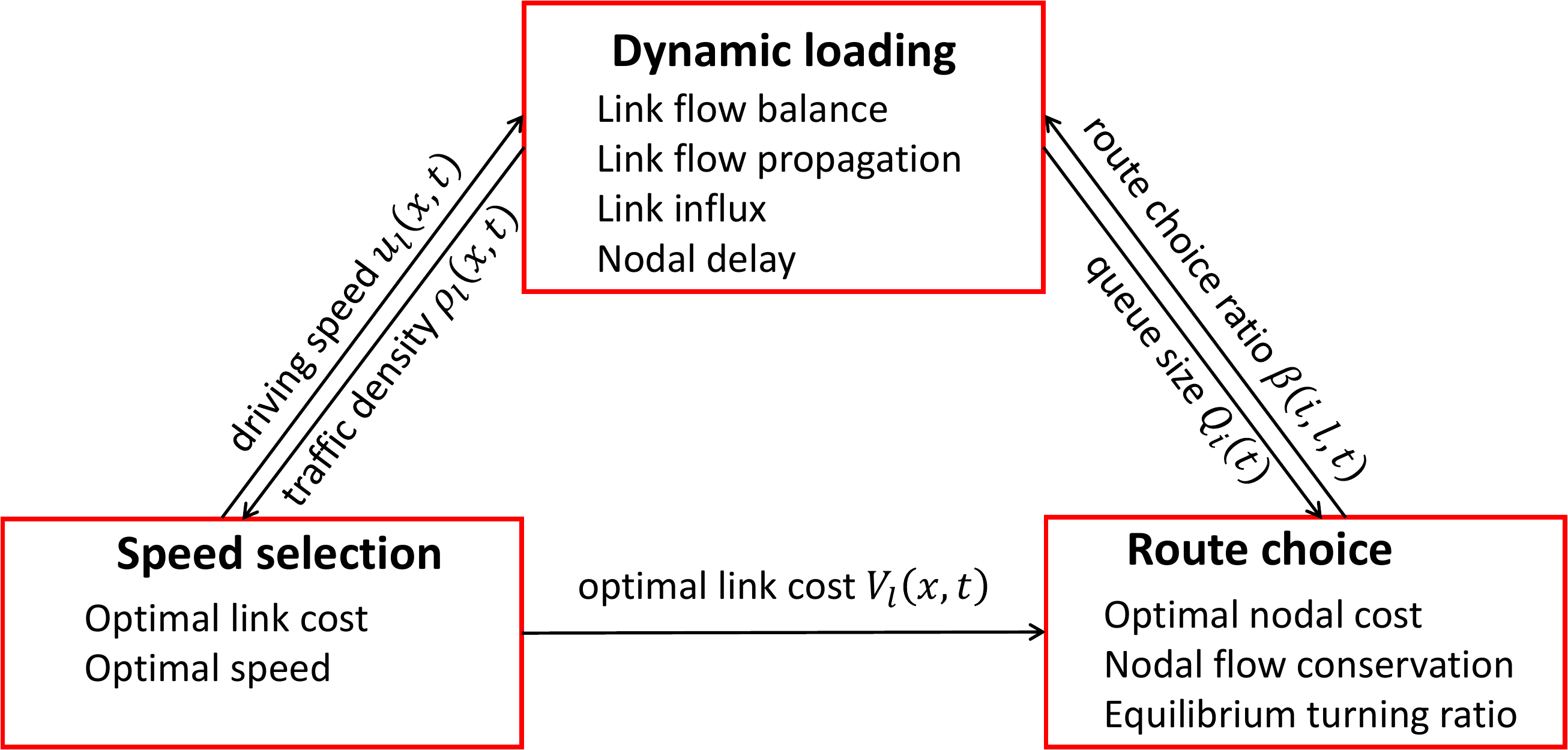}
	\caption{Components of MFG on networks}\label{fig:connection}
\end{figure}

The system [MFG-MiCP] can be divided into three coupled components: dynamic loading, route choice, and speed selection. 
The information exchange among these three components is demonstrated in Figure~\ref{fig:connection}.
In the ``dynamic loading" procedure, 
given cars' driving speeds $u_l(x,t), \forall l$ and next-go-to link choices $\beta(i,l,t), \forall i, l$, the traffic densities $\rho_l(x,t), \forall l$ and queue sizes $Q_i(t), \forall i$ are solved from the equations [link flow balance][link flow propagation][link influx][nodal delay], i.e., Eqs.~(\ref{eq:due_propagation}-\ref{eq:due_vickery}); 
In the ``route choice'' procedure, given the queue sizes $Q_i(t), \forall i$ and cars' optimal link costs $V_l(x,t), \forall l$, the cars' next-go-to link choices $\beta(i,l,t), \forall i, l$ are solved from the equations [optimal nodal cost][nodal flow conservation][equilibrium turning ratio], i.e., Eqs.~(\ref{eq:due_cost_node}-\ref{eq:due_equilibrium}); 
In the ``speed selection" procedure, 
given the traffic densities $\rho_l(x,t), \forall l$, 
the cars' driving speeds $u_l(x,t), \forall l$ and optimal link costs $V_l(x,t), \forall l$ are solved from the equations [optimal link cost][optimal speed], i.e., Eqs.~(\ref{eq:due_cost_link},\ref{eq:due_speed}). 

The classical DUE concept characterizes a large number of cars' equilibrium route choices. It does not incorporate cars' driving speed controls so only the dynamic loading and route choice components are coupled.

If we assume that cars do not optimize their driving speeds but choose the equilibrium speeds characterized by the LWR model \citep{lighthill1955kinematic,richards1956shock}:
\begin{align}
    u_l(x,t)=U(\rho_l(x,t)),\quad l\in\mathcal{L},\ x\in(0,\text{len}(l)),\ t\in(0,T).\label{eq:lwr_speed}
\end{align}
That is, the driving speed of every car on the network is totally determined by its nearby local density. 
Accordingly, we present a special case of MFG using the LWR speed, denoted as [MFG-LWR-MiCP], which is formulated below.

\begin{subequations}
\textbf{[MFG-LWR-MiCP]}
\begin{align}
	\mbox{[link flow balance]}\quad &\partial_t\rho_l(x,t)+\partial_xq_l(x,t)=0,\quad l\in\mathcal{L};\label{eq:lwr_propagation}\\
	\mbox{[link flow propagation]}\quad &q_l(x,t)=\rho_l(x,t)u_l(x,t),\quad l\in\mathcal{L};\label{eq:lwr_flow}\\
	\mbox{[link influx]}\quad
	&q_l(0,t)=\beta(i,l,t)\left(\sum_{h\in\text{IN}(i)}q_h(\text{len}(h),t)+d_i(t)-\frac{dQ_i(t)}{dt}\right),\notag\\ 
	&i\in\mathcal{N}\backslash\{s\},
	\ l\in\text{OUT}(i); \label{eq:lwr_assign}\\
	\mbox{[nodal delay]}\quad &0\leq \frac{dQ_i(t)}{dt}+M_i-\sum_{h\in\text{IN}(i)}q_h(\text{len}(h),t)-d_i(t)\perp Q_i(t)\geq0,\quad i\in\mathcal{N}\backslash\{s\};\label{eq:lwr_vickery}\\
	\mbox{[optimal link cost]}\quad &\partial_tV_{l}(x,t)+\min_{u_{\text{min}}\leq u\leq u_{\text{max}}}\left\{u\partial_xV_{l}(x,t)+f_{\text{run}}(u,\rho_{l}(x,t))\right\}=0,\quad l\in\mathcal{L};\label{eq:lwr_cost_link}\\
	\mbox{[LWR speed]}\quad &u_l(x,t)=U\left(\rho_l(x,t)\right),\quad l\in\mathcal{L};\label{eq:lwr_optimal_speed}\\
    \mbox{[optimal nodal cost]} \quad & V_h(\text{len}(h),t)=\pi_{i}\left(\min\{t+\frac{Q_{i}(t)}{M_{i}},T\}\right)+f_{\text{que}}\left(\min\{\frac{Q_{i}(t)}{M_{i}},T-t\}\right),\notag\\ 
    & i\in\mathcal{N}_I,\  h\in\text{IN}(i);\label{eq:lwr_cost_node}\\
	\mbox{[nodal flow conservation]}\quad &0\leq\pi_i(t)\perp\sum_{l\in\text{OUT}(i)}\beta(i,l,t)-1\geq0,\quad i\in\mathcal{N}\backslash\{s\};\label{eq:lwr_conservation}\\
	\mbox{[equilibrium turning ratio]}\quad &0\leq\beta(i,l,t)\perp V_l(0,t)-\pi_i(t)\geq0,\quad i\in\mathcal{N}\backslash\{s\},\  l\in\text{OUT}(i).\label{eq:lwr_equilibrium}
\end{align}\label{eq:mfg-lwr-micp}
\end{subequations}

In [MFG-LWR-MiCP], one can choose any running cost function $f_{\text{run}}(\cdot)$, not restricted to the total travel time as in classical DUE. Even under the same cost function, the optimal speed solved from the HJB equation may have a large deviation from the LWR speed defined in Eq.~\eqref{eq:lwr_speed}. 
We may solve very different equilibria from [MFG-MiCP] and [MFG-LWR-MiCP]. The two equilibria will be compared on Braess networks in Section~\ref{sec:numerical}.

\section{Solution algorithm}\label{sec:algo}

In this section, we will introduce the solution algorithm for solving the mean field equilibrium. We will first discretize the [MFG-MiCP] system on a mesh grid. The discretized system is a finite dimensional MiCP. Then we can solve this discretized MiCP using the PATH solver built in GAMS \citep{rosen2007gams}.



Denote $\Delta x$ and $\Delta t$ the spatial and temporal mesh sizes. 
For the time discretization, we divide the whole time horizon $[0,T]$ into $N_t+1$ time steps:
\begin{align}
	0=t_0<t_1<\cdots<t_{N_t}=T,\label{eq:time_grid}
\end{align}
where $t_k=k\Delta t$ ($k=0,1,\dots,N_t)$, and we assume that $N_t=T/\Delta t$ is a positive integer.

\begin{figure}[htbp]
	\centering
	\begin{subfigure}{.48\textwidth}
	\centering
	\begin{tikzpicture}[scale=.4]
        \begin{scope}[every node/.style={circle,thick,draw}]
            \node (s) at (-6,0) {1};
            \node (d) at (6,0) {4};
            \node (1) at (0,3) {2};
            \node (2) at (0,-3) {3};
            \node (3) at (-3,1.5) {5};
            \node (4) at (-3,-1.5) {6};
            \node (5) at (3,1.5) {7};
            \node (6) at (3,-1.5) {8};
        \end{scope}
        \begin{scope}[>=Stealth,auto,
                      every edge/.style={draw=black,very thick}]
            \path [->] (s) edge node {} (3);
            \path [->] (3) edge node {} (1);
            \path [->] (s) edge node {} (4);
            \path [->] (4) edge node {} (2);
            \path [->] (1) edge node {} (5);
            \path [->] (5) edge node {} (d);
            \path [->] (2) edge node {} (6);
            \path [->] (6) edge node {} (d);
        \end{scope}
    \end{tikzpicture}
    \label{fig:disnet1}
    \end{subfigure}
    \begin{subfigure}{.48\textwidth}
    \centering
	\begin{tikzpicture}[scale=.4]
        \begin{scope}[every node/.style={circle,thick,draw}]
            \node (s) at (-6,0) {1};
            \node (d) at (6,0) {4};
            \node (1) at (0,3) {2};
            \node (2) at (0,-3) {3};
            \node (3) at (-3,1.5) {5};
            \node (4) at (-3,-1.5) {6};
            \node (5) at (3,1.5) {7};
            \node (6) at (3,-1.5) {8};
            \node (7) at (0,0) {9};
        \end{scope}
        \begin{scope}[>=Stealth,auto,
                      every edge/.style={draw=black,very thick}]
            \path [->] (s) edge node {} (3);
            \path [->] (3) edge node {} (1);
            \path [->] (s) edge node {} (4);
            \path [->] (4) edge node {} (2);
            \path [->] (1) edge node {} (5);
            \path [->] (5) edge node {} (d);
            \path [->] (2) edge node {} (6);
            \path [->] (6) edge node {} (d);
            \path [->] (1) edge node {} (7);
            \path [->] (7) edge node {} (2);
        \end{scope}
    \end{tikzpicture}
    \label{fig:disnet2}
    \end{subfigure}
    \caption{Discretized Braess networks, all links have length 1 and $\Delta x=0.5$}
    \label{fig:disnet}
\end{figure}
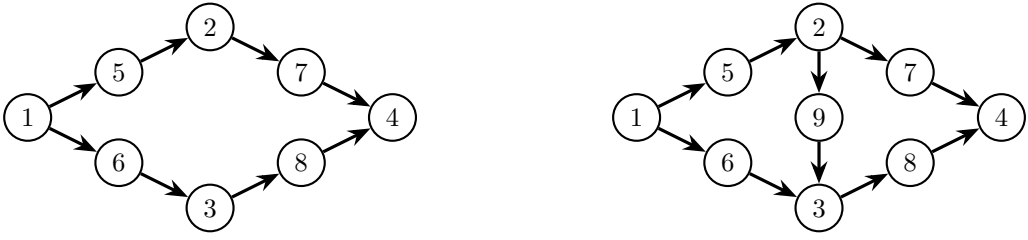

For the spatial discretization, we divide each link $l=(i,j)\in\mathcal{L}$ into $N_x^l=\text{len}(l)/\Delta x$ sublinks (we will always assume $N_x^l$ is a positive integer):
\begin{align}
	(i,i_1),\ (i_1,i_2),\ \cdots,\ (i_{N_x^l-2},i_{N_x^l-1}),\ (i_{N_x^l-1},j).
\end{align}
As an example, Figure~\ref{fig:disnet} shows the discretization of the Braess networks with four and five links. 
It is assumed that all links in the networks have length 1 and $\Delta x=0.5$.

We observe that the spatial discretization of a network ${\cal G}=(\mathcal{N},\mathcal{L})$ is equivalent to adding auxiliary nodes into $\mathcal{N}$ and replacing each link $l\in\mathcal{L}$ by a group of sublinks. Denote the set of 
auxiliary nodes by $\mathcal{N}'$ 
and the set of all nodes $\mathcal{N}^D=\mathcal{N}\cup{\mathcal{N}'}$. Denote the set of all sublinks by $\mathcal{L}^D$. The discretized network is then represented by ${\cal G}^D=(\mathcal{N}^D,\mathcal{L}^D)$. The discretized network $\mathcal{G}^D$ keeps the same topological structure as the original network $\mathcal{G}$ but has more nodes and links.

Now we are ready to discretize all variables and equations in [MFG-MiCP] on the discretized network $\mathcal{G}^D$ and at the time steps defined in Eq.~\eqref{eq:time_grid}. We define:
\begin{itemize}
	\item $\rho_{ij}^k$: the average traffic density on link $(i,j)\in\mathcal{L}^D$ at time $t_k$, $k=0,1,\dots,N_t$.
	\item $p_{ij}^k$: the entry flow rate entering link $(i,j)\in\mathcal{L}^D$ from node $i\in\mathcal{N}^D\backslash\{s\}$ at time $t_k$, $k=0,1,\dots,N_t-1$.
	\item $q_{ij}^k$: the exit flow rate leaving link $(i,j)\in\mathcal{L}^D$ from node $j\in\mathcal{N}^D$ at time $t_k$, $k=0,1,\dots,N_t-1$.
	\item $u_{ij}^k$: the average driving speed of cars on link $(i,j)\in\mathcal{L}^D$ at time $t^k$, $k=0,1,\dots,N_t-1$.
	\item $\beta_{ij}^k$: the percentage of cars selecting link $(i,j)\in\mathcal{L}^D$ among all cars arriving at node $i\in\mathcal{N}^D\backslash\{s\}$ at time $t^k$, $k=0,1,\dots,N_t-1$.
	\item $V_{ij}^k$: the optimal travel cost if a car exits the queue at node $i\in\mathcal{N}^D\backslash\{s\}$ and selects the link $(i,j)\in\mathcal{L}^D$ at time $t^k$, $k=0,1,\dots,N_t$.
	\item $\pi_i^k$: the optimal travel cost if a car exits the queue at node $i\in\mathcal{N}^D\backslash\{s\}$ and selects the minimum-cost link from node $i$ at time $t^k$, $k=0,1,\dots,N_t$.
	\item $\lambda_i^k$: the optimal travel cost if a car enters the queue at node $i\in\mathcal{N}^D$ at time $t^k$, $k=0,1,\dots,N_t$.
	\item $Q_i^k$: the queue size at node $i\in\mathcal{N}^D\backslash\{s\}$ at time $t^k$, $k=0,1,\dots,N_t$.
\end{itemize}

Provided the listed discrete variables, we then discretize each equation in [MFG-MiCP]:
\begin{itemize}
    \item The equations [link flow balance][optimal link cost][optimal speed] are discretized by the upwinding scheme \citep{huang2019game}. Note that the driving speed is in the range $[u_{\text{min}},u_{\text{max}}]$, the choice of the spatial and temporal mesh sizes should satisfy the CFL condition $u_{\text{max}}\Delta t\leq\Delta x$ to guarantee numerical stability \citep{leveque2002finite}.
    \item In the equation [link influx], the term $\frac{dQ_i(t)}{dt}$ is discretized as $\frac{Q_i^{k+1}-Q_i^k}{\Delta t}$ at time $t^k$.
    \item For the equation [nodal delay], the forward Euler discretization is not well-defined \citep{han2013partial1}.
    Here we use the following backward Euler scheme to discretize the equation [nodal delay], which is suggested by \cite{ban2012continuous}:
    \begin{align}
        0\leq \frac{Q_i^{k+1}-Q_i^k}{\Delta t}+M_i-\sum\nolimits_{(m,i)\in\mathcal{L}^D}q_{mi}^k-d_i^k\perp Q_i^{k+1}\geq0.
    \end{align}
    \item In [optimal nodal cost], the term $\pi_i\left(\min\{t+\frac{Q_i(t)}{M_i},T\}\right)+f_{\text{que}}\left(\min\{\frac{Q_i(t)}{M_i},T-t\}\right)$ at time $t^k$ is discretized using the following interpolation scheme \citep{ban2008link}:
    \small{
    \begin{align}
        \sum_{k'=k}^{N_t-1}&\mathbf{1}_{k'\leq k+Q_i^k/(M_i\Delta t)<k'+1}\cdot\left[\left(k'+1-\frac{Q^k_i}{M_i\Delta t}-k\right)\pi_i^{k'}+\left(\frac{Q^k_i}{M_i\Delta t}+k-k'\right)\pi_i^{k'+1}\right]+\mathbf{1}_{k+Q_i^k/(M_i\Delta t)\geq N_t}\cdot\pi_i^{N_t}\notag\\
        &+f_{\text{que}}\left(\min\{\frac{Q_i^k}{M_i},T-k\Delta t\}\right).
    \end{align}
    }
    \item The equations [link flow propagation][nodal flow conservation][equilibrium turning ratio] do not involve space or time derivatives, therefore they can be directly discretized.
\end{itemize}

Summarizing all above, the discretized system of [MFG-MiCP] is:
\begin{subnumcases}{}
    &$\rho_{ij}^{k+1}=\rho_{ij}^k+\frac{\Delta t}{\Delta x}(p_{ij}^k-q_{ij}^k)$,\quad $(i,j)\in\mathcal{L}^D$,\, $0\leq k\leq N_t-1$;\label{eq:micp1}\\
    &$q_{ij}^k=\rho_{ij}^ku_{ij}^k$,\quad $(i,j)\in\mathcal{L}^D$,\, $0\leq k\leq N_t-1$;\label{eq:mip2}\\
    &$p_{ij}^k=\beta_{ij}^k\left(\sum_{(m,i)\in\mathcal{L}^D}q_{mi}^k+d_i^k-\frac{Q_i^{k+1}-Q_i^k}{\Delta t}\right)$,\quad $(i,j)\in\mathcal{L}^D$,\, $0\leq k\leq N_t-1$;\label{eq:micp3}\\
    &$0\leq \frac{Q_i^{k+1}-Q_i^k}{\Delta t}+M_i-\sum_{(m,i)\in\mathcal{L}^D}q_{mi}^k-d_i^k\perp Q_i^{k+1}\geq0,\quad i\in\mathcal{N}^D\backslash\{s\}$,\, $0\leq k\leq N_t-1$;\label{eq:micp4}\\
    &$V_{ij}^{k+1}-V_{ij}^k+\Delta t\left(u_{ij}^k\frac{\lambda_j^{k+1}-V_{ij}^{k+1}}{\Delta x}+f(u_{ij}^k,\rho_{ij}^k)\right)=0$,\quad $(i,j)\in\mathcal{L}^D$,\, $0\leq k\leq N_t-1$;\label{eq:micp5}\\
    &$u_{ij}^k=\argmin_{u_{\text{min}}\leq u\leq u_{\text{max}}} \left\{u\frac{\lambda_j^{k+1}-V_{ij}^{k+1}}{\Delta x}+f(u,\rho_{ij}^k)\right\}=0$,\quad $(i,j)\in\mathcal{L}^D$,\, $0\leq k\leq N_t-1$;\label{eq:micp6}\\
    &$\lambda_i^k=\sum_{k'=k}^{N_t-1}\mathbf{1}_{k'\leq k+Q_i^k/(M_i\Delta t)<k'+1}\cdot\left[\left(k'+1-\frac{Q^k_i}{M_i\Delta t}-k\right)\pi_i^{k'}+\left(\frac{Q^k_i}{M_i\Delta t}+k-k'\right)\pi_i^{k'+1}\right]$\notag\\
    &\quad\quad$+\mathbf{1}_{k+Q_i^k/(M_i\Delta t)\geq N_t}\cdot\pi_i^{N_t}+f_{\text{que}}\left(\min\{\frac{Q_i^k}{M_i},T-k\Delta t\}\right)$,\quad $i\in\mathcal{N}^D\backslash\{s\}$,\, $0\leq k\leq N_t$ ;\label{eq:micp7}\\
    &$0\leq \pi_i^k\perp \sum_{(i,j)\in\mathcal{L}^D}\beta_{ij}^k-1\geq0$,\quad $i\in\mathcal{N}^D\backslash\{s\}$,\, $0\leq k\leq N_t-1$;\label{eq:micp8}\\
    &$0\leq\beta_{ij}^k\perp V_{ij}^k-\pi_i^k\geq0$,\quad
    $(i,j)\in\mathcal{L}^D$,\, $0\leq k\leq N_t-1$.\label{eq:micp9}
\end{subnumcases}


The given conditions are: $\rho_{ij}^0=0$ for all $(i,j)\in\mathcal{L}^D$, $Q_i^0=0$ for all $i\in\mathcal{N}^D\backslash\{s\}$, and $\lambda_s^k=0$ for all $0\leq k\leq N_t$.
To solve the discretized system we also need to specify the following: (i) the traffic demand $d_i^k=d_i(t_k)$ for all $i\in\mathcal{N}^D$ and all time steps $t_k$ ($k=0,\dots,N_t)$;
(ii) the terminal cost $V_{ij}^{N_t}$ for all $(i,j)\in\mathcal{L}^D$ at the final time;
(iii) the terminal cost $\pi_i^{N_t}$ for all $i\in\mathcal{N}^D$ at the final time.

Now we are going to present a solution algorithm to solve 
the discretized [MFG-MiCP] system.
The algorithm employs fixed-point iterations on the queue sizes $\mathbf{Q}=(Q_i^k)_{\forall i,k}$, following the approach in \cite{ban2008link}.
We will start from any feasible initialization of queue sizes $\mathbf{Q}^{(0)}$. 
In each iteration, we substitute the current queue sizes $\mathbf{Q}^{(n)}$ into the summation term in Eq.~\eqref{eq:micp7}, which gives a relaxed MiCP. The solution of the relaxed MiCP gives a set of updated queue sizes $\mathbf{Q}^{(n+1)}$. We keep doing the iteration and updating the queue sizes until the iteration error:
\begin{align}
    E^{(n)}_{\mathbf{Q}}=\left|\mathbf{Q}^{(n+1)}-\mathbf{Q}^{(n)}\right|_2^2\leq\varepsilon,\label{eq:iter_error}
\end{align}
where $|\cdot|_2$ is the 2-norm and $\varepsilon$ is a predefined threshold.

In each iteration of the solution algorithm, the relaxed MiCP is solved by the PATH solver built in GAMS. The version of GAMS used in this paper is GAMS 24.9.  The [MFG-LWR-MiCP] system can be discretized and solved in the same way.

\section{Numerical Experiments}\label{sec:numerical}
In this section we will demonstrate our mean field game model and solution algorithm on Braess networks and the OW network. 
The basic experimental settings are illustrated in Section~\ref{sec:exp_setting}. 
Then four experiments are designed for different purposes. 
\begin{enumerate}
    \item In the first experiment, we will validate numerical convergence of our solution algorithm 
(Section~\ref{sec:convergence}). 
    \item In the second experiment, we will compare the equilibria solved from [MFG-MiCP] and [MFG-LWR-MiCP] with the same running and queuing cost functions on a three-path Braess network to demonstrate the role of driving speed control (Section~\ref{sec:comp_mfg_lwr}). 
    \item In the third experiment, we will demonstrate the occurrence of Braess paradox in the context of AVs by comparing their travel costs solved from [MFG-MiCP] on Braess networks with two and three paths (Section~\ref{sec:paradox}).
    
    \item In the fourth experiment, we will test the performance of the presented solution algorithm on the OW network, and compare the equilibria solved from [MFG-MiCP] and [MFG-LWR-MiCP] with the same running and 
    cost functions  (Section~\ref{sec:num_ow}).
\end{enumerate}

\subsection{Experimental settings}\label{sec:exp_setting}
\begin{figure}[htbp]
	\centering
	\begin{subfigure}{.48\textwidth}
	\centering
	\begin{tikzpicture}[scale=.4]
        \begin{scope}[every node/.style={circle,thick,draw}]
            \node (s) at (-6,0) {1};
            \node (d) at (6,0) {4};
            \node (1) at (0,3) {2};
            \node (2) at (0,-3) {3};
        \end{scope}
        \begin{scope}[>=Stealth,auto,
                      every edge/.style={draw=black,very thick}]
            \path [->] (s) edge node {} (1);
            \path [->] (s) edge node {} (2);
            \path [->] (1) edge node {} (d);
            \path [->] (2) edge node {} (d);
        \end{scope}
    \end{tikzpicture}
    \caption{Two-path}
    \label{fig:net1}
    \end{subfigure}
    \begin{subfigure}{.48\textwidth}
    \centering
    \begin{tikzpicture}[scale=.4]
        \begin{scope}[every node/.style={circle,thick,draw}]
            \node (s) at (-6,0) {1};
            \node (d) at (6,0) {4};
            \node (1) at (0,3) {2};
            \node (2) at (0,-3) {3};
        \end{scope}
        \begin{scope}[>=Stealth,auto,
                      every edge/.style={draw=black,very thick}]
            \path [->] (s) edge node {} (1);
            \path [->] (s) edge node {} (2);
            \path [->] (1) edge node {} (d);
            \path [->] (2) edge node {} (d);
            \path [->] (1) edge node {} (2);
        \end{scope}
    \end{tikzpicture}
    \caption{Three-path}
    \label{fig:net2}
    \end{subfigure}
    \caption{Braess networks}
    \label{fig:net}
\end{figure}
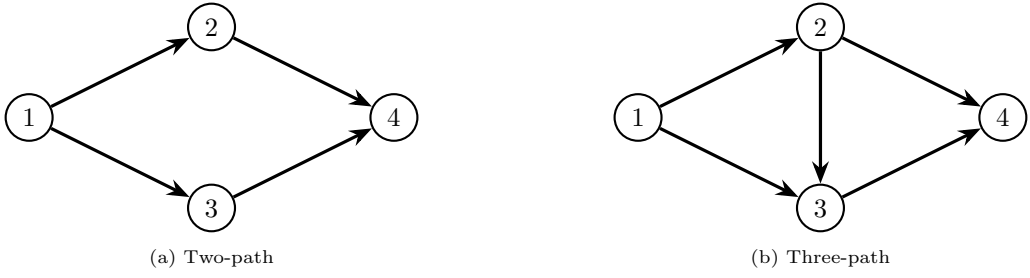

We will work on three networks: the two-path network with four links, as shown in Figure~\ref{fig:net1}, the three-path network with five links, as shown in Figure~\ref{fig:net2}, and the OW network that will be shown in Section~\ref{sec:num_ow}. The length of every link is assumed to be one if not particularly indicated, which is $\text{len}(l)=1,\,\forall l\in\mathcal{L}$. 

We define the AV's link running cost function as:
\begin{align}
	f_{\text{run}}(u,\rho)=\underbrace{\frac{c_1}2\left(\frac{u}{u_{\text{max}}}\right)^2}_{\text{kinetic energy}}+\underbrace{c_2\frac{\rho}{\rho_{\text{jam}}}}_{\text{safety}}+\underbrace{c_3}_{\text{efficiency}},\label{eq:costfct1}
\end{align}
where $c_1$, $c_2$ and $c_3$ are positive coefficients. 
The first term represents the kinetic energy; the second term quantifies driving safety by taking traffic densities as a penalty term, meaning that AVs tend to avoid congested areas; the third term quantifies driving efficiency, 
whose integration over one AV's travel time horizon $[t_0,t_f]$ on a network gives:
\begin{align}
    \int_{t_0}^{t_f}c_3\,dt=c_3(t_f-t_0),
\end{align}
which is proportional to the AV's total travel time on a network.

We define the AV's node queuing cost function as:
\begin{align}
    f_{\text{que}}(t_{\text{delay}})=c_4t_{\text{delay}},\label{eq:costfct2}
\end{align}
which is proportional to the AV's queuing delay $t_{\text{delay}}$ at a node.

The speed-density function for [MFG-LWR-MiCP] is set to be:
\begin{align}
    U(\rho)=u_{\text{max}}\left(1-\frac{\rho}{\rho_{\text{jam}}}\right).
\end{align}
We will assume $u_{\text{min}}=0$, $u_{\text{max}}=1$ and $\rho_{\text{jam}}=1$ if not particularly indicated.
To simulate a numerical experiment, we  still need to specify the following:
\begin{itemize}
    \item The simulation time $T$;
    \item The coefficients $c_1$, $c_2$, $c_3$ and $c_4$ in the cost functions defined in Eq.~\eqref{eq:costfct1} and Eq.~\eqref{eq:costfct2};
    \item The bottleneck capacities $M_i$ for all $i\in\mathcal{N}\backslash\{s\}$;
    \item the traffic demand $d_i(t)$ for all $i\in\mathcal{N}_O$ and $t\in[0,T]$;
    \item The terminal cost $V_l(x,T)$ for all $l\in\mathcal{L}$ and $\pi_i(T)$ for all $i\in\mathcal{N}\backslash\{s\}$ at the final time. 
\end{itemize}


\subsection{Numerical convergence}\label{sec:convergence}

To ensure the credibility of numerical solutions solved from the algorithm presented in Section~\ref{sec:algo}, it is necessary to validate numerical convergence. That is, when the discretization mesh sizes $\Delta x$ and $\Delta t$ go to zero, the numerical solution converges to the true solution of the continuum system.

In this subsection, we will validate numerical convergence with the [MFG-MiCP] system on the two-path network. The simulation time $T=3$ and the cost function coefficients $c_1=c_2=1$, $c_3=0.5$, $c_4=1$. The bottleneck capacity $M_i=1$ for every $i=1,2,3$. Cars enter the network through the single origin $o=1$ and the traffic demand is given by:
\begin{align}
    d_1(t)=\begin{cases}
          0.5,\quad &t\in[0,0.5];\\
          0,\quad &t\in(0.5,3].
    \end{cases}
    \label{equ:traffic_demand}
\end{align}
At the final time, the terminal costs $V_l(x,T)$ and $\pi_i(T)$ are set to be zero for all $l\in\mathcal{L}$ and $i\in\mathcal{N}\backslash\{s\}$.

\begin{figure}[H]
    \centering
    \begin{subfigure}{.48\textwidth}
        \includegraphics[width=\textwidth]{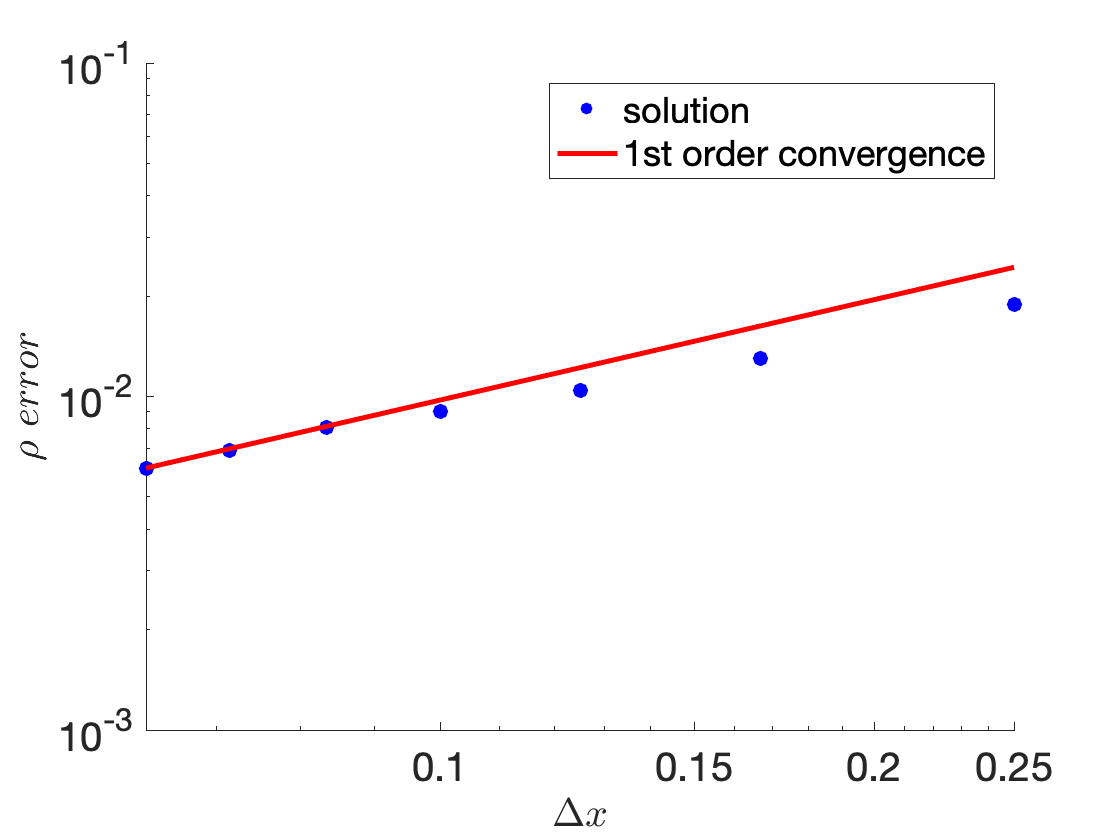}
        \caption{Error on $\rho$}
    \end{subfigure}
    \begin{subfigure}{.48\textwidth}
        \includegraphics[width=\textwidth]{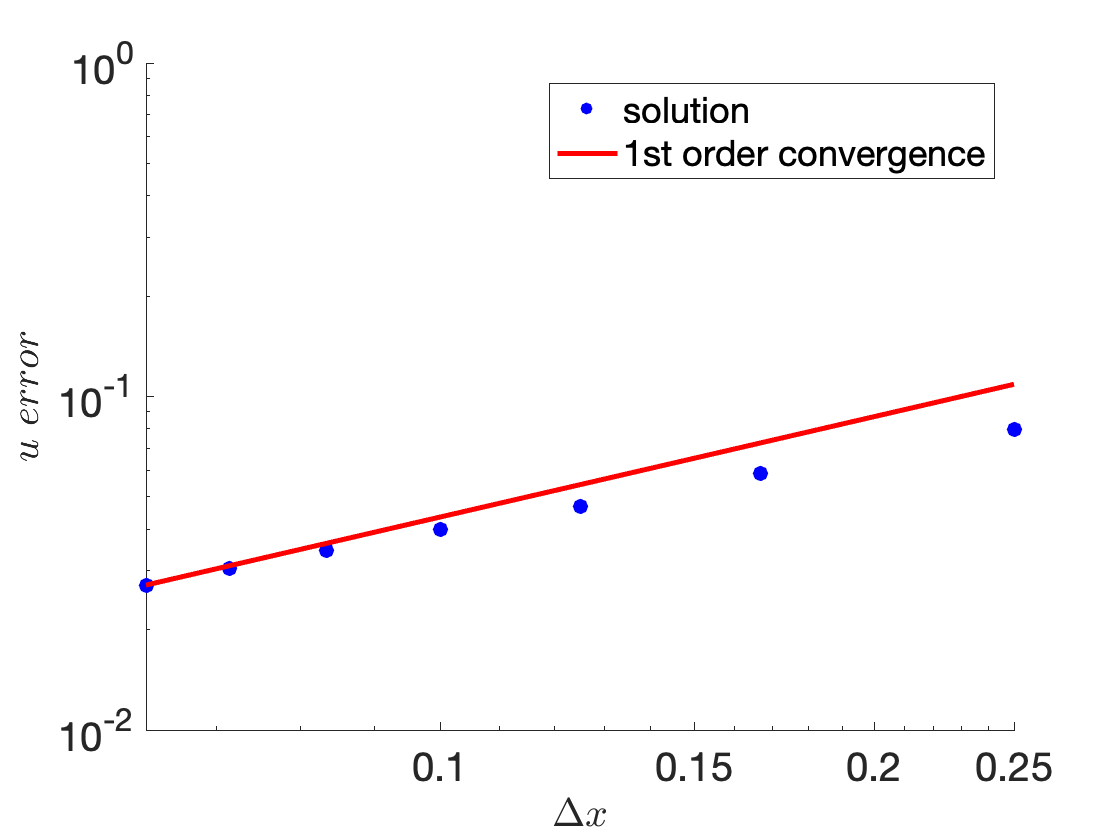}
        \caption{Error on $u$}
    \end{subfigure}
    
    \begin{subfigure}{.48\textwidth}
        \includegraphics[width=\textwidth]{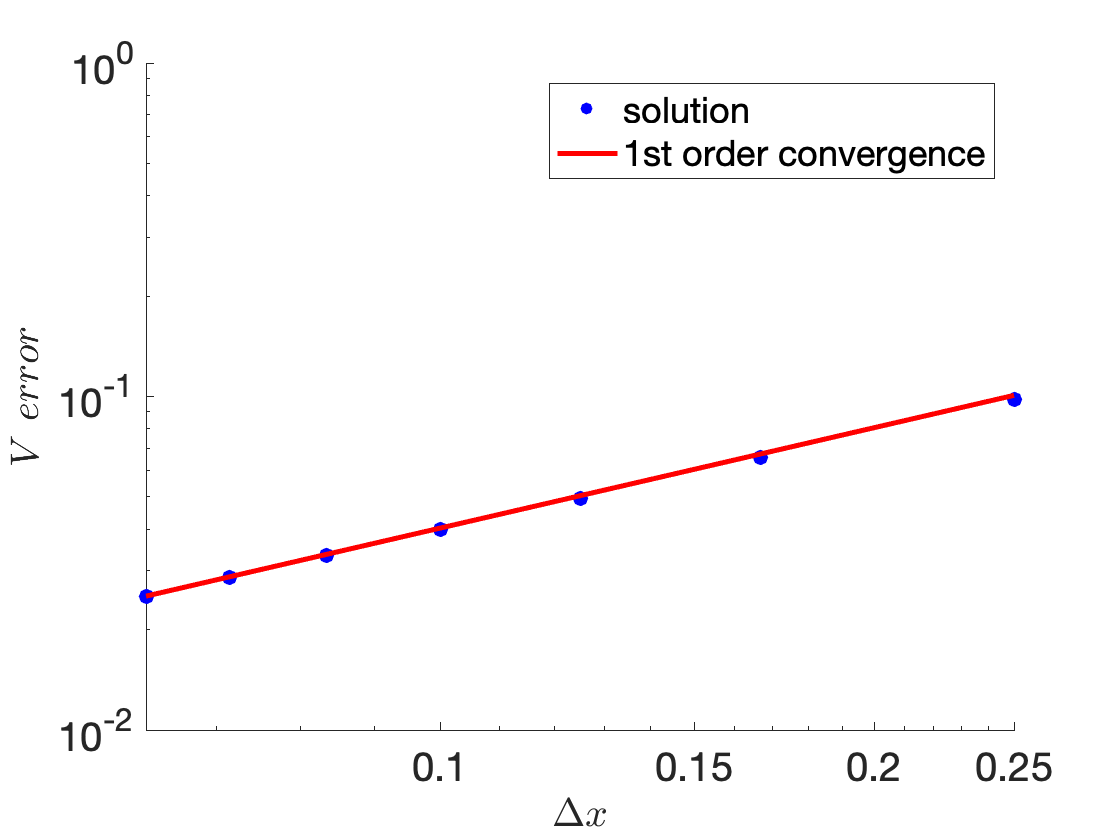}
        \caption{Error on $V$}
    \end{subfigure}
    \begin{subfigure}{.48\textwidth}
        \includegraphics[width=\textwidth]{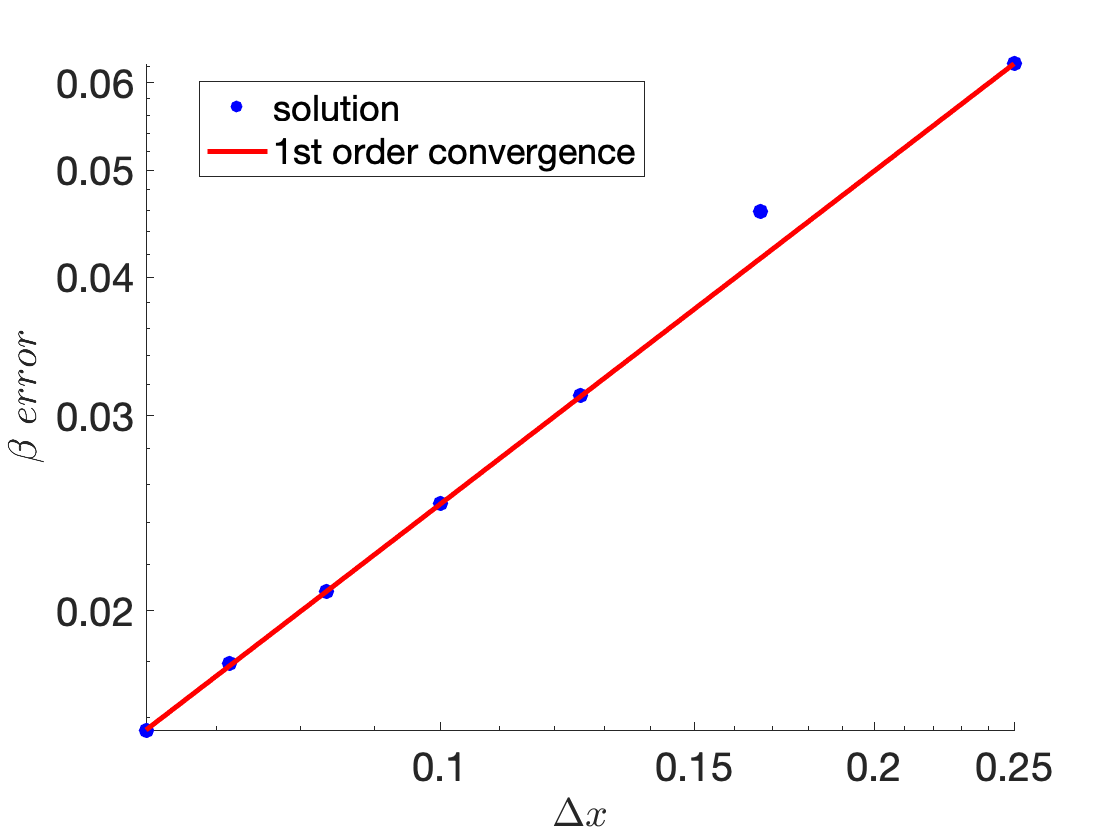}
        \caption{Error on $\beta$}
    \end{subfigure}
    \caption{Numerical convergence}
    \label{fig:convergence}
\end{figure}

Since the explicit solution of [MFG-MiCP] is not available, we will validate numerical convergence by comparing numerical solutions on multi-grid. We solve [MFG-MiCP] with different mesh sizes $\Delta x$ and $\Delta t$, but keep the ratio $\Delta x/\Delta t=1$. 
For each $\Delta x$, we compare the coarse solution with mesh size $\Delta x$ and the refined solution with mesh size $\Delta x/2$, which are denoted as SOL([MFG-MiCP], $\Delta x$) and SOL([MFG-MiCP], $\Delta x/2$), respectively. Note that these two numerical solutions are defined on different mesh grids, we will first project the coarse solution SOL([MFG-MiCP], $\Delta x$) onto fine grids with mesh size $\Delta x/2$ using the piecewise constant interpolation. 
Denote the projected solution as PROJ(SOL([MFG-MiCP], $\Delta x$), $\Delta x/2$). Then the solutions SOL([MFG-MiCP], $\Delta x/2$) and PROJ(SOL([MFG-MiCP], $\Delta x$), $\Delta x/2$) are on the same mesh grids. We compute the mean absolute error between the two solutions on variables $\rho$, $u$, $V$ and $\beta$. The errors are plotted with respect to $\Delta x$ in Figure~\ref{fig:convergence}.


In Figure~\ref{fig:convergence}, the $x$-axis is the mesh size $\Delta x$ and the $y$-axis is the numerical error. The blue dots represent errors computed from different values of $\Delta x$. The red line represents the perfect first order convergence, i.e., the error is proportional to $\Delta x$. 
We observe that the errors on all variables go to zero as $\Delta x\to0$. This means the change on the numerical solution will go to zero if we keep refining the mesh grids, which validates the numerical convergence \citep{leveque2007finite}. In particular, it is shown in the figure that the numerical solution has first order convergence.

\subsection{Comparison between MFG and LWR equilibria on the three-path network}\label{sec:comp_mfg_lwr}

In this subsection, we would like to compare the MFG and LWR-based optimal controls of AVs across the three-path network, and demonstrate the role of driving speed control. We solve both [MFG-MiCP] and [MFG-LWR-MiCP] with the same running and queuing cost functions on the three-path network. The solved equilibria are compared in terms of density evolution, velocity profiles and queuing delays.


The settings of these experiments are as follows:

\begin{itemize}
    \item The simulation time $T=12$. 
    \item The links $(1,2)$ and $(1,3)$ have the same running cost function coefficients $c_1=c_2=c_3=1$; the link $(2,4)$ has the running cost function coefficients $c_1=5$, $c_2=5$ and $c_3=5$; the link $(3,4)$ has the running cost function coefficients $c_1=0.5$, $c_2=0$ and $c_3=0.25$; the link $(2,3)$ has the running cost function coefficients $c_1=0.5$, $c_2=0$ and $c_3=0.1$. With these coefficients, the links $(1,2)$ and $(1,3)$ are symmetric.
    \item The bottleneck capacity $M_1=M_2=1$ and $M_3=0.35$. The queuing cost function coefficient $c_4=1$.
    \item AVs enter the network through origin nodes $1$ and $2$. The traffic demands are given by:
    \begin{align}
    d_1(t)=\begin{cases}
          0.6,\quad &t\leq0.5;\\
          0,\quad &t>0.5;
    \end{cases}, \quad
    d_2(t)=\begin{cases}
          0.5,\quad &t\leq0.5;\\
          0,\quad &t>0.5.
    \end{cases}
    \end{align}
    \item The AVs' terminal costs at time $T$ are given in the following way: we first specify the nodal terminal costs $\pi_1(T)=2$, $\pi_2(T)=1.33$, $\pi_3(T)=0.67$; the link terminal cost $V_l(x,T)$ is given by linear interpolation between the nodal terminal costs.
\end{itemize}

 

For both the solved MFG and LWR equilibria, the link $(2,4)$ turns out not to be selected, likely because its running cost coefficients are large. Hence cars entering the network through node $1$ drive along the path $1\to3\to4$ while cars entering the network through node $2$ drive along the path $2\to3\to4$. The flows from the two paths merge at node $3$ and form a non-empty queue.

\begin{figure}[htbp]
    \centering
    \begin{subfigure}{.35\textwidth}
        \includegraphics[width=\textwidth]{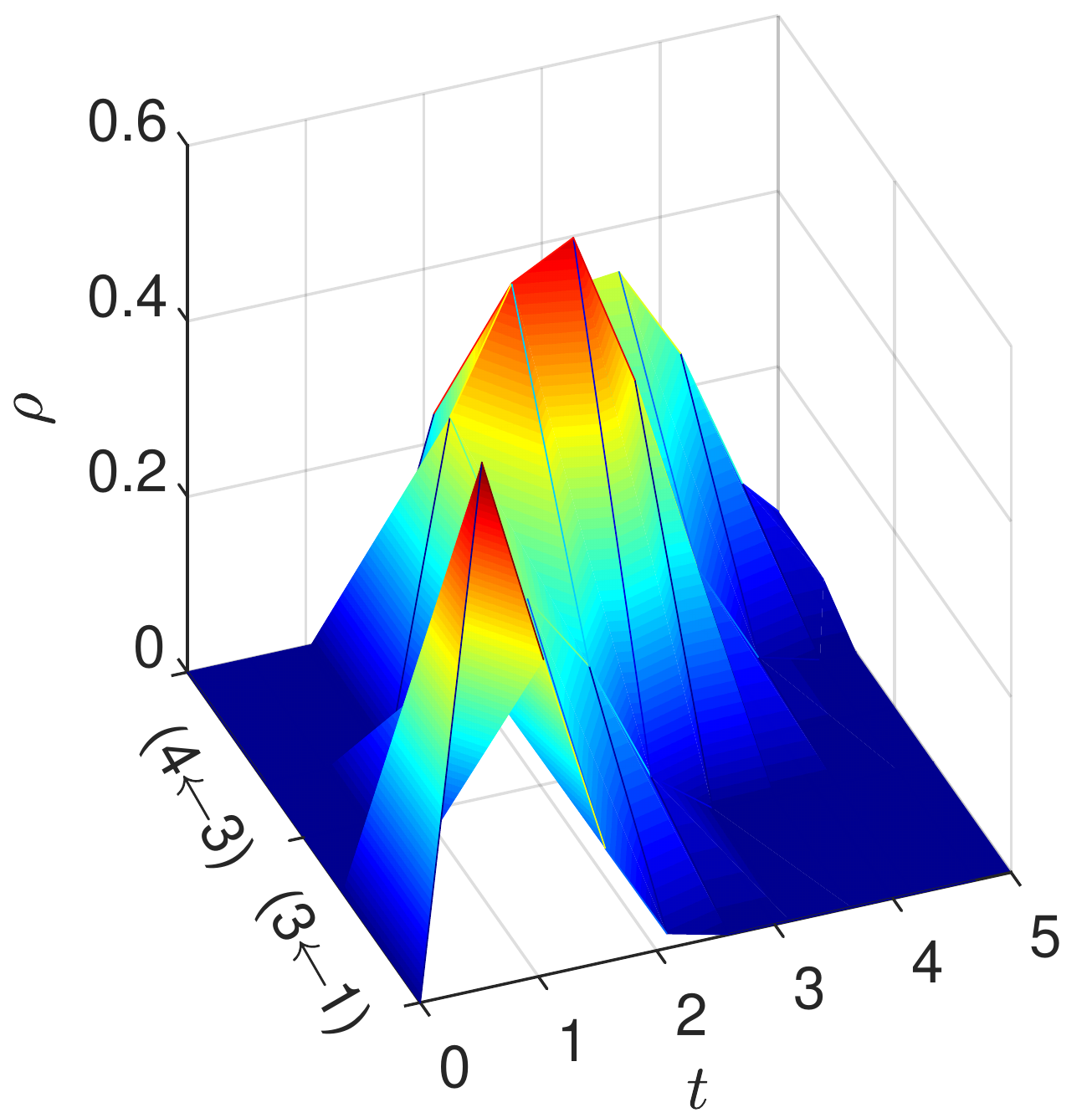}
        \caption{LWR Density: Path $1\rightarrow 3\rightarrow 4$}
        \label{fig:den_lwr_1}
    \end{subfigure}
    \hspace{0.6in}
    \begin{subfigure}{.35\textwidth}
        \includegraphics[width=\textwidth]{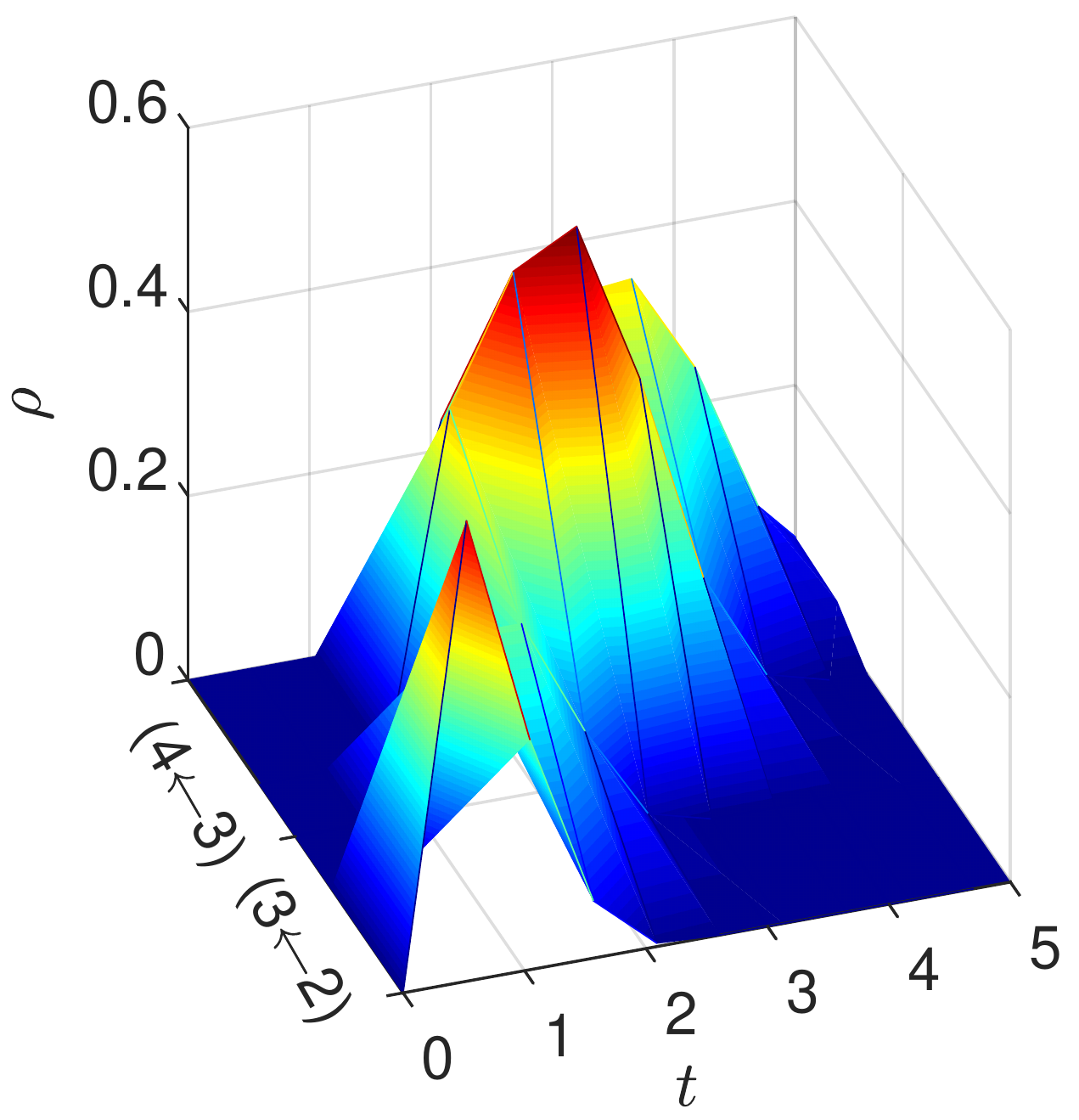}
        \caption{LWR Density: Path $2\rightarrow 3\rightarrow 4$}
        \label{fig:den_lwr_2}
    \end{subfigure}
    
    \begin{subfigure}{.35\textwidth}
        \includegraphics[width=\textwidth]{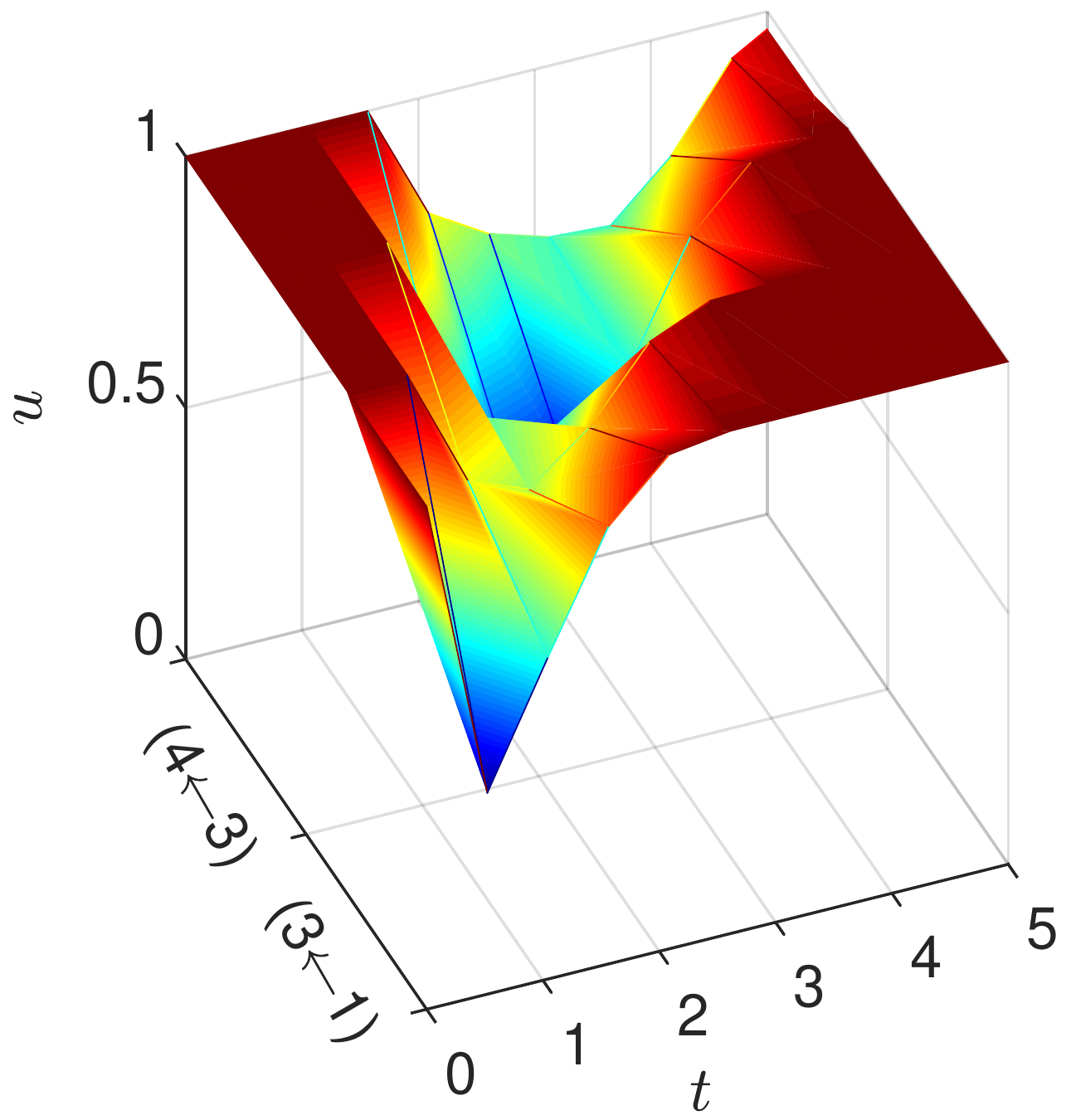}
        \caption{LWR Velocity: Path $1\rightarrow 3\rightarrow 4$}
        \label{fig:v_lwr}
    \end{subfigure}
    \hspace{0.6in}
    \begin{subfigure}{.35\textwidth}
        \includegraphics[width=\textwidth]{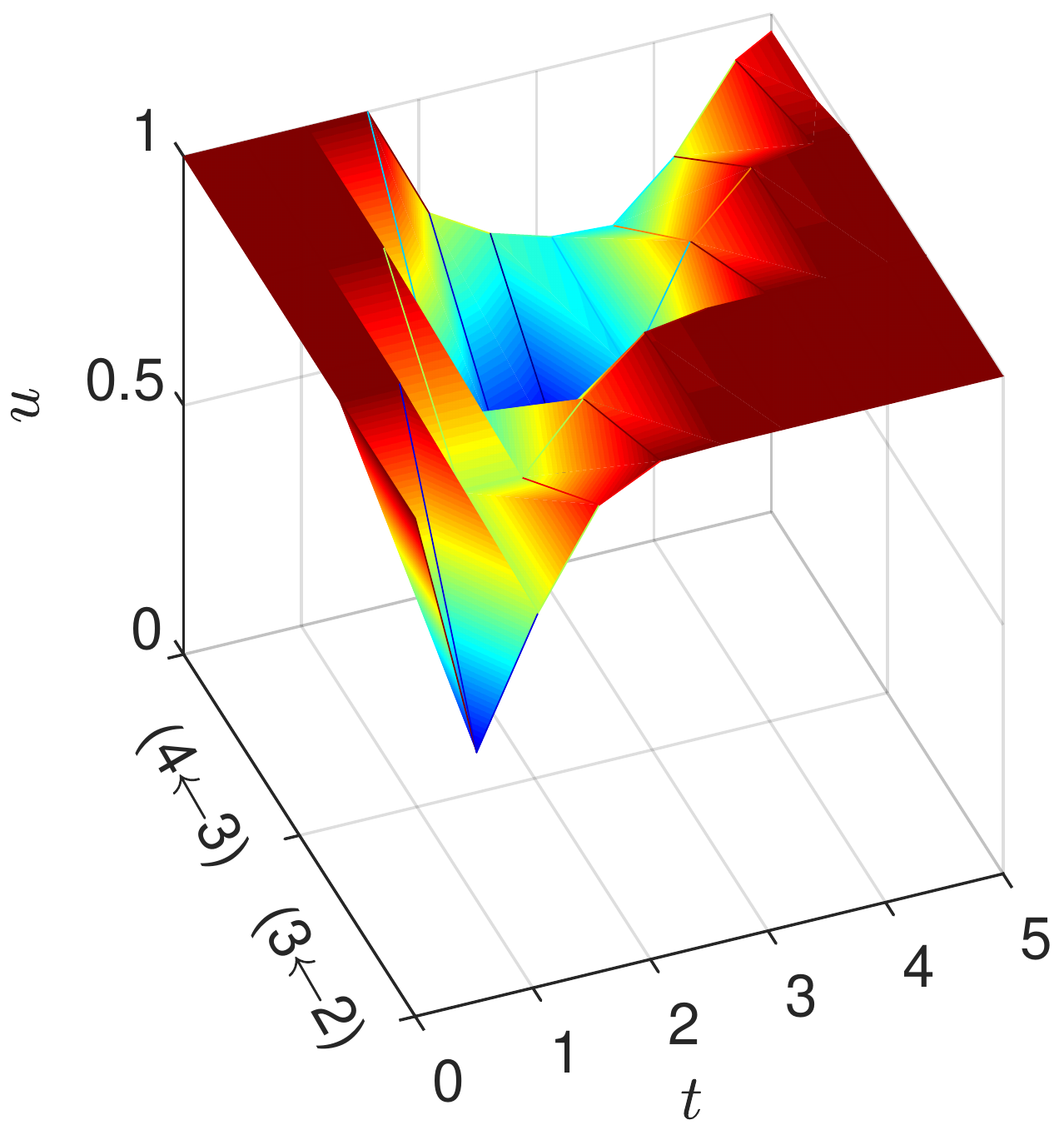}
        \caption{LWR Velocity: Path $2\rightarrow 3\rightarrow 4$}
        \label{fig:density_evolution_4_link_lwr}
    \end{subfigure}
    \caption{LWR equilibrium traffic on the three-path network}
    \label{fig:equ_lwr}
\end{figure}

\begin{figure}[htbp]
    \centering
    \begin{subfigure}{.34\textwidth}
        \includegraphics[width=\textwidth]{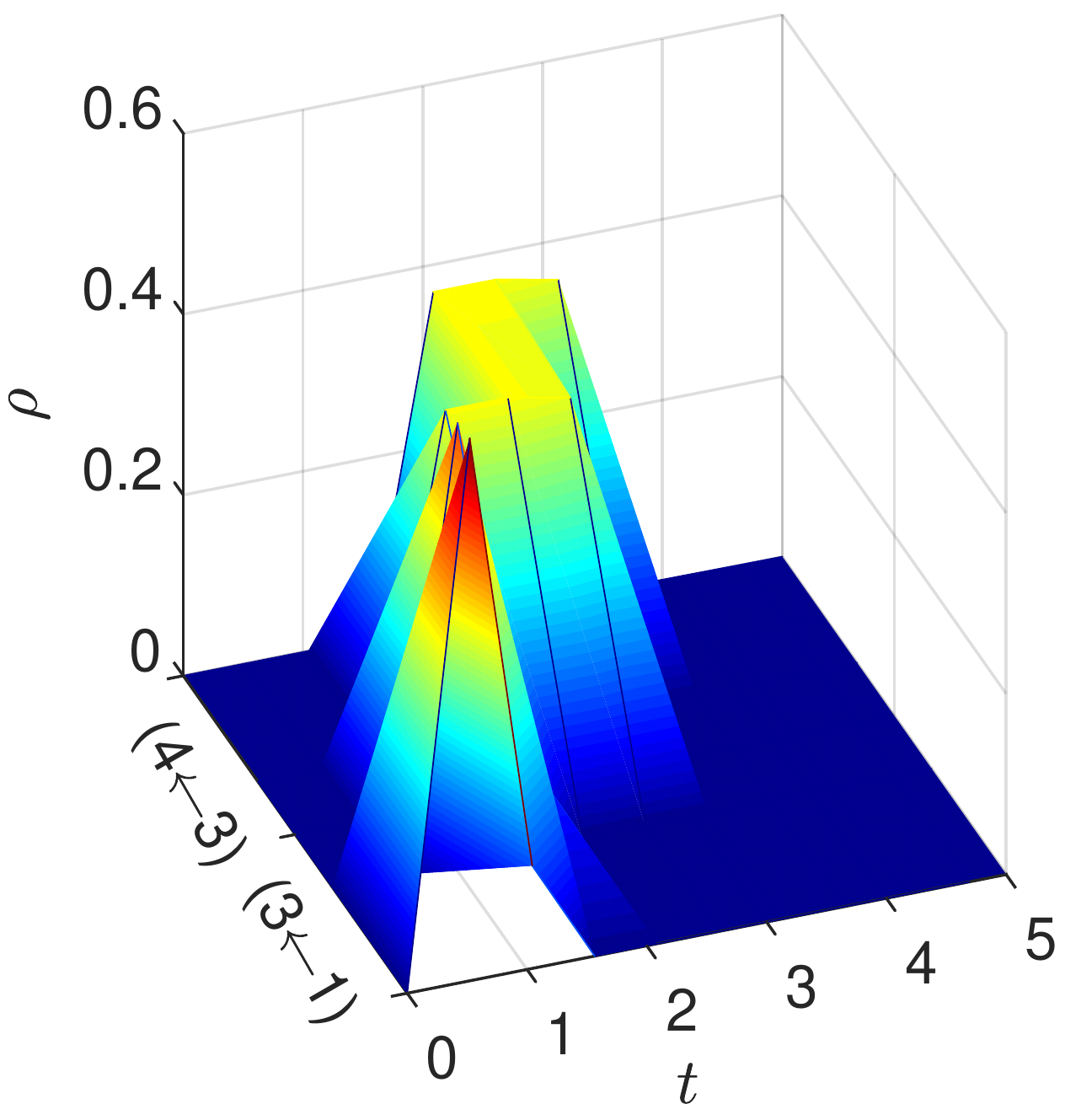}
        \caption{MFG Density: Path $1\rightarrow 3\rightarrow 4$}
        \label{fig:den_mfg_1}
    \end{subfigure}
    \hspace{0.6in}
    \begin{subfigure}{.35\textwidth}
        \includegraphics[width=\textwidth]{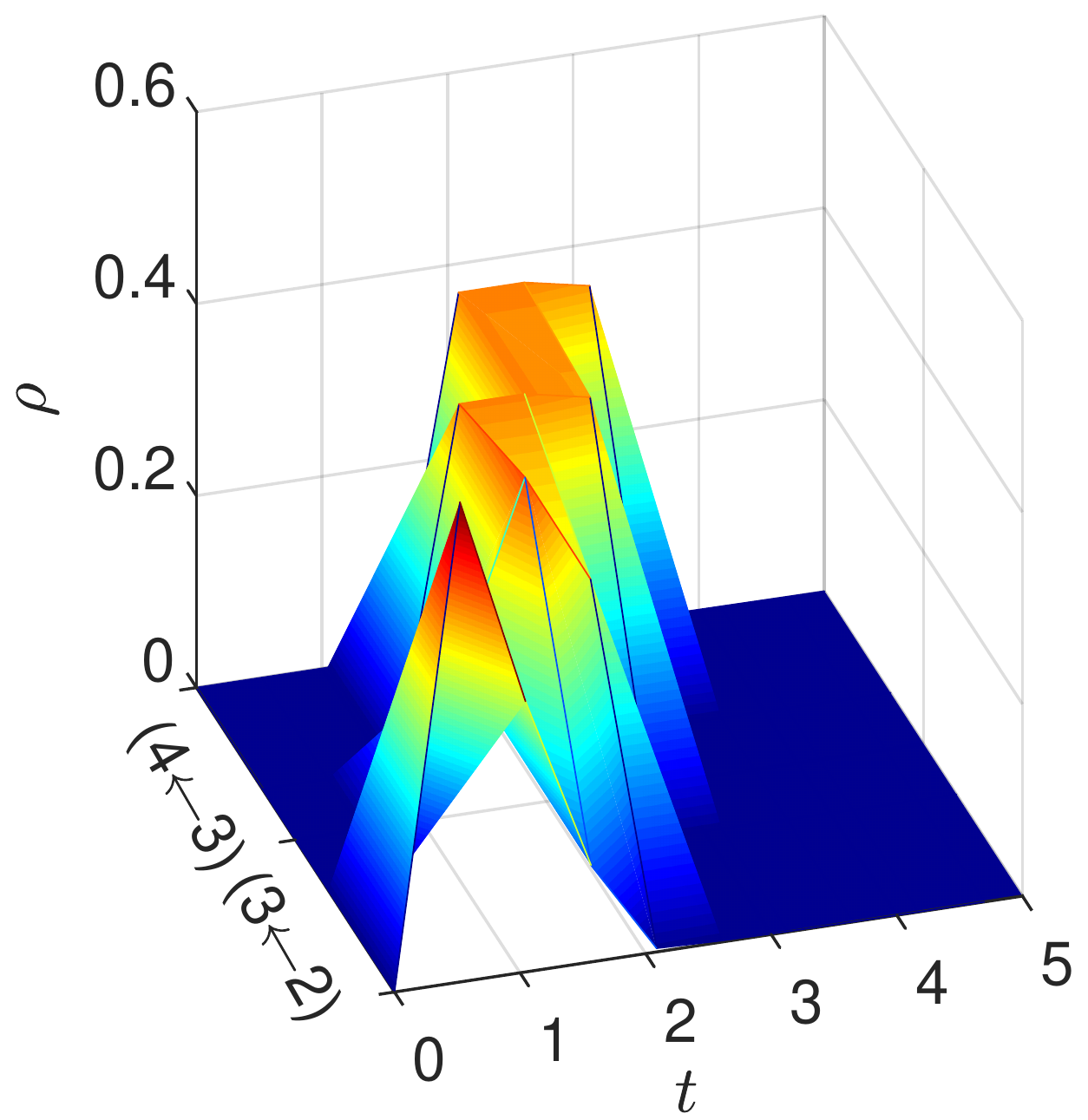}
        \caption{MFG Density: Path $2\rightarrow 3\rightarrow 4$}
        \label{fig:den_mfg_2}
    \end{subfigure}
    
    \begin{subfigure}{.35\textwidth}
        \includegraphics[width=\textwidth]{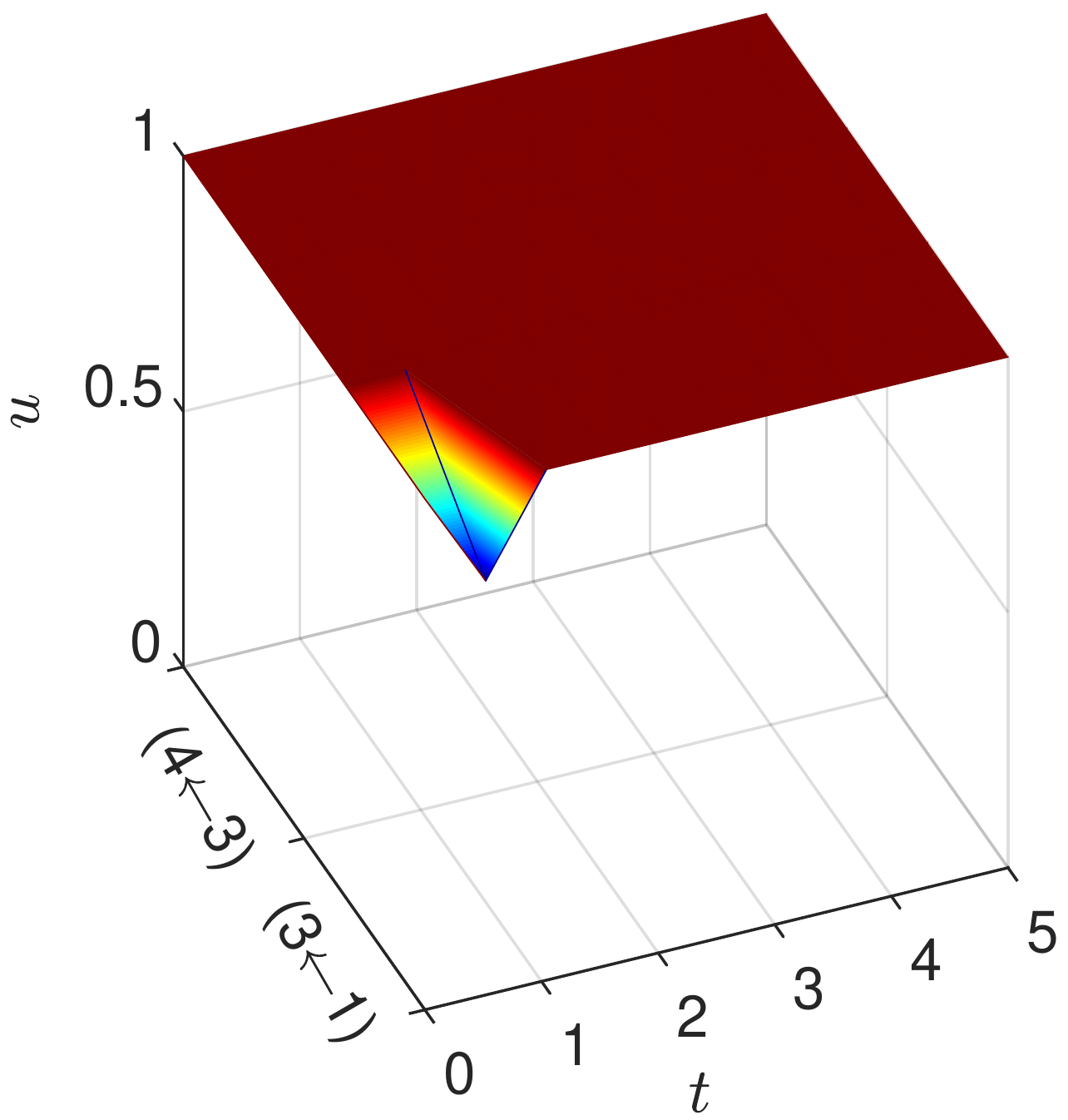}
        \caption{MFG Velocity: Path $1\rightarrow 3\rightarrow 4$}
        \label{fig:v_mfg}
    \end{subfigure}
    \hspace{0.6in}
    \begin{subfigure}{.35\textwidth}
        \includegraphics[width=\textwidth]{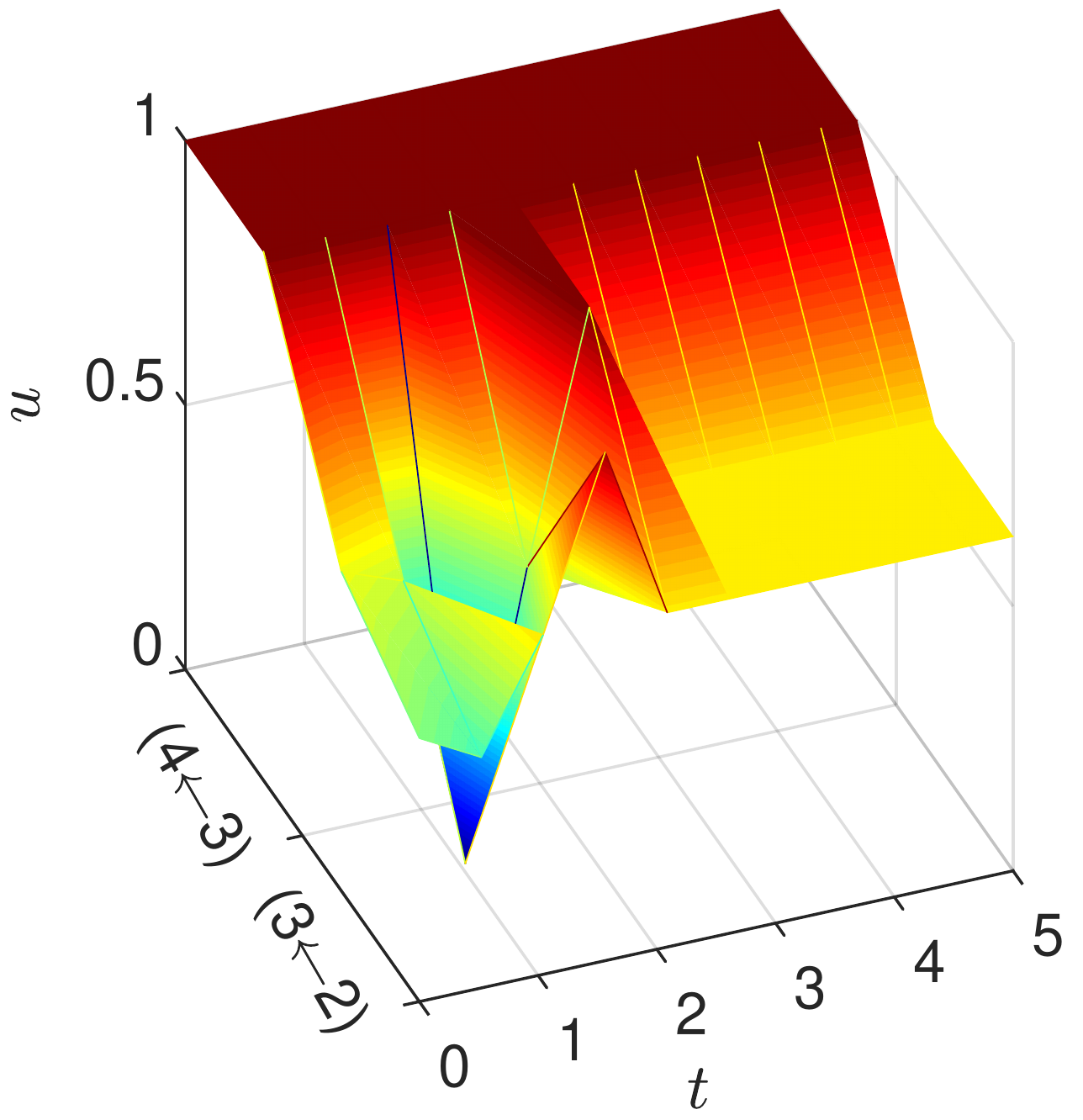}
        \caption{MFG Velocity: Path $2\rightarrow 3\rightarrow 4$}
        \label{fig:density_evolution_4_link_mfg}
    \end{subfigure}
    \caption{MFG equilibrium traffic on the three-path network}
    \label{fig:equ_mfg}
\end{figure}

In Figure~\ref{fig:equ_lwr} the traffic density and velocity for the LWR equilibrium along the two paths is plotted in a 3D diagram. The $x$-axis represents the path as a continuous road of length 2,
the $y$-axis is time and the $z$-axis represents traffic density or velocity. In Figure~\ref{fig:equ_mfg} the traffic density and velocity for the MFG equilibrium is plotted in the same way.

By comparing Figure~\ref{fig:equ_mfg} to Figure~\ref{fig:equ_lwr},
we observe the difference between MFG and LWR equilibria on driving speeds: In the LWR equilibrium, AVs drive faster in low density areas and slower in high density areas. In other words, the LWR speed is totally determined by its local traffic density and thus this strategy is ``myopic". 
By contrast, in the MFG equilibrium, AVs may drive at relatively high speeds even if the local density is relatively high. That is, AVs deployed with the MFG equilibrium controls can look ``farther'' and adjust their driving speeds in a way to optimize their total travel costs over the network. In other words, even the cost incurred by high speed and high density is relatively large for the current link where an AV is, the AV can anticipate the future cost incurred on the subsequent links with the goal of minimizing the cumulative total cost.

\begin{figure}[htbp]
    \centering
    \begin{subfigure}{.39\textwidth}
        \includegraphics[width=\textwidth]{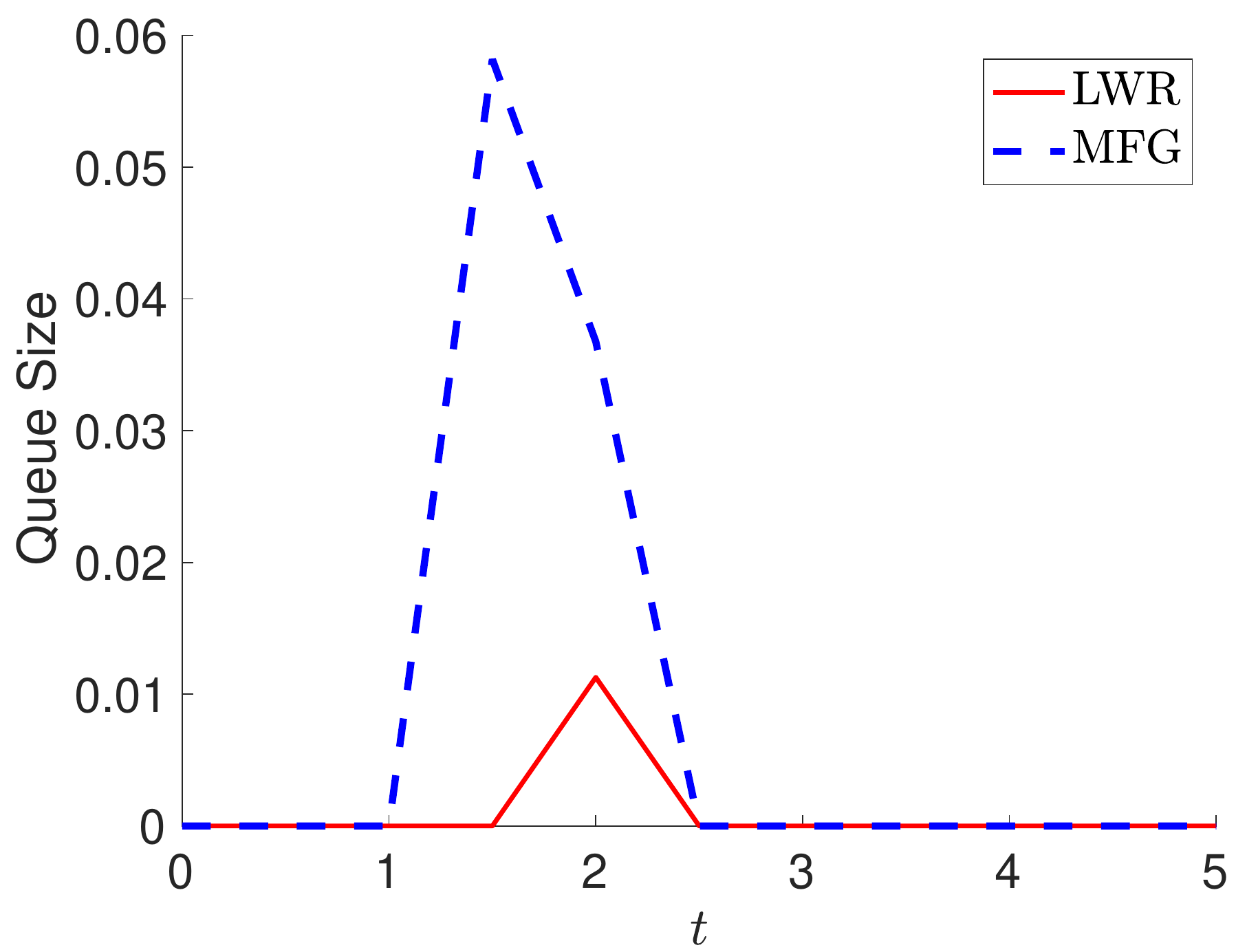}
        \caption{Queue Size at Node $3$}
        \label{fig:que_change}
    \end{subfigure}
    \hspace{0.6in}
    \begin{subfigure}{.35\textwidth}
        \includegraphics[width=\textwidth]{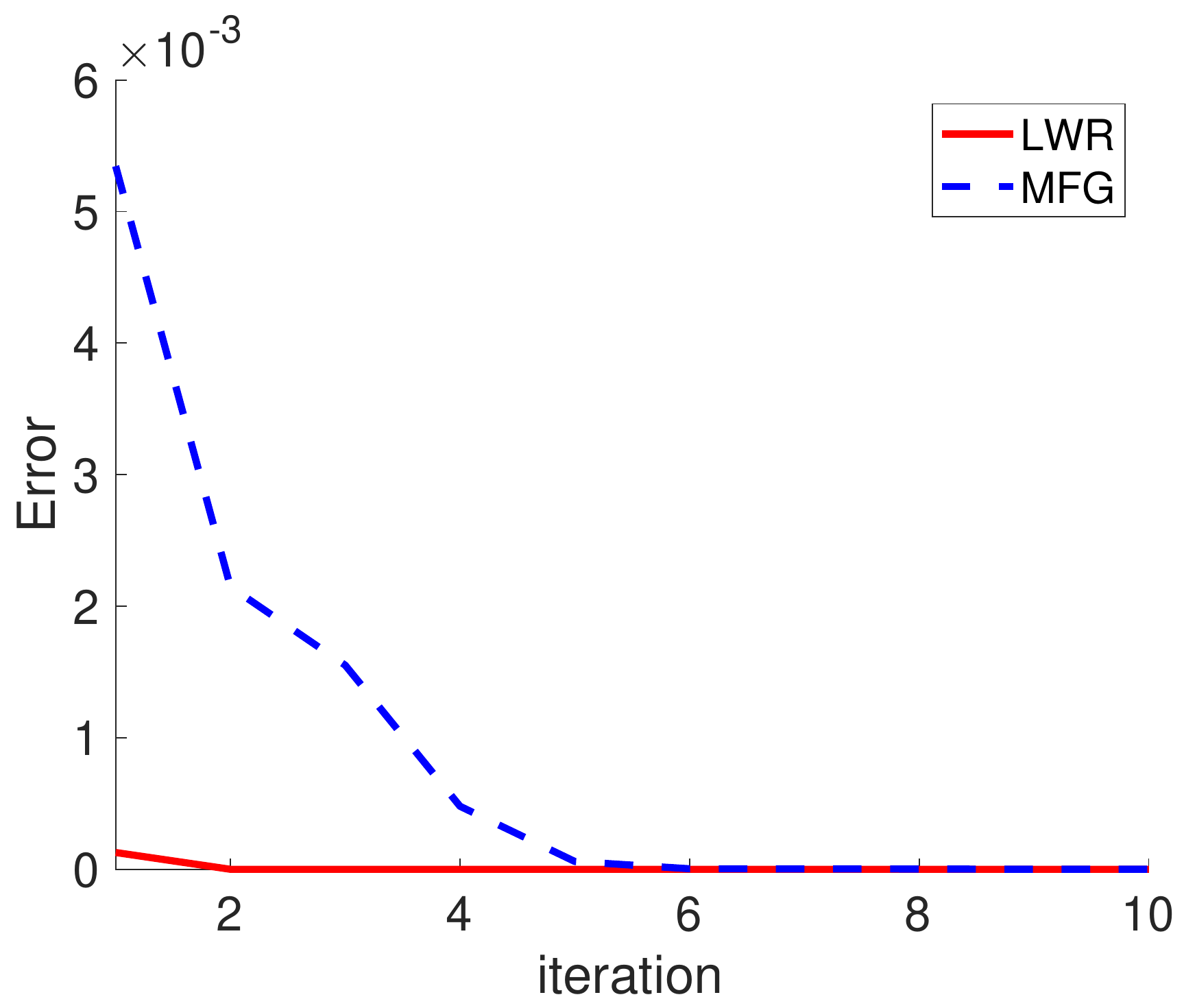}
        \caption{Iteration Error}
        \label{fig:queue_gap}
    \end{subfigure}
    \caption{Queuing information at node $3$ on the three-path network}
    \label{fig:compare}
\end{figure}

It is also interesting to compare the queue size at node $3$ of MFG and LWR equilibria during the time horizon $t\in[0,5]$, which is shown in Figure~\ref{fig:que_change}. We do not plot the queue size for later times because the network is already empty at $t=5$ for both equilibria. We observe from Figure~\ref{fig:que_change} that the queue size of the MFG equilibrium is larger than that of the LWR equilibrium. This is because AVs drive at higher speeds on both links $(1,3)$ and $(2,3)$ when $t\leq2$ in the MFG equilibrium than in the LWR equilibrium. Hence the node $3$ is more congested and more cars stuck in the queue at the node in the MFG equilibrium.
However, although AVs have longer queuing delays at node $3$, we observe from Figure~\ref{fig:den_mfg_1} and Figure~\ref{fig:den_mfg_2} that all AVs arrive at the destination at approximately $t=4$ in the MFG equilibrium, while that arrival time for AVs in the LWR equilibrium is approximately $t=5$ from Figure~\ref{fig:den_lwr_1} and Figure~\ref{fig:den_lwr_2}.

In Figure~\ref{fig:queue_gap}, the solution algorithm presented in Section~\ref{sec:algo} is validated by showing convergence of the fixed-point iteration on the queue size. The iteration error, which is defined in Eq.~\eqref{eq:iter_error}, converges to zero in two iterations for the LWR equilibrium and in six iterations for the MFG equilibrium. It verifies the efficiency of the solution algorithm.

\color{black}

\subsection{Braess paradox}\label{sec:paradox}

This subsection aims to demonstrate the occurrence of the Braess paradox. We will first solve a MFE on the two-path network. Then the middle link is added and a new MFE is solved on the three-path network. 
We artificially create the occurrence of the paradox,  which is when travel cost with the middle link is higher than that without.

The settings of this experiment are as follows: 
\begin{itemize}
    \item The simulation time $T=6$.
    \item 
    The links $(1,2)$ and $(3,4)$ have the same running cost function coefficients $c_1=1$, $c_2=5$ and $c_3=0$; The links $(1,3)$ and $(2,4)$ have the same running cost function coefficients $c_1=1$, $c_2=0$ and $c_3=3$. It means that $(1,2)$ and $(3,4)$ are bottleneck links on which AVs' travel costs highly depend on traffic density. In contrast, AVs can move freely on links $(1,3)$ and $(2,4)$ to minimize their travel times no matter how congested the links are.
    \item The three-path network is constructed by adding the middle link $(2,3)$ into the two-path network. The middle link $(2,3)$ has a length of $0.25$ and its running cost function coefficients are $c_1=1$, $c_2=5$ and $c_3=0$. Since this link is very short, AVs can traverse it in a short time and incur a smaller travel cost than other links.
    \item The bottleneck capacity $M_i=0.8$ for every $i=1,2,3$ on both two-path and three-path networks. The queuing cost function coefficient $c_4=1$.
    \item For both cases, AVs enter the network through the single origin $o=1$ and the traffic demand is given by:
    \begin{align}
        d_1(t)=\begin{cases}
          0.75,\quad &t\in[0,1];\\
          0,\quad &t\in(1,6].
    \end{cases}
    \end{align}
    \item The AVs' terminal costs at time $T$ are given by:
    \begin{align}
        V_l(x,T)=\begin{cases}
            2-x,\quad &l=(1,2),\,(1,3),\,(2,3);\\
            1-x,\quad &l=(2,4),\,(3,4),
        \end{cases}
    \end{align}
    and:
    \begin{align}
        \pi_i(T)=\begin{cases}
            2,\quad &i=1;\\
            1,\quad &i=2,3.
        \end{cases}
    \end{align}
    It means that there will be a penalty if an AV does not reach the destination before time $T$, and the penalty cost depends on the distance between the AV's position on the network and the destination.
\end{itemize}

\begin{figure}[H]
    \centering
    \begin{subfigure}{.48\textwidth}
        \includegraphics[width=\textwidth]{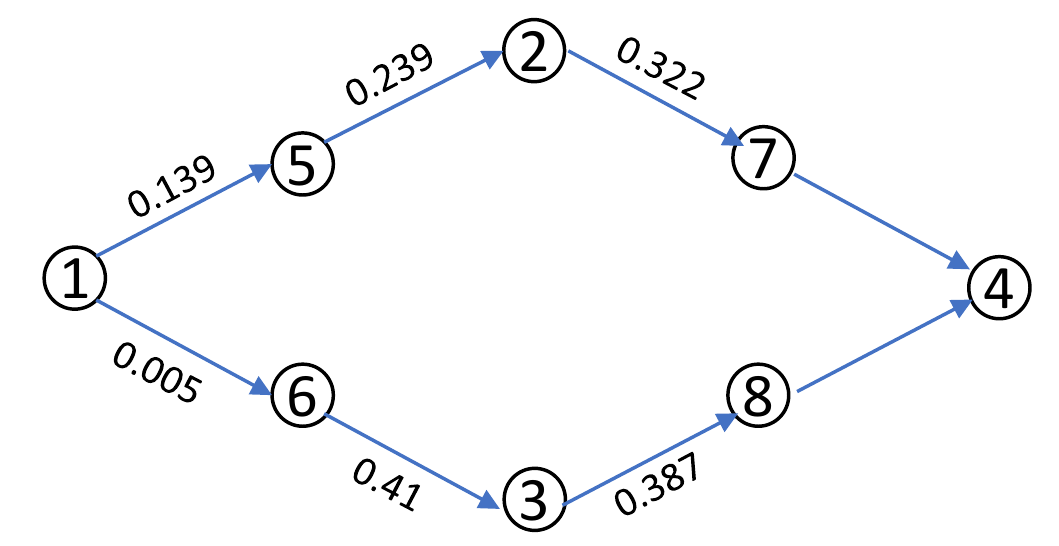}
        \caption{Two-path, $t=1.75$}\label{fig:paradox_flow_1}
    \end{subfigure}
    \begin{subfigure}{.48\textwidth}
        \includegraphics[width=\textwidth]{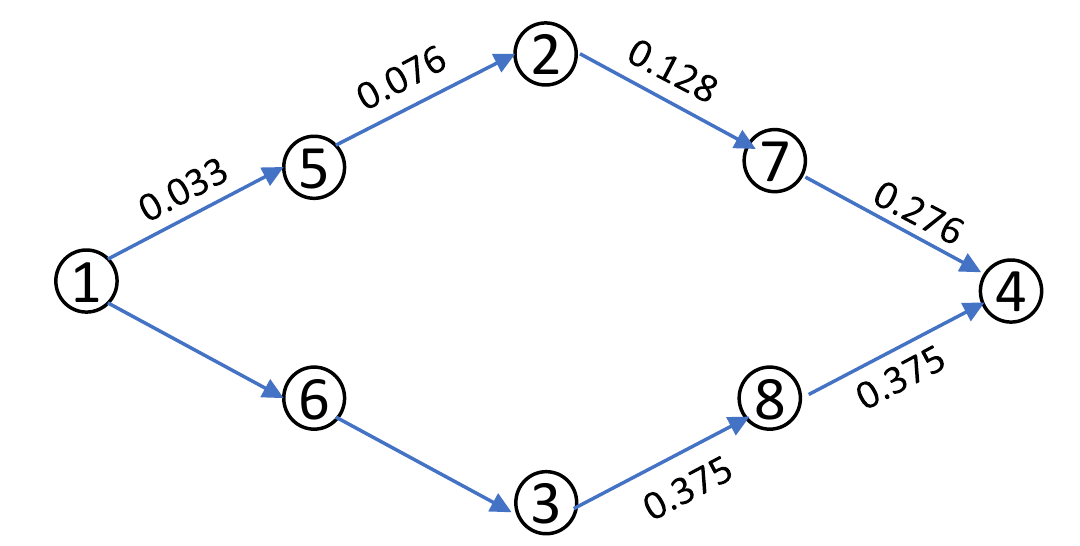}
        \caption{Two-path, $t=2.5$}\label{fig:paradox_flow_2}
    \end{subfigure}
    
    \begin{subfigure}{.48\textwidth}
        \includegraphics[width=\textwidth]{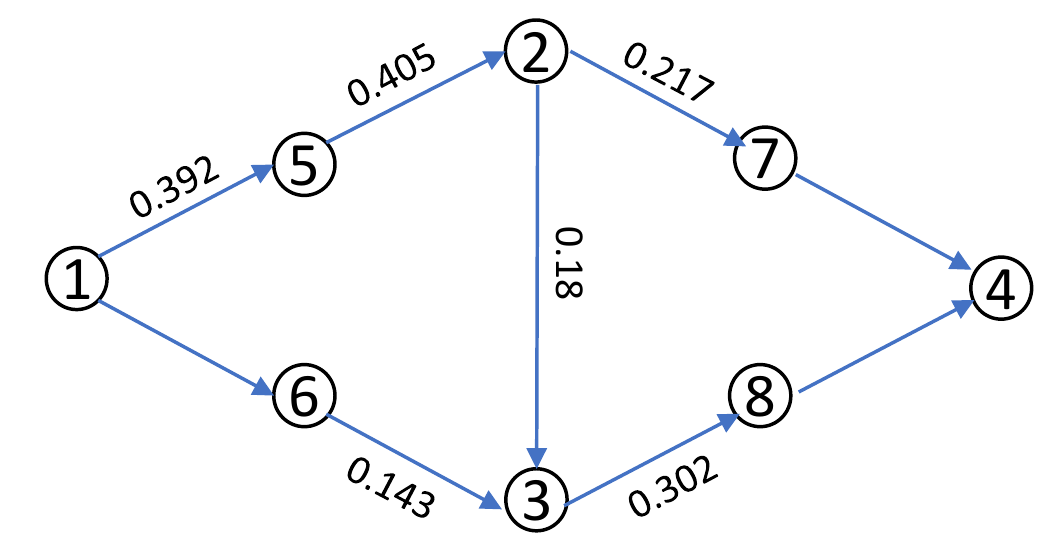}
        \caption{Three-path, $t=1.75$}\label{fig:paradox_flow_3}
    \end{subfigure}
    \begin{subfigure}{.48\textwidth}
        \includegraphics[width=\textwidth]{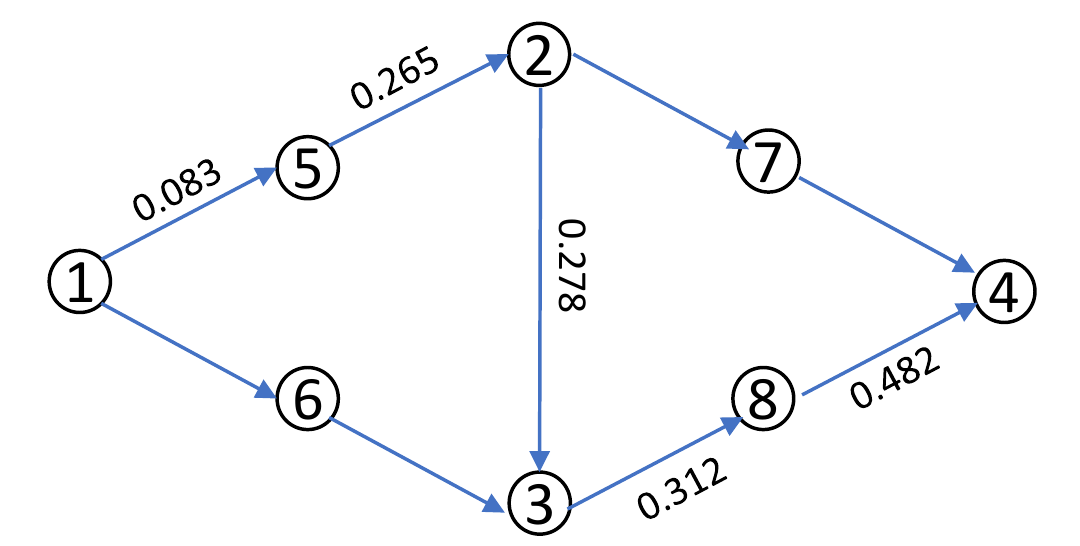}
        \caption{Three-path, $t=2.5$}\label{fig:paradox_flow_4}
    \end{subfigure}
    \caption{Density evolution on the two-path and three-path networks}
    \label{fig:paradox_flow}
\end{figure}

In Figure~\ref{fig:paradox_flow}, the traffic density is plotted on all links of the two networks at the time snapshots $t=1.75$ and $t=2.5$, respectively, to demonstrate the AVs' dynamic route choices on the two networks. 
We observe from Figure~\ref{fig:paradox_flow_1} and Figure~\ref{fig:paradox_flow_2} that both paths are used in the two-path network, but more AVs choose the top path $1\to2\to4$. This is different from the classical static Braess paradox in which cars are split equally on the two paths. In the latter situation, cars have the same travel time on links $(1,2)$ and $(3,4)$ because the two links have the same link travel time function and the total flow passing the two links is the same. 
Same for links $(1,3)$ and $(2,4)$.
So the total travel time is the same on the two paths. However, equally split flow is no longer an equilibrium in the dynamic case. Even if the links $(1,2)$ and $(3,4)$ have the same cost function coefficients and the total flow passing the two links is the same, the travel costs on the two links are in general different because the costs also depend on the distribution of the total flow with respect to time. The flow distribution entering the link $(3,4)$ at different times depends on the AVs' speed controls on the link $(1,3)$, which is in general different from the flow distribution entering the link $(1,2)$. For the same reason, AVs have different travel costs on links $(1,3)$ and $(2,4)$, and the two paths 
give different travel costs. 
To summarize, the existence of AVs' speed controls break the symmetry between the two paths.


We observe from Figure~\ref{fig:paradox_flow_3} and Figure~\ref{fig:paradox_flow_4} that a majority of AVs choose to travel along the path $1\to2\to3\to4$ on the three-path network. At $t=1.75$, only 22 percent of AVs are on link $(1,3)$ and $(2,4)$; at $t=2.5$, all AVs are on the path $1\to2\to3\to4$. This is because AVs experience a small travel cost on the middle link $(2,3)$, so most of them choose to go along the path $1\to2\to3\to4$ and benefit from the middle link.

\begin{figure}[H]
    \centering
    \includegraphics[width=.5\textwidth]{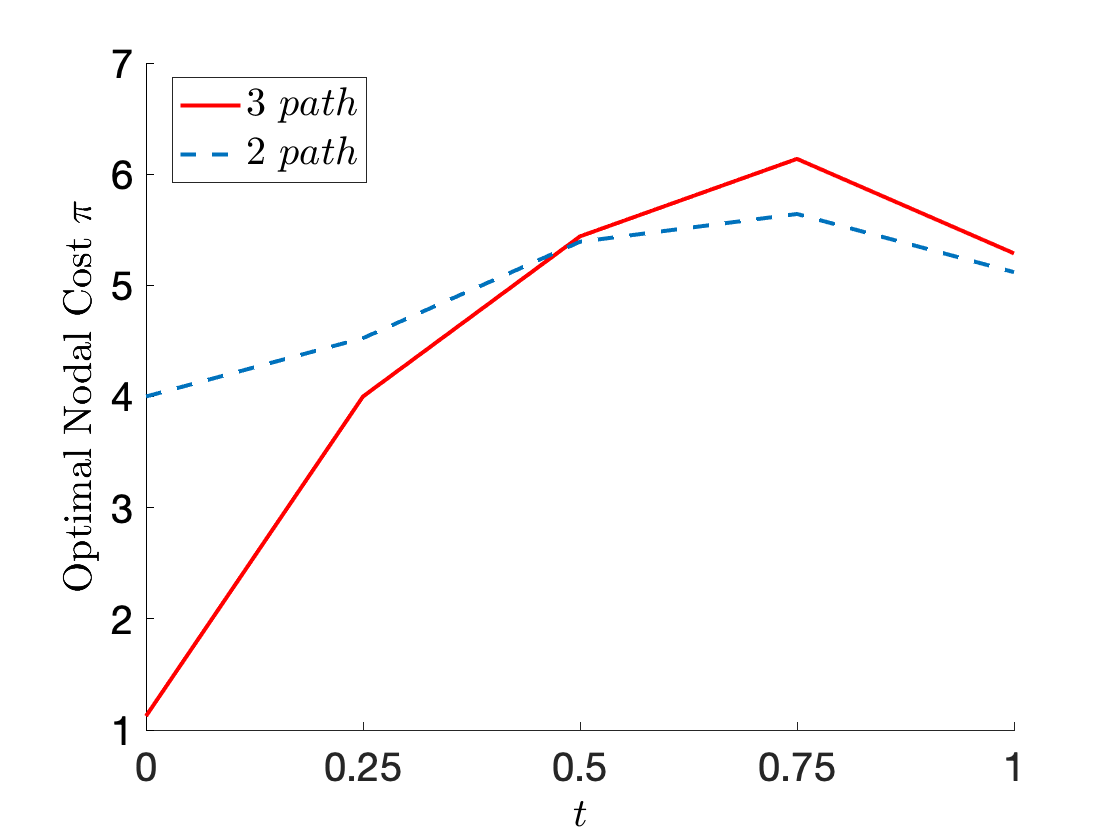}
    \caption{Optimal travel cost at the origin node on the two-path and three-path networks}
    \label{fig:paradox_cost}
\end{figure}

In Figure~\ref{fig:paradox_cost}, we plot the optimal nodal cost $\pi_1(t)$, which is the optimal travel cost of AVs entering the network through node 1 at time $t$, for $t\in[0,1]$. We observe that when $t<0.5$, AVs will have a lower travel cost on the three-path network than on the two-path network; but when $t>0.5$, AVs will have a higher cost on the three-path network. In other words, AVs can reduce their travel costs by utilizing the middle link if they enter the network early; but AVs entering the network at later times fall victim to the seemingly low-cost middle link. It is because at early times, there are a few cars on the network and the traffic density $\rho$ on each link is very small, so 
the links $(1,2)$ and $(3,4)$ are superior to links $(1,3)$ and $(2,4)$; at later times, there are more cars on the network and the links $(1,2)$ and $(3,4)$ become bottlenecks since AVs' travel costs on these two links are highly sensitive to their respective traffic densities.



\subsection{OW Network}\label{sec:num_ow}

In this subsection, we would like to test the solution algorithm presented in Section~\ref{sec:algo} on the OW network, which has 13 nodes and 24 links. Both [MFG-MiCP] and [MFG-LWR-MiCP] are solved on the OW network with the same running and queuing cost functions, and we compare the solved equilibria in terms of average velocity, occupied links and nodal optimal costs.

The settings of this experiment are as follows:

\begin{itemize}
    \item The simulation time $T=12$. 
    \item All links other than $(4,5)$ and $(5,9)$ have the same running cost function coefficients $c_1=c_2=1$ and $c_3=0.5$; The links $(4,5)$ and $(5,9)$ have the same running cost function coefficients $c_1=c_2=2$ and $c_3=5$.
    \item The bottleneck capacity $M_i=1$ for all nodes $i$. The queuing cost function coefficient $c_4=1$.
    \item AVs enter the network through the single origin $o=1$. The traffic demand is given by:
    \begin{align}
    d_1(t)&=\begin{cases}
          0.75,\quad &t\in[0,2];\\
          0,\quad &t\in(2,12].
    \end{cases}
    \end{align}
    \item The AVs' terminal costs at time $T$ for all nodes are given according to their distance from destination $s=13$ with $\pi_{13}(T)=0$. 
    
\end{itemize}

For both [MFG-MiCP] and [MFG-LWR-MiCP], we solve the equilibria by the developed solution algorithm. The initial values of of queue sizes are set to zero. 
The computation times for the MFG and LWR equilibria are 62.24s and 28.06s, respectively, on a laptop. That is, the presented solution algorithm is still efficient on the OW network.

\begin{figure}[H]
    \centering
    \begin{subfigure}{.8\textwidth}
        \includegraphics[width=\textwidth]{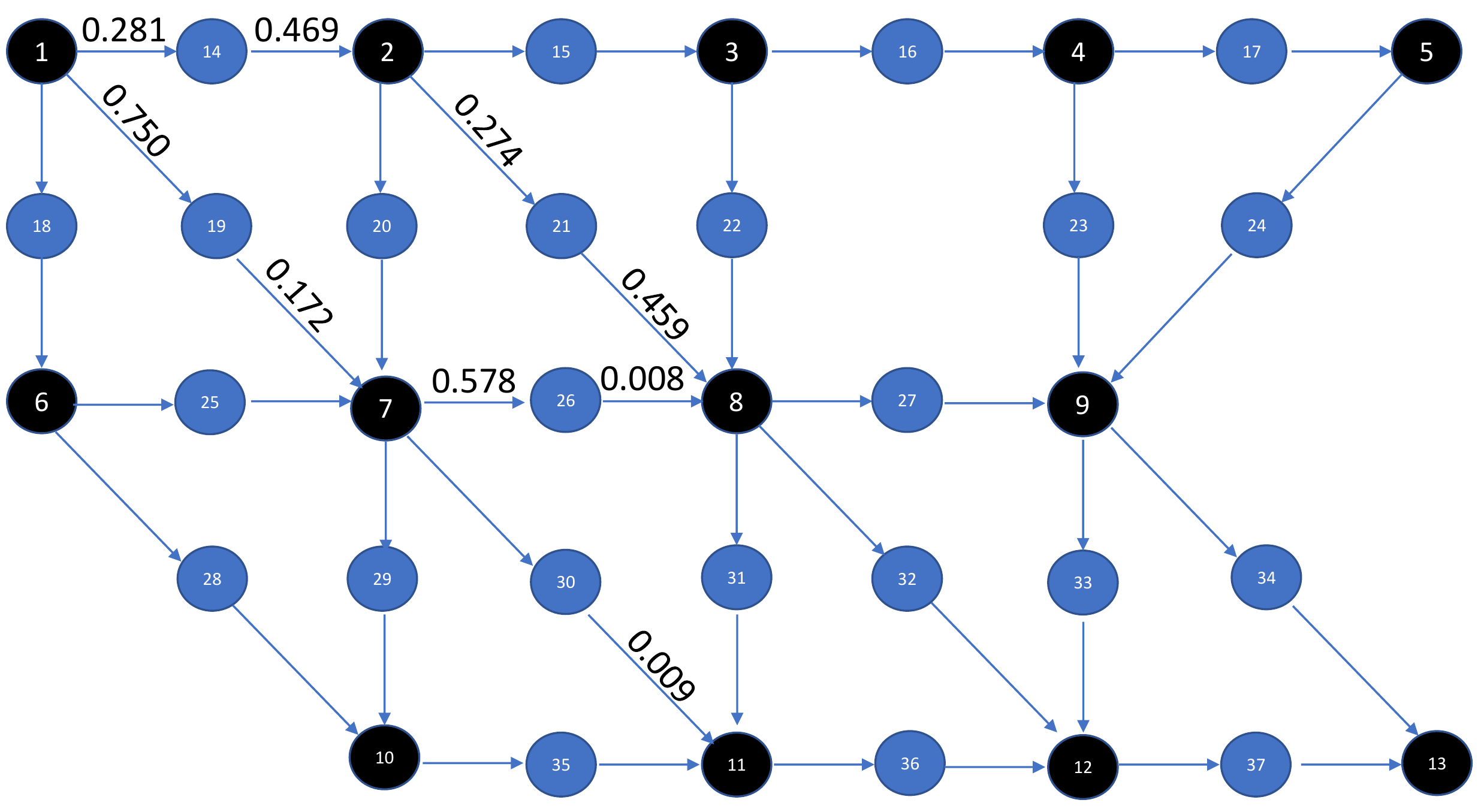}
        \caption{Traffic Density at $T=2$}\label{fig:ow_1}
    \end{subfigure}
    \hspace{3in}
    \begin{subfigure}{.8\textwidth}
        \includegraphics[width=\textwidth]{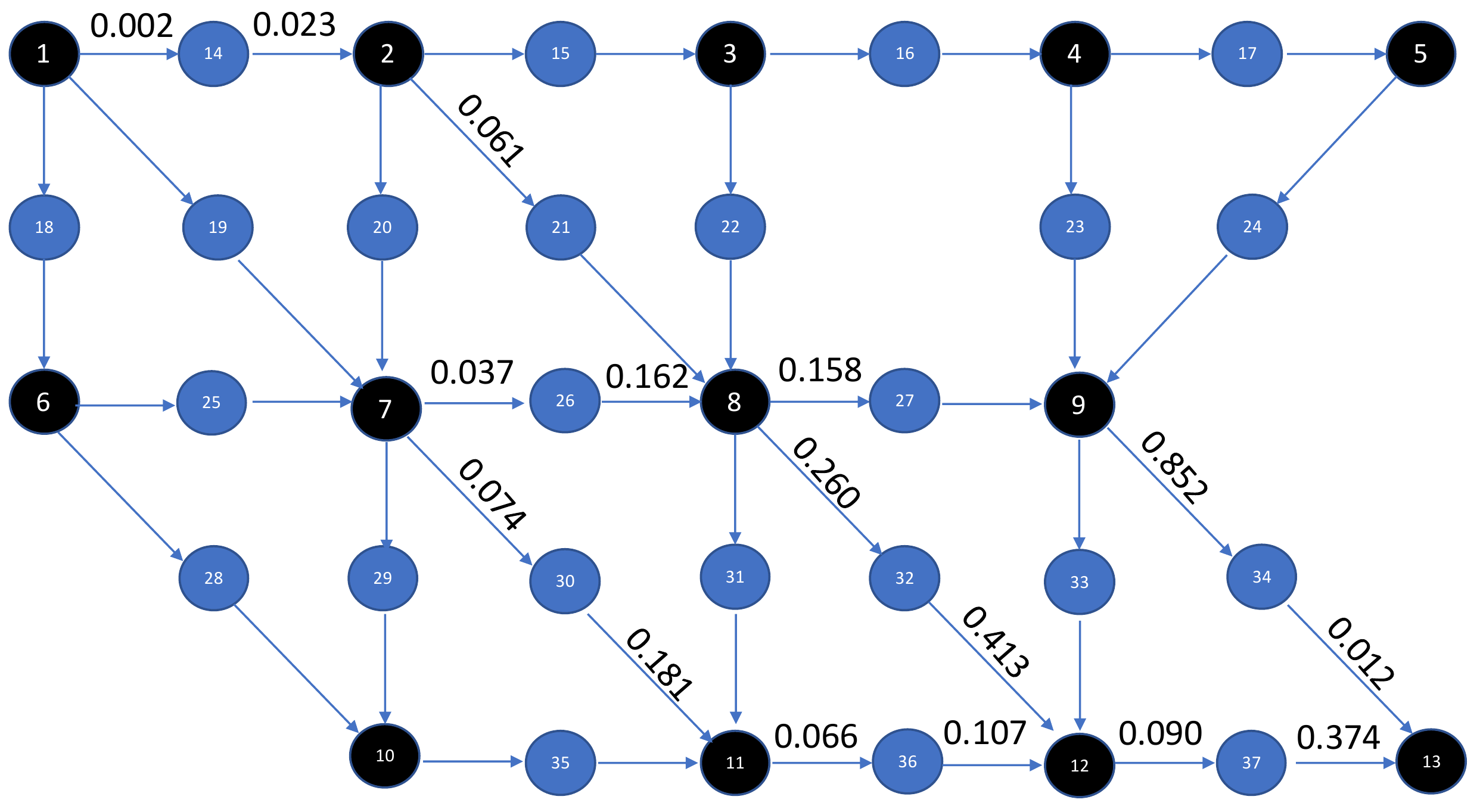}
        \caption{Traffic Density at $T=5$}\label{fig:ow_2}
    \end{subfigure}
    
    \caption{MFG equilibrium density evolution on the OW network}
    \label{fig:ow_flow}
\end{figure}

Figure~\ref{fig:ow_1} and Figure~\ref{fig:ow_2} show the density evolution of the solved MFG equilibrium at time $T=2$ and $T=5$, respectively. The black nodes are those in the OW network and the blue ones are auxiliary nodes introduced for link discretization. 

\begin{figure}[H]
    \centering
    \begin{subfigure}{.4\textwidth}
        \includegraphics[width=\textwidth]{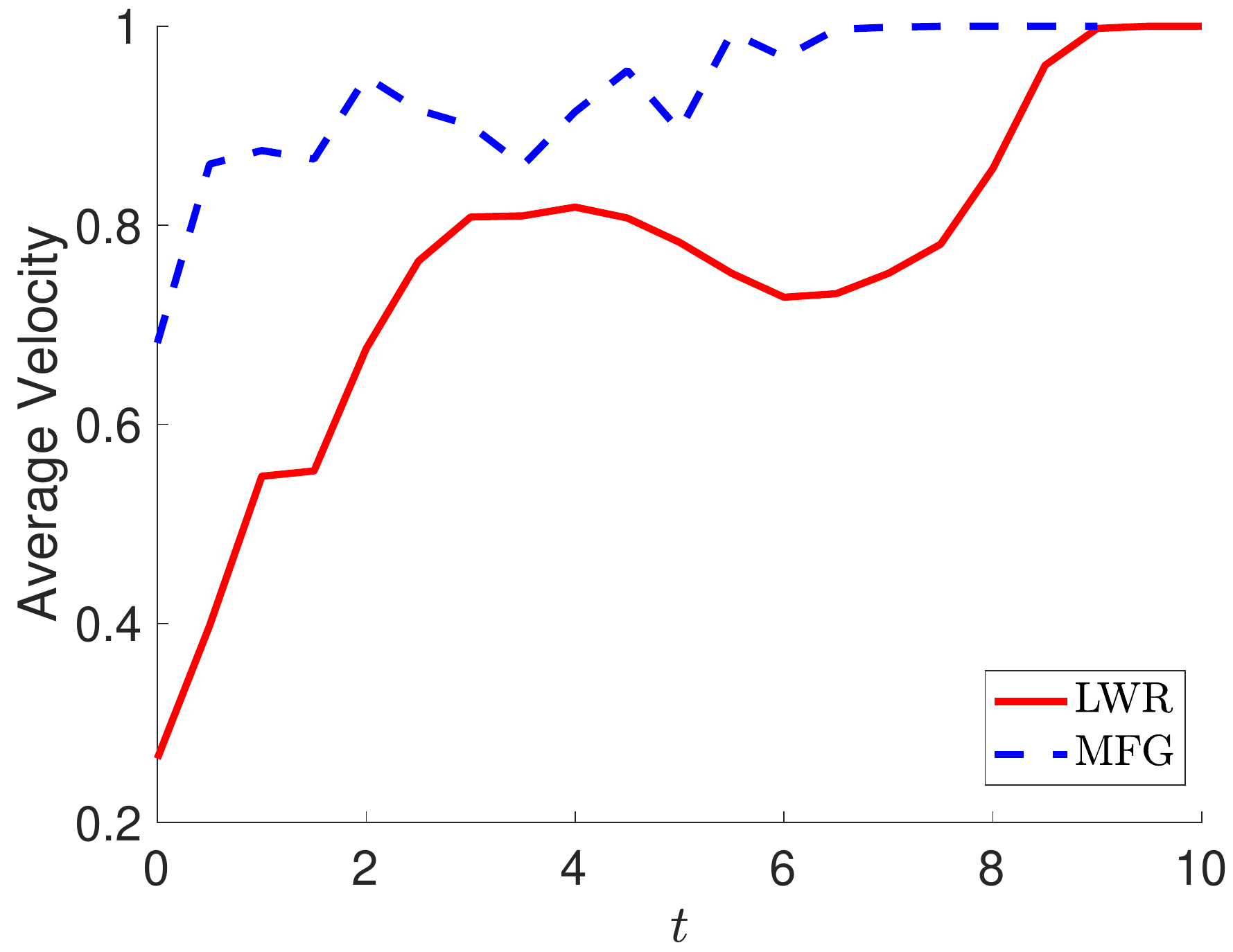}
        \caption{Average Velocity}\label{fig:ow_3}
    \end{subfigure}
    \begin{subfigure}{.4\textwidth}
        \includegraphics[width=\textwidth]{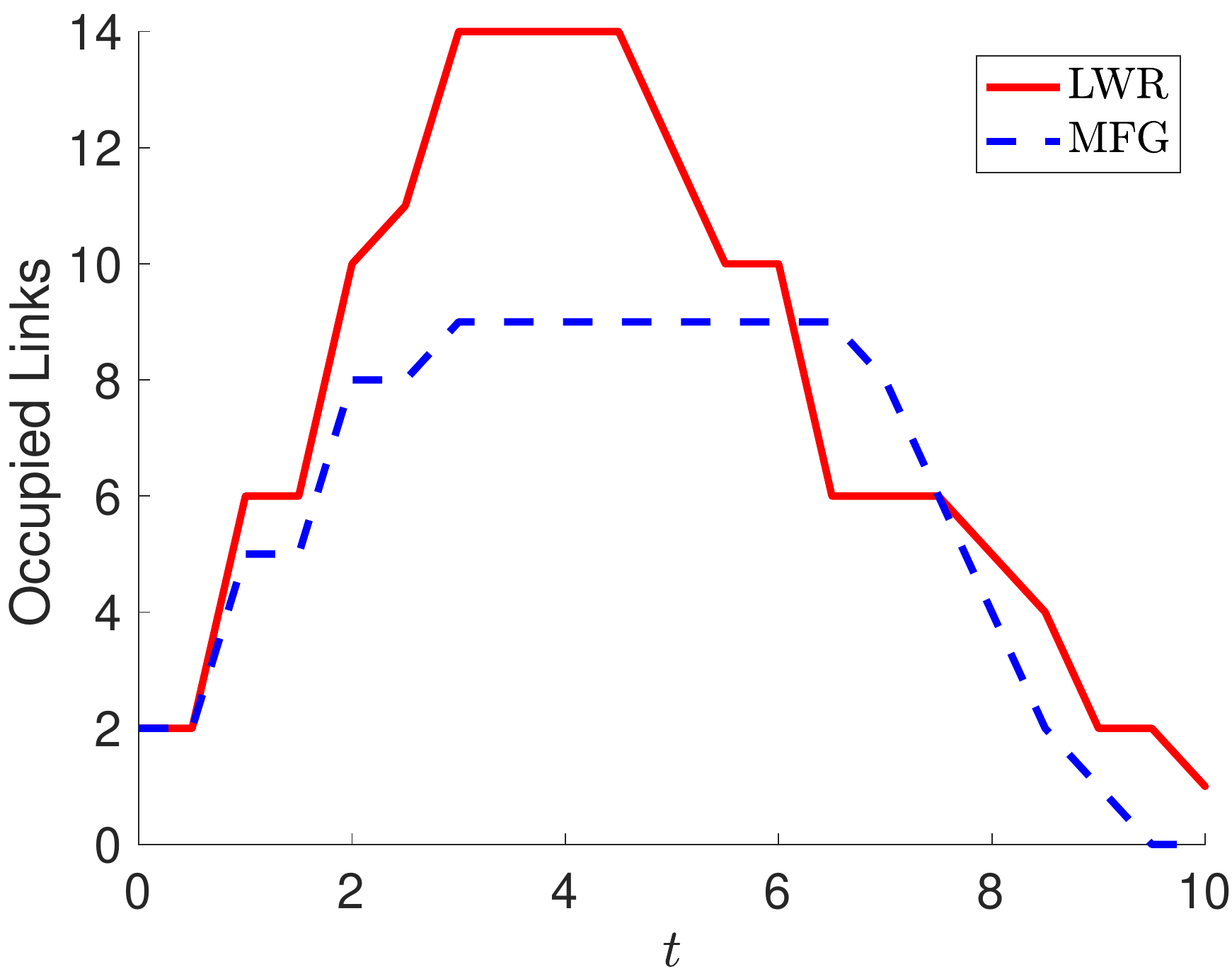}
        \caption{Occupied Links}\label{fig:ow_4}
    \end{subfigure}
    
    \begin{subfigure}{.4\textwidth}
        \includegraphics[width=\textwidth]{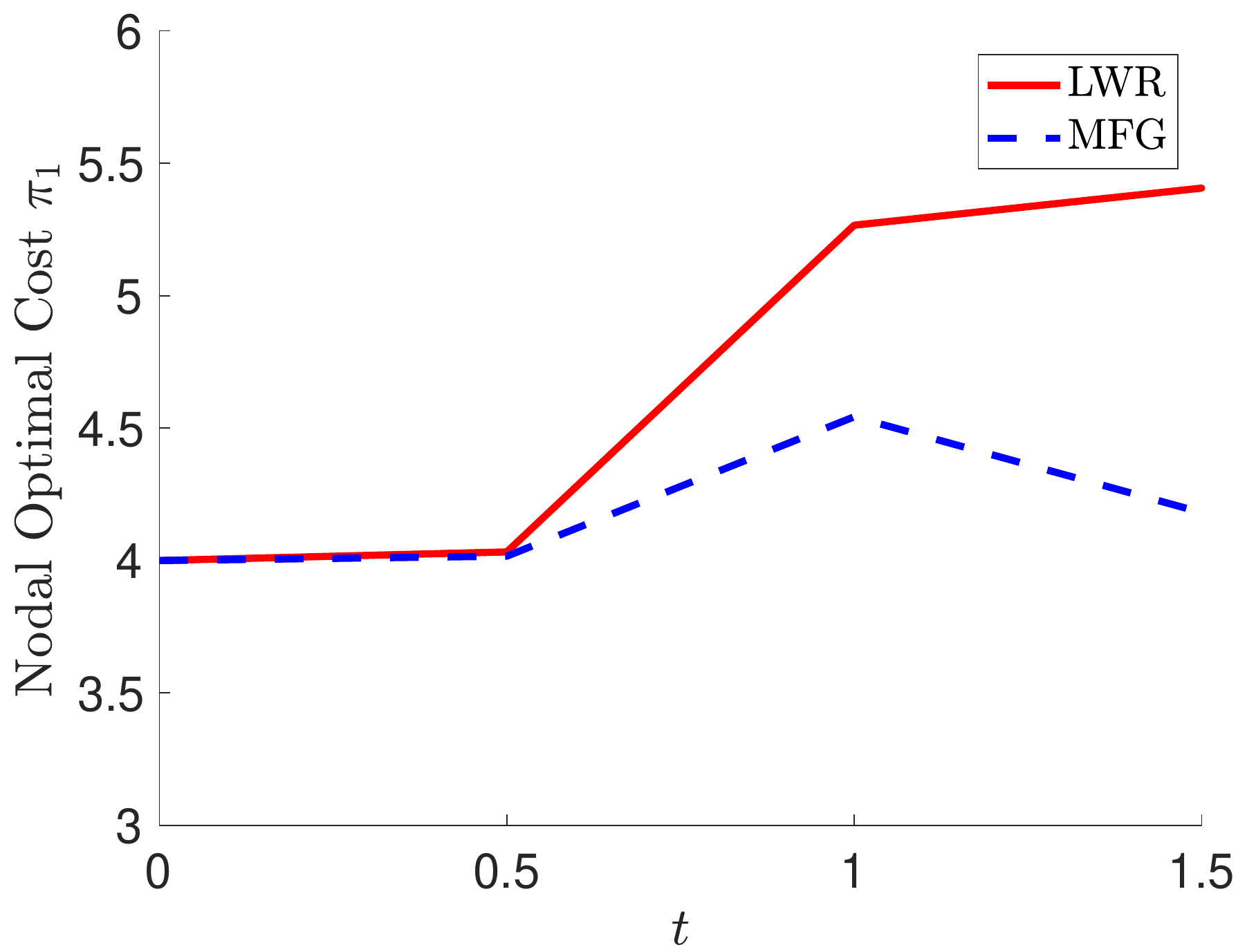}
        \caption{Nodal Optimal Cost}\label{fig:ow_5}
    \end{subfigure}
    
    \caption{Comparison between MFG and LWR equilibria on the OW network}
    \label{fig:ow_costs}
\end{figure}

Figure~\ref{fig:ow_3} shows the average velocity of all AVs on the network and Figure~\ref{fig:ow_4} shows the number of AVs' occupied links as a function of time $t$ for both MFG and LWR equilibria. We observe from the figures that AVs tend to select fewer links and drive at higher speeds in the MFG equilibrium than in the LWR equilibrium. Such a phenomenon is similar to what we observe on the three-path network: AVs deployed with MFG equilibrium controls can freely select their driving speeds to optimize their cumulative travel costs over a network. 

In Figure~\ref{fig:ow_5}, the optimal nodal cost $\pi_1(t)$ at the origin $1$ of the two equilibria is compared during the time horizon $t\in[0,1.5]$. We observe that before time $t=0.5$, the MFG and LWR equilibria have the same optimal nodal cost, possibly because AVs entering the network at early times experience little congestion and AVs in both MFG and LWR equilibria can have high speeds. While when $t>0.5$, cars entering the network at time $t$ may go into high density areas, and AVs in the LWR equilibrium are restricted to drive at very low speeds. Because AVs in the MFG equilibrium have more control power on the driving speed control, they can adjust their speeds and achieve lower travel costs. 

\color{black}

\section{Conclusions and future directions} \label{sec:conclud}

This paper illustrates the application of mean field games to the modeling of AVs' driving and routing game on a network. 
Mean field games (MFGs) have been widely used to design decision-making processes in finances and economics, with rather limited research in traffic modeling. 
Its application to AVs' decision choices on a network is a natural generalization of the classical DTA models. 
DTA models are popular tools to prescribe human drivers' en-route path choices and traffic flow evolution in the short term. 
However, they do not stipulate cars' velocity control while moving inside a road segment. 
Velocity control of AVs becomes increasingly important to navigate AVs across a congested traffic network. 
Therefore a new modeling paradigm is warranted. 
To model the movement of AVs in a transportation network, two choices are involved: driving speed control in the interior of a link, and next-go-to link choice at a junction node.   
To implement MFG, we decompose the problem into a backward HJB equation and a forward continuity equation. 
The derived MFG system is then reformulated as a MiCP. Based on the MiCP formulation, we illustrate the connection between MFG and classical DUE concept. An efficient numerical algorithm is developed to solve the MiCP, and the solved AVs' optimal controls are illustrated on the Braess network 
and the OW network. 
On the Braess network, by comparing the equilibria
under optimal speed selection and the LWR speed, the advantage of AVs' optimal velocity controls is demonstrated. Braess paradox is discovered in the AV traffic, indicating that we need to carefully design AV control algorithms to prevent the paradox from happening. 
On the OW network, by showing that our solution algorithm can successfully solve the equilibrium solution in a short computation time, we verify the efficiency of the solution algorithm.

This work can be extended as follows. 
(1) The proposed MFG system is defined for the many-to-one scenario.
A multi-class MFG can be developed to accommodate multiple destinations. 
(2) The equilibrium existence is given for a particular network and a family of cost functions. We would like to extend the theoretical analysis to more general networks and cost functions.
(3) We assume all cars to be AVs. When both AVs and human-driven vehicles co-exist on a transportation network, a multi-class mean field game can be developed to study the mixed traffic equilibrium.
Given the goal of this paper is the first-of-its-kind to lay the theoretical foundation for a scalable control framework of AVs' driving and routing game on a network, all the aforementioned analytical and computational issues will be left for future research.  

\section*{Acknowledgments}

The authors would like to thank Data Science Institute from
Columbia University for providing a seed grant for this research. The third author
acknowledges support from the National Science Foundation under the award number CMMI-1943998.
The fourth author
acknowledges support from the National Science Foundation under the award numbers CCF-1704833 and DMS-2012562.


\bibliographystyle{elsarticle-harv}

\bibliography{lit_MFG,survey_Di,survey_MFG,survey_DCL,DUE}

\begin{thebibliography}{71}
\expandafter\ifx\csname natexlab\endcsname\relax\def\natexlab#1{#1}\fi
\providecommand{\url}[1]{\texttt{#1}}
\providecommand{\href}[2]{#2}
\providecommand{\path}[1]{#1}
\providecommand{\DOIprefix}{doi:}
\providecommand{\ArXivprefix}{arXiv:}
\providecommand{\URLprefix}{URL: }
\providecommand{\Pubmedprefix}{pmid:}
\providecommand{\doi}[1]{\href{http://dx.doi.org/#1}{\path{#1}}}
\providecommand{\Pubmed}[1]{\href{pmid:#1}{\path{#1}}}
\providecommand{\bibinfo}[2]{#2}
\ifx\xfnm\relax \def\xfnm[#1]{\unskip,\space#1}\fi
\bibitem[{Ban et~al.(2008)Ban, Liu, Ferris and Ran}]{ban2008link}
\bibinfo{author}{Ban, X.J.}, \bibinfo{author}{Liu, H.X.},
  \bibinfo{author}{Ferris, M.C.}, \bibinfo{author}{Ran, B.},
  \bibinfo{year}{2008}.
\newblock \bibinfo{title}{A link-node complementarity model and solution
  algorithm for dynamic user equilibria with exact flow propagations}.
\newblock \bibinfo{journal}{Transportation Research Part B: Methodological}
  \bibinfo{volume}{42}, \bibinfo{pages}{823--842}.
\bibitem[{Ban et~al.(2012)Ban, Pang, Liu and Ma}]{ban2012continuous}
\bibinfo{author}{Ban, X.J.}, \bibinfo{author}{Pang, J.S.},
  \bibinfo{author}{Liu, H.X.}, \bibinfo{author}{Ma, R.}, \bibinfo{year}{2012}.
\newblock \bibinfo{title}{Continuous-time point-queue models in dynamic network
  loading}.
\newblock \bibinfo{journal}{Transportation Research Part B: Methodological}
  \bibinfo{volume}{46}, \bibinfo{pages}{360--380}.
\bibitem[{Batista et~al.(2021)Batista, Leclercq and
  MenÃ©ndez}]{batista2021dta}
\bibinfo{author}{Batista, S.}, \bibinfo{author}{Leclercq, L.},
  \bibinfo{author}{MenÃ©ndez, M.}, \bibinfo{year}{2021}.
\newblock \bibinfo{title}{Dynamic traffic assignment for regional networks with
  traffic-dependent trip lengths and regional paths}.
\newblock \bibinfo{journal}{Transportation Research Part C: Emerging
  Technologies} \bibinfo{volume}{127}, \bibinfo{pages}{103076}.
\bibitem[{Bauso et~al.(2016)Bauso, Zhang and
  Papachristodoulou}]{bauso2016density}
\bibinfo{author}{Bauso, D.}, \bibinfo{author}{Zhang, X.},
  \bibinfo{author}{Papachristodoulou, A.}, \bibinfo{year}{2016}.
\newblock \bibinfo{title}{Density flow in dynamical networks via mean-field
  games}.
\newblock \bibinfo{journal}{IEEE Transactions on Automatic Control}
  \bibinfo{volume}{62}, \bibinfo{pages}{1342--1355}.
\bibitem[{Bliemer et~al.(2017)Bliemer, Raadsen, Brederode, Bell, Wismans and
  Smith}]{bliemer2017genetics}
\bibinfo{author}{Bliemer, M.C.J.}, \bibinfo{author}{Raadsen, M.P.H.},
  \bibinfo{author}{Brederode, L.J.N.}, \bibinfo{author}{Bell, M.G.H.},
  \bibinfo{author}{Wismans, L.J.J.}, \bibinfo{author}{Smith, M.J.},
  \bibinfo{year}{2017}.
\newblock \bibinfo{title}{Genetics of traffic assignment models for strategic
  transport planning}.
\newblock \bibinfo{journal}{Transport Reviews} \bibinfo{volume}{37}.
\bibitem[{Boyce et~al.(1995)Boyce, Ran and Leblanc}]{boyce1995solving}
\bibinfo{author}{Boyce, D.E.}, \bibinfo{author}{Ran, B.},
  \bibinfo{author}{Leblanc, L.J.}, \bibinfo{year}{1995}.
\newblock \bibinfo{title}{Solving an instantaneous dynamic user-optimal route
  choice model}.
\newblock \bibinfo{journal}{Transportation Science} \bibinfo{volume}{29},
  \bibinfo{pages}{128--142}.
\bibitem[{Burger et~al.(2013)Burger, Di~Francesco, Markowich and
  Wolfram}]{burger2013mean}
\bibinfo{author}{Burger, M.}, \bibinfo{author}{Di~Francesco, M.},
  \bibinfo{author}{Markowich, P.}, \bibinfo{author}{Wolfram, M.T.},
  \bibinfo{year}{2013}.
\newblock \bibinfo{title}{Mean field games with nonlinear mobilities in
  pedestrian dynamics}.
\newblock \bibinfo{journal}{arXiv preprint arXiv:1304.5201} .
\bibitem[{Cardaliaguet(2010)}]{cardaliaguet2010notes}
\bibinfo{author}{Cardaliaguet, P.}, \bibinfo{year}{2010}.
\newblock \bibinfo{title}{Notes on mean field games}.
\newblock \bibinfo{type}{Technical Report}.
\bibitem[{Cardaliaguet(2015)}]{cardaliaguet2015weak}
\bibinfo{author}{Cardaliaguet, P.}, \bibinfo{year}{2015}.
\newblock \bibinfo{title}{Weak solutions for first order mean field games with
  local coupling}, in: \bibinfo{booktitle}{Analysis and geometry in control
  theory and its applications}. \bibinfo{publisher}{Springer}, pp.
  \bibinfo{pages}{111--158}.
\bibitem[{Carey and Ge(2003)}]{carey2003comparing}
\bibinfo{author}{Carey, M.}, \bibinfo{author}{Ge, Y.}, \bibinfo{year}{2003}.
\newblock \bibinfo{title}{Comparing whole-link travel time models}.
\newblock \bibinfo{journal}{Transportation Research Part B: Methodological}
  \bibinfo{volume}{37}, \bibinfo{pages}{905--926}.
\bibitem[{Carey and McCartney(2002)}]{carey2002behaviour}
\bibinfo{author}{Carey, M.}, \bibinfo{author}{McCartney, M.},
  \bibinfo{year}{2002}.
\newblock \bibinfo{title}{Behaviour of a whole-link travel time model used in
  dynamic traffic assignment}.
\newblock \bibinfo{journal}{Transportation Research Part B: Methodological}
  \bibinfo{volume}{36}, \bibinfo{pages}{83--95}.
\bibitem[{Chen et~al.(2016)Chen, He and Yin}]{chen2016optimal}
\bibinfo{author}{Chen, Z.}, \bibinfo{author}{He, F.}, \bibinfo{author}{Yin,
  Y.}, \bibinfo{year}{2016}.
\newblock \bibinfo{title}{Optimal deployment of charging lanes for electric
  vehicles in transportation networks}.
\newblock \bibinfo{journal}{Transportation Research Part B: Methodological}
  \bibinfo{volume}{91}, \bibinfo{pages}{344--365}.
\bibitem[{Chevalier et~al.(2015)Chevalier, Le~Ny and
  Malham{\'e}}]{chevalier2015micro}
\bibinfo{author}{Chevalier, G.}, \bibinfo{author}{Le~Ny, J.},
  \bibinfo{author}{Malham{\'e}, R.}, \bibinfo{year}{2015}.
\newblock \bibinfo{title}{A micro-macro traffic model based on mean-field
  games}, in: \bibinfo{booktitle}{2015 American Control Conference (ACC)},
  \bibinfo{organization}{IEEE}. pp. \bibinfo{pages}{1983--1988}.
\bibitem[{Chiu et~al.(2013)Chiu, Khani, Noh, Bustillos and
  Hickman}]{chiu2013technical}
\bibinfo{author}{Chiu, Y.C.}, \bibinfo{author}{Khani, A.},
  \bibinfo{author}{Noh, H.}, \bibinfo{author}{Bustillos, B.},
  \bibinfo{author}{Hickman, M.}, \bibinfo{year}{2013}.
\newblock \bibinfo{title}{Technical report on SHRP 2 C10B version of DynusT and
  FAST-TrIPs}.
\newblock \bibinfo{type}{Technical Report}.
\bibitem[{Couillet et~al.(2012)Couillet, Perlaza, Tembine and
  Debbah}]{couillet2012electrical}
\bibinfo{author}{Couillet, R.}, \bibinfo{author}{Perlaza, S.M.},
  \bibinfo{author}{Tembine, H.}, \bibinfo{author}{Debbah, M.},
  \bibinfo{year}{2012}.
\newblock \bibinfo{title}{Electrical vehicles in the smart grid: A mean field
  game analysis}.
\newblock \bibinfo{journal}{IEEE Journal on Selected Areas in Communications}
  \bibinfo{volume}{30}, \bibinfo{pages}{1086--1096}.
\bibitem[{Degond et~al.(2014)Degond, Liu and Ringhofer}]{degond2014large}
\bibinfo{author}{Degond, P.}, \bibinfo{author}{Liu, J.G.},
  \bibinfo{author}{Ringhofer, C.}, \bibinfo{year}{2014}.
\newblock \bibinfo{title}{Large-scale dynamics of mean-field games driven by
  local nash equilibria}.
\newblock \bibinfo{journal}{Journal of Nonlinear Science} \bibinfo{volume}{24},
  \bibinfo{pages}{93--115}.
\bibitem[{Di et~al.(2014)Di, He, Guo and Liu}]{di2014braess}
\bibinfo{author}{Di, X.}, \bibinfo{author}{He, X.}, \bibinfo{author}{Guo, X.},
  \bibinfo{author}{Liu, H.X.}, \bibinfo{year}{2014}.
\newblock \bibinfo{title}{Braess paradox under the boundedly rational user
  equilibria}.
\newblock \bibinfo{journal}{Transportation Research Part B}
  \bibinfo{volume}{67}, \bibinfo{pages}{86--108}.
\bibitem[{Di and Liu(2016)}]{di2016boundedly}
\bibinfo{author}{Di, X.}, \bibinfo{author}{Liu, H.X.}, \bibinfo{year}{2016}.
\newblock \bibinfo{title}{Boundedly rational route choice behavior: A review of
  models and methodologies}.
\newblock \bibinfo{journal}{Transportation Research Part B}
  \bibinfo{volume}{85}, \bibinfo{pages}{142--179}.
\bibitem[{Di et~al.(2013)Di, Liu, Pang and Ban}]{di2013boundedlyTRB}
\bibinfo{author}{Di, X.}, \bibinfo{author}{Liu, H.X.}, \bibinfo{author}{Pang,
  J.S.}, \bibinfo{author}{Ban, X.J.}, \bibinfo{year}{2013}.
\newblock \bibinfo{title}{Boundedly rational user equilibria ({BRUE}):
  Mathematical formulation and solution sets}.
\newblock \bibinfo{journal}{Transportation Research Part B} ,
  \bibinfo{pages}{300–--313}.
\bibitem[{Djehiche et~al.(2016)Djehiche, Tcheukam and
  Tembine}]{djehiche2016mean}
\bibinfo{author}{Djehiche, B.}, \bibinfo{author}{Tcheukam, A.},
  \bibinfo{author}{Tembine, H.}, \bibinfo{year}{2016}.
\newblock \bibinfo{title}{Mean-field-type games in engineering}.
\newblock \bibinfo{journal}{arXiv preprint arXiv:1605.03281} .
\bibitem[{Doan and Ukkusuri(2015)}]{doan2015queue}
\bibinfo{author}{Doan, K.}, \bibinfo{author}{Ukkusuri, S.V.},
  \bibinfo{year}{2015}.
\newblock \bibinfo{title}{Dynamic system optimal model for multi-od traffic
  networks with an advanced spatial queuing model}.
\newblock \bibinfo{journal}{Transportation Research Part C: Emerging
  Technologies} \bibinfo{volume}{51}, \bibinfo{pages}{41--65}.
\bibitem[{Dockner et~al.(2000)Dockner, Jorgensen, Van~Long and
  Sorger}]{dockner2000differential}
\bibinfo{author}{Dockner, E.J.}, \bibinfo{author}{Jorgensen, S.},
  \bibinfo{author}{Van~Long, N.}, \bibinfo{author}{Sorger, G.},
  \bibinfo{year}{2000}.
\newblock \bibinfo{title}{Differential games in economics and management
  science}.
\newblock \bibinfo{publisher}{Cambridge University Press}.
\bibitem[{Facchinei and Pang(2007)}]{facchinei2007finite}
\bibinfo{author}{Facchinei, F.}, \bibinfo{author}{Pang, J.S.},
  \bibinfo{year}{2007}.
\newblock \bibinfo{title}{Finite-dimensional variational inequalities and
  complementarity problems}.
\newblock \bibinfo{publisher}{Springer Science \& Business Media}.
\bibitem[{Festa and G{\"o}ttlich(2017)}]{festa2017mean}
\bibinfo{author}{Festa, A.}, \bibinfo{author}{G{\"o}ttlich, S.},
  \bibinfo{year}{2017}.
\newblock \bibinfo{title}{A mean field games approach for multi-lane traffic
  management}.
\newblock \bibinfo{journal}{arXiv preprint arXiv:1711.04116} .
\bibitem[{Friesz et~al.(1993)Friesz, Bernstein, Smith, Tobin and
  Wie}]{friesz1993variational}
\bibinfo{author}{Friesz, T.L.}, \bibinfo{author}{Bernstein, D.},
  \bibinfo{author}{Smith, T.E.}, \bibinfo{author}{Tobin, R.L.},
  \bibinfo{author}{Wie, B.W.}, \bibinfo{year}{1993}.
\newblock \bibinfo{title}{A variational inequality formulation of the dynamic
  network user equilibrium problem}.
\newblock \bibinfo{journal}{Operations research} \bibinfo{volume}{41},
  \bibinfo{pages}{179--191}.
\bibitem[{Friesz and Han(2019)}]{friesz2019mathematical}
\bibinfo{author}{Friesz, T.L.}, \bibinfo{author}{Han, K.},
  \bibinfo{year}{2019}.
\newblock \bibinfo{title}{The mathematical foundations of dynamic user
  equilibrium}.
\newblock \bibinfo{journal}{Transportation research part B: methodological}
  \bibinfo{volume}{126}, \bibinfo{pages}{309--328}.
\bibitem[{Friesz et~al.(2013)Friesz, Han, Neto, Meimand and
  Yao}]{friesz2013dynamic}
\bibinfo{author}{Friesz, T.L.}, \bibinfo{author}{Han, K.},
  \bibinfo{author}{Neto, P.A.}, \bibinfo{author}{Meimand, A.},
  \bibinfo{author}{Yao, T.}, \bibinfo{year}{2013}.
\newblock \bibinfo{title}{Dynamic user equilibrium based on a hydrodynamic
  model}.
\newblock \bibinfo{journal}{Transportation Research Part B: Methodological}
  \bibinfo{volume}{47}, \bibinfo{pages}{102--126}.
\bibitem[{Gentile(2015)}]{gentile2015using}
\bibinfo{author}{Gentile, G.}, \bibinfo{year}{2015}.
\newblock \bibinfo{title}{Using the general link transmission model in a
  dynamic traffic assignment to simulate congestion on urban networks}.
\newblock \bibinfo{journal}{Transportation Research Procedia}
  \bibinfo{volume}{5}, \bibinfo{pages}{66--81}.
\bibitem[{Gentile et~al.(2007)Gentile, Meschini and
  Papola}]{gentile2007spillback}
\bibinfo{author}{Gentile, G.}, \bibinfo{author}{Meschini, L.},
  \bibinfo{author}{Papola, N.}, \bibinfo{year}{2007}.
\newblock \bibinfo{title}{Spillback congestion in dynamic traffic assignment: a
  macroscopic flow model with time-varying bottlenecks}.
\newblock \bibinfo{journal}{Transportation Research Part B: Methodological}
  \bibinfo{volume}{41}, \bibinfo{pages}{1114--1138}.
\bibitem[{Gu{\'e}ant et~al.(2011)Gu{\'e}ant, Lasry and Lions}]{gueant2011mean}
\bibinfo{author}{Gu{\'e}ant, O.}, \bibinfo{author}{Lasry, J.M.},
  \bibinfo{author}{Lions, P.L.}, \bibinfo{year}{2011}.
\newblock \bibinfo{title}{Mean field games and applications}, in:
  \bibinfo{booktitle}{Paris-Princeton lectures on mathematical finance 2010}.
  \bibinfo{publisher}{Springer}, pp. \bibinfo{pages}{205--266}.
\bibitem[{Gummadi et~al.(2012)Gummadi, Johari and Yu}]{gummadi2012mean}
\bibinfo{author}{Gummadi, R.}, \bibinfo{author}{Johari, R.},
  \bibinfo{author}{Yu, J.Y.}, \bibinfo{year}{2012}.
\newblock \bibinfo{title}{Mean field equilibria of multi armed bandit games},
  in: \bibinfo{booktitle}{2012 50th Annual Allerton Conference on
  Communication, Control, and Computing (Allerton)},
  \bibinfo{organization}{IEEE}. pp. \bibinfo{pages}{1110--1110}.
\bibitem[{Han et~al.(2019)Han, Eve and Friesz}]{han2019computing}
\bibinfo{author}{Han, K.}, \bibinfo{author}{Eve, G.}, \bibinfo{author}{Friesz,
  T.L.}, \bibinfo{year}{2019}.
\newblock \bibinfo{title}{Computing dynamic user equilibria on large-scale
  networks with software implementation}.
\newblock \bibinfo{journal}{Networks and Spatial Economics} ,
  \bibinfo{pages}{1--34}.
\bibitem[{Han et~al.(2013)Han, Friesz and Yao}]{han2013partial1}
\bibinfo{author}{Han, K.}, \bibinfo{author}{Friesz, T.L.},
  \bibinfo{author}{Yao, T.}, \bibinfo{year}{2013}.
\newblock \bibinfo{title}{A partial differential equation formulation of
  vickrey’s bottleneck model, part i: Methodology and theoretical analysis}.
\newblock \bibinfo{journal}{Transportation Research Part B: Methodological}
  \bibinfo{volume}{49}, \bibinfo{pages}{55--74}.
\bibitem[{Han et~al.(2016)Han, Piccoli and Friesz}]{han2016continuity}
\bibinfo{author}{Han, K.}, \bibinfo{author}{Piccoli, B.},
  \bibinfo{author}{Friesz, T.L.}, \bibinfo{year}{2016}.
\newblock \bibinfo{title}{Continuity of the path delay operator for dynamic
  network loading with spillback}.
\newblock \bibinfo{journal}{Transportation Research Part B: Methodological}
  \bibinfo{volume}{92}, \bibinfo{pages}{211--233}.
\bibitem[{Huang et~al.(2019)Huang, Di, Du and Chen}]{huang2019stable}
\bibinfo{author}{Huang, K.}, \bibinfo{author}{Di, X.}, \bibinfo{author}{Du,
  Q.}, \bibinfo{author}{Chen, X.}, \bibinfo{year}{2019}.
\newblock \bibinfo{title}{Stabilizing traffic via autonomous vehicles: A
  continuum mean field game approach}, in: \bibinfo{booktitle}{2019 IEEE
  Intelligent Transportation Systems Conference (ITSC)},
  \bibinfo{organization}{IEEE}. pp. \bibinfo{pages}{3269--3274}.
\bibitem[{Huang et~al.(2020a)Huang, Di, Du and Chen}]{huang2019game}
\bibinfo{author}{Huang, K.}, \bibinfo{author}{Di, X.}, \bibinfo{author}{Du,
  Q.}, \bibinfo{author}{Chen, X.}, \bibinfo{year}{2020}a.
\newblock \bibinfo{title}{A game-theoretic framework for autonomous vehicles
  velocity control: Bridging microscopic differential games and macroscopic
  mean field games}.
\newblock \bibinfo{journal}{Discrete and Continuous Dynamical Systems - Series
  B} \bibinfo{volume}{25}, \bibinfo{pages}{4869--4903}.
\bibitem[{Huang et~al.(2020b)Huang, Di, Du and Chen}]{huang2019mixed}
\bibinfo{author}{Huang, K.}, \bibinfo{author}{Di, X.}, \bibinfo{author}{Du,
  Q.}, \bibinfo{author}{Chen, X.}, \bibinfo{year}{2020}b.
\newblock \bibinfo{title}{Scalable traffic stability analysis in mixed-autonomy
  using continuum models}.
\newblock \bibinfo{journal}{Transportation Research Part C: Emerging
  Technologies} \bibinfo{volume}{111}, \bibinfo{pages}{616--630}.
\bibitem[{Iyer et~al.(2014)Iyer, Johari and Sundararajan}]{iyer2014mean}
\bibinfo{author}{Iyer, K.}, \bibinfo{author}{Johari, R.},
  \bibinfo{author}{Sundararajan, M.}, \bibinfo{year}{2014}.
\newblock \bibinfo{title}{Mean field equilibria of dynamic auctions with
  learning}.
\newblock \bibinfo{journal}{Management Science} \bibinfo{volume}{60},
  \bibinfo{pages}{2949--2970}.
\bibitem[{Kachroo et~al.(2017)Kachroo, Agarwal, Piccoli and
  {\"O}zbay}]{kachroo2017multiscale}
\bibinfo{author}{Kachroo, P.}, \bibinfo{author}{Agarwal, S.},
  \bibinfo{author}{Piccoli, B.}, \bibinfo{author}{{\"O}zbay, K.},
  \bibinfo{year}{2017}.
\newblock \bibinfo{title}{Multiscale modeling and control architecture for v2x
  enabled traffic streams}.
\newblock \bibinfo{journal}{IEEE Transactions on Vehicular Technology}
  \bibinfo{volume}{66}, \bibinfo{pages}{4616--4626}.
\bibitem[{Kachroo et~al.(2016)Kachroo, Agarwal and Sastry}]{kachroo2016inverse}
\bibinfo{author}{Kachroo, P.}, \bibinfo{author}{Agarwal, S.},
  \bibinfo{author}{Sastry, S.}, \bibinfo{year}{2016}.
\newblock \bibinfo{title}{Inverse problem for non-viscous mean field control:
  Example from traffic}.
\newblock \bibinfo{journal}{IEEE Transactions on Automatic Control}
  \bibinfo{volume}{61}, \bibinfo{pages}{3412--3421}.
\bibitem[{Kuwahara and Akamatsu(2001)}]{kuwahara2001dynamic}
\bibinfo{author}{Kuwahara, M.}, \bibinfo{author}{Akamatsu, T.},
  \bibinfo{year}{2001}.
\newblock \bibinfo{title}{Dynamic user optimal assignment with physical queues
  for a many-to-many od pattern}.
\newblock \bibinfo{journal}{Transportation Research Part B: Methodological}
  \bibinfo{volume}{35}, \bibinfo{pages}{461--479}.
\bibitem[{Lachapelle et~al.(2010)Lachapelle, Salomon and
  Turinici}]{lachapelle2010computation}
\bibinfo{author}{Lachapelle, A.}, \bibinfo{author}{Salomon, J.},
  \bibinfo{author}{Turinici, G.}, \bibinfo{year}{2010}.
\newblock \bibinfo{title}{Computation of mean field equilibria in economics}.
\newblock \bibinfo{journal}{Mathematical Models and Methods in Applied
  Sciences} \bibinfo{volume}{20}, \bibinfo{pages}{567--588}.
\bibitem[{Lachapelle and Wolfram(2011)}]{lachapelle2011mean}
\bibinfo{author}{Lachapelle, A.}, \bibinfo{author}{Wolfram, M.T.},
  \bibinfo{year}{2011}.
\newblock \bibinfo{title}{On a mean field game approach modeling congestion and
  aversion in pedestrian crowds}.
\newblock \bibinfo{journal}{Transportation research part B: methodological}
  \bibinfo{volume}{45}, \bibinfo{pages}{1572--1589}.
\bibitem[{Lam and Huang(1995)}]{lam1995dynamic}
\bibinfo{author}{Lam, W.H.}, \bibinfo{author}{Huang, H.J.},
  \bibinfo{year}{1995}.
\newblock \bibinfo{title}{Dynamic user optimal traffic assignment model for
  many to one travel demand}.
\newblock \bibinfo{journal}{Transportation Research Part B: Methodological}
  \bibinfo{volume}{29}, \bibinfo{pages}{243--259}.
\bibitem[{Lebacque and Khoshyaran(1999)}]{lebacque1999modelling}
\bibinfo{author}{Lebacque, J.P.}, \bibinfo{author}{Khoshyaran, M.},
  \bibinfo{year}{1999}.
\newblock \bibinfo{title}{Modelling vehicular traffic flow on networks using
  macroscopic models}.
\newblock \bibinfo{journal}{Finite volumes for complex applications II} ,
  \bibinfo{pages}{551--558}.
\bibitem[{LeVeque(2002)}]{leveque2002finite}
\bibinfo{author}{LeVeque, R.J.}, \bibinfo{year}{2002}.
\newblock \bibinfo{title}{Finite volume methods for hyperbolic problems}.
  volume~\bibinfo{volume}{31}.
\newblock \bibinfo{publisher}{Cambridge university press}.
\bibitem[{LeVeque(2007)}]{leveque2007finite}
\bibinfo{author}{LeVeque, R.J.}, \bibinfo{year}{2007}.
\newblock \bibinfo{title}{Finite difference methods for ordinary and partial
  differential equations: steady-state and time-dependent problems}.
\newblock \bibinfo{publisher}{SIAM}.
\bibitem[{Levin and Boyles(2016)}]{levin2016multiclass}
\bibinfo{author}{Levin, M.W.}, \bibinfo{author}{Boyles, S.D.},
  \bibinfo{year}{2016}.
\newblock \bibinfo{title}{A multiclass cell transmission model for shared human
  and autonomous vehicle roads}.
\newblock \bibinfo{journal}{Transportation Research Part C: Emerging
  Technologies} \bibinfo{volume}{62}, \bibinfo{pages}{103--116}.
\bibitem[{Lighthill(1952)}]{lighthill1952sound}
\bibinfo{author}{Lighthill, M.J.}, \bibinfo{year}{1952}.
\newblock \bibinfo{title}{On sound generated aerodynamically {I}. general
  theory}.
\newblock \bibinfo{journal}{Proc. R. Soc. Lond. A} \bibinfo{volume}{211},
  \bibinfo{pages}{564--587}.
\bibitem[{Lighthill and Whitham(1955)}]{lighthill1955kinematic}
\bibinfo{author}{Lighthill, M.J.}, \bibinfo{author}{Whitham, G.B.},
  \bibinfo{year}{1955}.
\newblock \bibinfo{title}{On kinematic waves {II}. {A} theory of traffic flow
  on long crowded roads}.
\newblock \bibinfo{journal}{Proc. R. Soc. Lond. A} \bibinfo{volume}{229},
  \bibinfo{pages}{317--345}.
\bibitem[{Lo(1999)}]{lo1999dynamic}
\bibinfo{author}{Lo, H.}, \bibinfo{year}{1999}.
\newblock \bibinfo{title}{A dynamic traffic assignment formulation that
  encapsulates the cell-transmission model}, in: \bibinfo{booktitle}{14th
  International Symposium on Transportation and Traffic TheoryTransportation
  Research Institute}.
\bibitem[{Lo and Szeto(2002)}]{lo2002cell}
\bibinfo{author}{Lo, H.K.}, \bibinfo{author}{Szeto, W.Y.},
  \bibinfo{year}{2002}.
\newblock \bibinfo{title}{A cell-based variational inequality formulation of
  the dynamic user optimal assignment problem}.
\newblock \bibinfo{journal}{Transportation Research Part B: Methodological}
  \bibinfo{volume}{36}, \bibinfo{pages}{421--443}.
\bibitem[{Mahmassani(2001)}]{mahmassani2001dynamic}
\bibinfo{author}{Mahmassani, H.S.}, \bibinfo{year}{2001}.
\newblock \bibinfo{title}{Dynamic network traffic assignment and simulation
  methodology for advanced system management applications}.
\newblock \bibinfo{journal}{Networks and spatial economics}
  \bibinfo{volume}{1}, \bibinfo{pages}{267--292}.
\bibitem[{Merchant and Nemhauser(1978a)}]{merchant1978model}
\bibinfo{author}{Merchant, D.K.}, \bibinfo{author}{Nemhauser, G.L.},
  \bibinfo{year}{1978}a.
\newblock \bibinfo{title}{A model and an algorithm for the dynamic traffic
  assignment problems}.
\newblock \bibinfo{journal}{Transportation science} \bibinfo{volume}{12},
  \bibinfo{pages}{183--199}.
\bibitem[{Merchant and Nemhauser(1978b)}]{merchant1978optimality}
\bibinfo{author}{Merchant, D.K.}, \bibinfo{author}{Nemhauser, G.L.},
  \bibinfo{year}{1978}b.
\newblock \bibinfo{title}{Optimality conditions for a dynamic traffic
  assignment model}.
\newblock \bibinfo{journal}{Transportation Science} \bibinfo{volume}{12},
  \bibinfo{pages}{200--207}.
\bibitem[{Nie and Zhang(2005a)}]{nie2005delay}
\bibinfo{author}{Nie, X.}, \bibinfo{author}{Zhang, H.}, \bibinfo{year}{2005}a.
\newblock \bibinfo{title}{Delay-function-based link models: their properties
  and computational issues}.
\newblock \bibinfo{journal}{Transportation Research Part B: Methodological}
  \bibinfo{volume}{39}, \bibinfo{pages}{729--751}.
\bibitem[{Nie and Zhang(2005b)}]{nie2005comparative}
\bibinfo{author}{Nie, X.}, \bibinfo{author}{Zhang, H.M.},
  \bibinfo{year}{2005}b.
\newblock \bibinfo{title}{A comparative study of some macroscopic link models
  used in dynamic traffic assignment}.
\newblock \bibinfo{journal}{Networks and Spatial Economics}
  \bibinfo{volume}{5}, \bibinfo{pages}{89--115}.
\bibitem[{Osorio et~al.(2011)Osorio, Fl{\"o}tter{\"o}d and
  Bierlaire}]{osorio2011dynamic}
\bibinfo{author}{Osorio, C.}, \bibinfo{author}{Fl{\"o}tter{\"o}d, G.},
  \bibinfo{author}{Bierlaire, M.}, \bibinfo{year}{2011}.
\newblock \bibinfo{title}{Dynamic network loading: a stochastic differentiable
  model that derives link state distributions}.
\newblock \bibinfo{journal}{Transportation Research Part B: Methodological}
  \bibinfo{volume}{45}, \bibinfo{pages}{1410--1423}.
\bibitem[{Peeta and Ziliaskopoulos(2001)}]{peeta2001foundations}
\bibinfo{author}{Peeta, S.}, \bibinfo{author}{Ziliaskopoulos, A.K.},
  \bibinfo{year}{2001}.
\newblock \bibinfo{title}{Foundations of dynamic traffic assignment: The past,
  the present and the future}.
\newblock \bibinfo{journal}{Networks and spatial economics}
  \bibinfo{volume}{1}, \bibinfo{pages}{233--265}.
\bibitem[{Ran and Boyce(1996)}]{ran1996link1}
\bibinfo{author}{Ran, B.}, \bibinfo{author}{Boyce, D.E.}, \bibinfo{year}{1996}.
\newblock \bibinfo{title}{A link-based variational inequality formulation of
  ideal dynamic user-optimal route choice problem}.
\newblock \bibinfo{journal}{Transportation Research Part C: Emerging
  Technologies} \bibinfo{volume}{4}, \bibinfo{pages}{1--12}.
\bibitem[{Ran et~al.(1996)Ran, Hall and Boyce}]{ran1996link2}
\bibinfo{author}{Ran, B.}, \bibinfo{author}{Hall, R.W.},
  \bibinfo{author}{Boyce, D.E.}, \bibinfo{year}{1996}.
\newblock \bibinfo{title}{A link-based variational inequality model for dynamic
  departure time/route choice}.
\newblock \bibinfo{journal}{Transportation Research Part B: Methodological}
  \bibinfo{volume}{30}, \bibinfo{pages}{31--46}.
\bibitem[{Richards(1956)}]{richards1956shock}
\bibinfo{author}{Richards, P.I.}, \bibinfo{year}{1956}.
\newblock \bibinfo{title}{Shock waves on the highway}.
\newblock \bibinfo{journal}{Operations research} \bibinfo{volume}{4},
  \bibinfo{pages}{42--51}.
\bibitem[{Rosenthal(2007)}]{rosen2007gams}
\bibinfo{author}{Rosenthal, R.E.}, \bibinfo{year}{2007}.
\newblock \bibinfo{title}{Gams: A user’s guide}, \bibinfo{publisher}{GAMS
  Development Corporation}.
\bibitem[{Szeto and Lo(2004)}]{szeto2004cell}
\bibinfo{author}{Szeto, W.}, \bibinfo{author}{Lo, H.K.}, \bibinfo{year}{2004}.
\newblock \bibinfo{title}{A cell-based simultaneous route and departure time
  choice model with elastic demand}.
\newblock \bibinfo{journal}{Transportation Research Part B: Methodological}
  \bibinfo{volume}{38}, \bibinfo{pages}{593--612}.
\bibitem[{Ukkusuri et~al.(2012)Ukkusuri, Han and Doan}]{ukkusuri2012dta}
\bibinfo{author}{Ukkusuri, S.V.}, \bibinfo{author}{Han, L.},
  \bibinfo{author}{Doan, K.}, \bibinfo{year}{2012}.
\newblock \bibinfo{title}{Dynamic user equilibrium with a path based cell
  transmission model for general traffic networks}.
\newblock \bibinfo{journal}{Transportation Research Part B: Methodological}
  \bibinfo{volume}{46}, \bibinfo{pages}{1657--1684}.
\bibitem[{Vickrey(1969)}]{vickrey1969congestion}
\bibinfo{author}{Vickrey, W.S.}, \bibinfo{year}{1969}.
\newblock \bibinfo{title}{Congestion theory and transport investment}.
\newblock \bibinfo{journal}{The American Economic Review} \bibinfo{volume}{59},
  \bibinfo{pages}{251--260}.
\bibitem[{Wie et~al.(1990)Wie, Friesz and Tobin}]{wie1990dynamic}
\bibinfo{author}{Wie, B.W.}, \bibinfo{author}{Friesz, T.L.},
  \bibinfo{author}{Tobin, R.L.}, \bibinfo{year}{1990}.
\newblock \bibinfo{title}{Dynamic user optimal traffic assignment on congested
  multidestination networks}.
\newblock \bibinfo{journal}{Transportation Research Part B: Methodological}
  \bibinfo{volume}{24}, \bibinfo{pages}{431--442}.
\bibitem[{Xu et~al.(1999)Xu, Wu, Florian, Marcotte and Zhu}]{xu1999advances}
\bibinfo{author}{Xu, Y.}, \bibinfo{author}{Wu, J.H.}, \bibinfo{author}{Florian,
  M.}, \bibinfo{author}{Marcotte, P.}, \bibinfo{author}{Zhu, D.},
  \bibinfo{year}{1999}.
\newblock \bibinfo{title}{Advances in the continuous dynamic network loading
  problem}.
\newblock \bibinfo{journal}{Transportation Science} \bibinfo{volume}{33},
  \bibinfo{pages}{341--353}.
\bibitem[{Yu et~al.(2020)Yu, Han and Ochieng}]{yu2020day}
\bibinfo{author}{Yu, Y.}, \bibinfo{author}{Han, K.}, \bibinfo{author}{Ochieng,
  W.}, \bibinfo{year}{2020}.
\newblock \bibinfo{title}{Day-to-day dynamic traffic assignment with imperfect
  information, bounded rationality and information sharing}.
\newblock \bibinfo{journal}{Transportation Research Part C: Emerging
  Technologies} \bibinfo{volume}{114}.
\bibitem[{Zhang et~al.(2013)Zhang, Nie and Qian}]{zhang2013modelling}
\bibinfo{author}{Zhang, H.}, \bibinfo{author}{Nie, Y.}, \bibinfo{author}{Qian,
  Z.}, \bibinfo{year}{2013}.
\newblock \bibinfo{title}{Modelling network flow with and without link
  interactions: the cases of point queue, spatial queue and cell transmission
  model}.
\newblock \bibinfo{journal}{Transportmetrica B: Transport Dynamics}
  \bibinfo{volume}{1}, \bibinfo{pages}{33--51}.
\bibitem[{Zhu and Ukkusuri(2015)}]{zhu2015dta}
\bibinfo{author}{Zhu, F.}, \bibinfo{author}{Ukkusuri, S.V.},
  \bibinfo{year}{2015}.
\newblock \bibinfo{title}{A linear programming formulation for autonomous
  intersection control within a dynamic traffic assignment and connected
  vehicle environment}.
\newblock \bibinfo{journal}{Transportation Research Part C: Emerging
  Technologies} \bibinfo{volume}{55}, \bibinfo{pages}{363--378}.

\end{thebibliography}







\end{document}